\documentclass [12pt]{article}

\newtheorem{theorem}{Theorem}[section]
\newtheorem{lemma}{Lemma}[section]

\newtheorem{remark}{Remark}[section]
\newcommand{\eqnsection}{
   \renewcommand{\theequation}{\thesection.\arabic{equation}}
   \makeatletter
   \csname @addtoreset\endcsname{equation}{section}
   \makeatother}

\def \ov{\overline}

\def \be{\begin{equation}}
\def \ee{\end{equation}}
\def \bt{\begin{theorem}}
\def \et{\end{theorem}}
\def \bea{\begin{eqnarray}}
\def \eea{\end{eqnarray}}
\def \bas{\begin{eqnarray*}}
\def \eas{\end{eqnarray*}}
\def \bl{\begin{lemma}} 
\def \el{\end{lemma}}

\def \al{\alpha}
\def \bb{\beta}
\def \ga{\gamma}

\def \de{\delta}
\def \De{\Delta}
\def \ep{\epsilon}

\def \la{\lambda}

\def \si{\sigma}

\def \th{\theta}

\def \ff{\infty}
\def \wh{\widehat}
\def \wt{\widetilde}

\def \rar{\rightarrow}

\def \cd{\,\cdot\,}

\def \II{{\cal I}}

\def \MM{{\cal M}}
\def \NN{{\cal N}}

\def \PP{{\cal P}}

\def \RR{{\cal R}}

\def \ZZ{{\cal Z}}

\def \({\left(}
\def \){\right)}

\def \lc{\left\{}
\def \rc{\right\}}

\def \nn{\nonumber}
\def \Proof{\noindent{\bf Proof $\,$ }}

\def \bc{\begin{center} }
\def \ec{\end{center} }
\def \bs{\begin{slide} }
\def \es{\end{slide} }

\def\square{{\vcenter{\vbox{\hrule height.3pt
        \hbox{\vrule width.3pt height5pt \kern5pt
           \vrule width.3pt}
        \hrule height.3pt}}}}
\def\qed{{\hfill $\square$ \bigskip}}

\eqnsection
\begin{document}

\def\wh{\widehat}
\def\ol{\overline}

\title{A CLT  for the $L^{2}$ moduli of continuity \\of local times of L\'evy processes }

\author{   Michael B. Marcus 
\hspace{ .2in}  Jay Rosen\thanks
     {The research of both authors was  supported, in part, by grants from the National Science
Foundation and PSC-CUNY.}}


\maketitle

\bibliographystyle{amsplain}

\begin{abstract}   
Let $X=\{X_{t},t\in R_{+}\}$   be a symmetric L\'{e}vy process 
with  local time   $\{L^{ x }_{ t}\,;\,(x,t)\in R^{ 1}\times  R^{  1}_{ +}\}$. 
When the L\'{e}vy exponent $\psi(\la)$  is regularly varying at infinity with   index $1<\beta\leq 2$ and satisfies some additional regularity conditions
\begin{eqnarray}
&& \sqrt{h\psi^{2}(1/h)} \lc \int ( L^{ x+h}_{1}- L^{ x}_{ 1})^{ 2}\,dx- E\( \int (  L^{ x+h}_{1}- L^{ x}_{ 1})^{ 2}\,dx\)\rc\nn\\
&&\hspace{1 in}
\stackrel{\mathcal{L}}{\Longrightarrow}  (8c_{\beta,1})^{1/2}\,\,\eta\,\, \( \int  (L_{1}^{x})^{2}\,dx\)^{1/2}\,\nn, 
\end{eqnarray}
as $h\rar 0$, where
$\eta$ is a normal random variable with mean zero and variance one that is independent of $L^{ x }_{ t}$, and $c_{\beta,1}$ is a known constant.

\end{abstract}

 \footnotetext{  { Key words and phrases:} Central Limit Theorem,   $L^{2}$ moduli of continuity,   local time,  L\'evy process.}

 \footnotetext{ {  AMS 2000 subject classification:}  Primary 60F05, 60J55, 60G51.}

\section{Introduction  }

 In   \cite{CLR}, with X.  Chen and W. Li, we obtain  a central limit theorem for the $L^{2}$ modulus of continuity of local times of Brownian motion. In this paper, this result is extended to symmetric stable processes with index   $1<\bb< 2$ and, in one respect, to a much larger class of L\'evy processes. 
 
 Let $X=\{X_{t},t\in R_{+}\}$ \label{page1} be a symmetric stable process 
of index $1<\bb<2$. We normalize $X $ so that
\begin{equation}
E\(e^{i\la X_{t}}\)=e^{-|\la|^{\bb}t}.\label{st.1}
\end{equation}
Set
\begin{equation}
c_{\bb,0}={2\over \pi}\int_{0}^{\ff} {\sin^{2}p/2  \over p^{\bb}}\,dp\label{st.22}
\end{equation}
and 
 \begin{equation}
c_{\bb,1}={16\over \pi}\int_{0}^{\ff} {\sin^{4}p/2  \over p ^{2\bb}}\,dp.\label{st.22m}
\end{equation}

 Let $\{L^{ x }_{ t}\,;\,(x,t)\in R^{ 1}\times  R^{  1}_{ +}\}$ denote the local time  of   $X$ and define
 \begin{equation}
\al_{t}=\int  (L_{t}^{x})^{2}\,dx.\label{1.4}
   \end{equation}
 
 An integral sign without limits is to be read as $\int_{-\ff}^{\ff}\,$.

\begin{theorem}\label{theo-clt2s} Let $\{L^{ x }_{ t}\,;\,(x,t)\in R^{ 1}\times  R^{  1}_{ +}\}$ be the local time  of the symmetric stable process $X $
of index   $1<\bb\le 2$. For each fixed $t$  
\begin{equation} { \int ( L^{ x+h}_{t}- L^{ x}_{ t})^{ 2}\,dx-4 c_{\bb,0}h^{\bb-1}t\over h^{(2\bb-1)/2}}
\stackrel{\mathcal{L}}{\Longrightarrow} (8c_{\bb,1})^{1/2}\sqrt{\al_{t}}\,\eta \label{5.0weaks}
\end{equation}  
as $h\rar 0$, where $\eta$ is a normal random variable with mean zero and variance one that is independent of $\al_{\cdot}$.  Equivalently 
\begin{equation} { \int ( L^{ x+1}_{t}- L^{ x}_{ t})^{ 2}\,dx- 4c_{\bb,0} t\over t^{ (2\bb-1)/(2\bb)}}
\stackrel{\mathcal{L}}{\Longrightarrow}(8c_{\bb,1})^{1/2}\sqrt{\al_{1}}\,\eta\label{5.0tweaks}
\end{equation}
as $t\rar\ff$.
\end{theorem} 

  In     \cite{lp} we show that under the hypotheses of Theorem \ref{theo-clt2s},    for all
$t\in R_+$,  
\begin{equation}
\lim_{ h\downarrow 0}  \int  {  (L^{ x+h}_{ t} -L^{ x }_{
t})^{ 2}\over h^{\bb-1}}\,dx =4c_{\bb,0}\,t \hspace{.2 in}\mbox{ a. s.}, \label{srp3.1}
\end{equation}
 and in $L^{1}$,  so (\ref{5.0weaks}) does have the form of a central limit theorem. 
 
  Theorem \ref{theo-clt2s}, with $\la=2$, is actually a   central limit theorem for   the local time of    $B(2t)$ where $B(t)$ is Brownian motion. This is because the characteristic function of a canonical 2--stable process  is not the same as the characteristic function of Brownian motion, as one can see in (\ref{st.1}). The analogue of (\ref{5.0weaks}) for Brownian motion is given in \cite[Theorem 1.1]{CLR}. (  Naturally, the constants are different.)   The proofs in \cite{CLR} have a superficial resemblance to the proofs in this paper, but use many special properties of Brownian motion,  particularly the scaling property. While this type of proof could  be extended to cover symmetric stable processes, we do not see how to extend them to the much larger class of L\'evy processes   considered in this paper, for which  there is no  scaling property.
 
The equivalence of (\ref{5.0weaks}) and (\ref{5.0tweaks}) 
 follows from the scaling relationship for stable processes, see e.g. \cite[Lemma 10.5.2]{book},
\begin{equation}
\{ L^{ x}_{  t/h^{\bb}}\,;\,( x,t)\in R^{ 1}\times R^{ 1}_{ +}\}
\stackrel{\mathcal{L}}{=}\{h^{ -(\bb-1)} L^{h x}_{ t}\,;\,( x,t)\in R^{ 1}\times
R^{ 1}_{ +}\},\label{scales}
\end{equation}
which implies that 
\begin{equation}
\int ( L^{ x+h}_{t}- L^{ x}_{ t})^{ 2}\,dx\stackrel{\mathcal{L}}{=}h^{2\bb-1} \int ( L^{ x+1}_{t/h^{\bb}}- L^{ x}_{ t/h^{\bb}})^{ 2}\,dx\label{scls}
\end{equation}
   and
\begin{equation}
 \al_{t}  \stackrel{\mathcal{L}}{=}t^{(2\bb-1)/\bb}\al_{1}.\label{1.10}
   \end{equation}
 
   In fact, because of (\ref{scales}), to prove Theorem \ref{theo-clt2s}, it suffices to prove (\ref{5.0weaks})   with $t=1$.   In this paper we do this for a much larger class of symmetric L\'evy processes than symmetric   stable processes.
 We extend the  definition on page \pageref{page1} so that now $X=\{X_{t},t\in R_{+}\}$   is a symmetric L\'{e}vy process 
with
\begin{equation}
E\(e^{i\la X_{t}}\)=e^{-\psi(\la)t}\label{2rt.1v}
\end{equation}  
where
 $\psi(\la)$   satisfies
 \begin{itemize}
 \item[1.]
 \be
 \mbox{ $ \psi(\la)$ is regularly varying at infinity with   index $1<\bb\leq 2$} ; \label{regcond}
 \ee
 \item[2.] $\psi$ is twice differentiable almost surely and there exist  constants $D_{1},D_{2}<\ff$ such that for all  $\la\ge 1$ 
   \be
   \la |\psi'(\la)|\le D_{1}\psi(\la)\quad\mbox{and}\quad   \la^{2} |\psi''(\la)|\le D_{2}\psi(\la)\hspace{.2 in}a.s.\label{88.m}
\ee   and  
   \begin{equation}
  \int_{0}^{1}  (\psi'(\la))^{2}\,d\la<\ff,  \quad \int_{0}^{1}  |\psi''(\la)|\,d\la<\ff.\label{114}
   \end{equation}
   \item[3.]  
   \begin{equation}
   \int_{0}^{1}\frac{\psi(\la)}{\la}\,d\la<\ff\label{1.16}
   \end{equation}
   \end{itemize}
    (This last condition is very weak.)

\medskip	  The next theorem is the main result in this paper

 \begin{theorem}\label{theo-clt2r} Let $\{L^{ x }_{ t}\,;\,(x,t)\in R^{ 1}\times  R^{  1}_{ +}\}$ be the local time  of the symmetric  L\'{e}vy process $X$
whose L\'{e}vy exponent
 $\psi(\la)$  satisfies (\ref{regcond})--(\ref{1.16}). Then
\bea
&& \sqrt{h\psi^{2}(1/h)} \lc \int ( L^{ x+h}_{1}- L^{ x}_{ 1})^{ 2}\,dx- E\( \int (  L^{ x+h}_{1}- L^{ x}_{ 1})^{ 2}\,dx\)\rc\nn\\
&&\hspace{2 in}
\stackrel{\mathcal{L}}{\Longrightarrow}  (8c_{\bb,1})^{1/2}\sqrt{\al_{  1}}\,\eta\label{r5.0weaks},
\eea
as $h\rar 0$. (Here $\al_{1}$ is as defined in (\ref{1.4}) but for the local time of   $X$.)
\end{theorem} 

 \medskip	  Theorem  \ref{theo-clt2r}  tells us something about Theorem \ref{theo-clt2s} that isn't obvious from the statement of  Theorem \ref{theo-clt2s}. For a symmetric stable process of index $\bb$, $\psi(p)=|p|^{\bb}$ for all $p$. However
we see from (\ref{88.m}) and (\ref{114}) that  (\ref{5.0weaks}) in Theorem \ref{theo-clt2s}  depends primarily on the fact that 
\begin{equation}
   \lim_{|p|\to \ff}\frac{\psi(p)}{|p|^{\bb}}=C 
   \end{equation}
 for some finite, nonzero constant $C$,  i.e. the behavior of $\psi(\cd)$ at zero is not relevant. The conditions in (\ref{114}) are very weak.

 It is interesting to consider the following rearrangements of (\ref{5.0weaks}) and (\ref{5.0tweaks}): 
 \begin{equation} { \int \( {L^{ x+h}_{t}- L^{ x}_{ t}\over h^{ ( \bb-1)/2 }}\)^{ 2}\,dx-4 c_{\bb,0} t\over h^{1/2}}
\stackrel{\mathcal{L}}{\Longrightarrow} (8c_{\bb,1})^{1/2}\sqrt{\al_{t}}\,\eta \label{5.0weaksdd},
\end{equation}  
as $h\rar 0$, and
\begin{equation} t^{1/(2\bb)}\( \frac{1}{t}\int ( L^{ x+1}_{t}- L^{ x}_{ t})^{ 2}\,dx- 4c_{\bb,0} \) 
\stackrel{\mathcal{L}}{\Longrightarrow}(8c_{\bb,1})^{1/2}\sqrt{\al_{1}}\,\eta\label{5.0tweaksdd},
\end{equation}
as $t\rar\ff$. Written this way, (\ref{5.0weaksdd}), with $h^{1/2}$ in the denominator, looks more like a classical CLT    and  (\ref{5.0tweaksdd}) displays the behavior of the long time average of the $L^{2}$ norm of the increments of the local time.

 The critical ingredient in the proof  of  Theorems \ref{theo-clt2s}   is Lemma \ref{lem-multiple} which gives moments for the  $L^{2}$ norm of increments of local times of L\'evy processes satisfying (\ref{regcond})--(\ref{1.16}). 
  In Section \ref{sec-indf} we state  Lemmas \ref{lem-vproprvt}--\ref{lem-4.3} which give properties of $p_{t}(x)$, $\De^{h}p_{t}(x)$ and $\De^{h}\De^{h}p_{t}(x)$, where $p_{t}(x)$ is the transition probability density of   these L\'evy processes. The lemmas are used in Section \ref{sec-3} in the proof of Lemma \ref{lem-multiple}. They are proved in Section \ref{sec-est}. Theorem \ref{theo-clt2r} is proved in Section \ref{sec-CLT}. In its proof we need good asymptotic estimates of   the mean and variance of
\begin{equation}
    \int ( L^{ x+h}_{1}- L^{ x}_{1})^{2}\,dx.\label{1.23}
   \end{equation} 
These are stated in   Section \ref{sec-CLT} and proved in Section \ref{sec-vare}.   Theorem  \ref{theo-clt2s} is just a special case of  Theorem  \ref{theo-clt2r} with obvious extensions which are possible for stable processes because of the scaling property of local times of stable processes. The few remarks that show how Theorem  \ref{theo-clt2s} follows from   Theorem  \ref{theo-clt2r} are given at the end of Section \ref{sec-vare}.  In Section \ref{sec-Kac}, which is an appendix, we prove the version of Kac Moment Formula that we use in this paper.

    \section{Estimates for the probability densities of certain L\'evy processes  }\label{sec-indf}

 Let $p_{s}(x)$ denote the density of the symmetric  L\'{e}vy process $X=\{X_{s}, s\in R_{+}\}$ with L\'{e}vy exponent $\psi (p)$ as described in  (\ref{2rt.1v}).  Let $\De_{ x}^{ h}$  denote
the finite difference operator on the variable $x$, i.e.
\begin{equation}
\De_{x}^{ h}\,f(x)=f(x+h)-f(x).\label{pot.3w}
\end{equation}
We write $\De^{ h}$ for $\De_{x}^{ h}$ when the variable $x$ is clear.

\medskip	
The following lemmas   provide the main estimates we use in this paper. Their  proofs are given in Section \ref{sec-est}. 

\begin{lemma}\label{lem-vproprvt} Let    $X $ be a symmetric L\'{e}vy process 
with  L\'{e}vy exponent
 $\psi(\la)$  that is regularly varying at infinity with index $1<\bb\leq 2$ and satisfies (\ref{88.m}) and   (\ref{114}).  Let $p_{s }(x)$ denote the transition probability density of $X$.
   Then   
   \begin{equation}
p_{s }(x)\leq C{\psi^{-1}(1/s)\vee 1\over 1+x^{2}} ,\qquad\forall\,x\in R^{1}, s\in (0,1];\label{2.dens}
   \end{equation}
   \begin{equation}
u(x):= \int_{0}^{1}\, \,p_{s }(x)\,ds\leq {C\over 1+x^{2}} ,\qquad\forall\,x\in R^{1};\label{2.2}
 \ee
 \be 
\int  \int_{0}^{t}\, \,p_{s }(x)\,ds\,dx=t,\label{bl.1}
   \end{equation}
  and for all $h$ sufficiently small  
      \bea
v(x):=  \int_{0}^{1}\, |\De ^{ h}\,p_{s }(x)|\,ds&\leq& C\({1 \over h\psi (1/h)}\wedge \frac{h}{|x|}\wedge\frac{h }{x^{2}}\)   \label{bl.2}\\
&\le &{C \over h\psi (1/h)}\({1\over 1+x^{2}}\)\nn,
   \eea
  \begin{equation} 
\int  v (x)\,dx= O\(h\log 1/h\),   \label{bl.01}
\end{equation}  
and  
 \begin{equation}
 \int  v ^{p}(x)\,dx= O\(\frac{h}{h^{p-1}\psi^{p-1} (1/h)}\),\quad p\ge 2,   \label{lastf}  
   \end{equation} 
   as $h\to 0$.  
    In addition,  
 \begin{equation} 
w(x):= \int_{0}^{1}\,|\De^{ h}\De^{ -h} \,p_{s }(x)|\,ds \leq C\({1 \over h\psi (1/h)}\wedge { 1\over \psi (1/h)|x| }\wedge {h^{2}\over |x|^{2}}\);\label{bl.3}
\end{equation}
\begin{equation} 
  \int w(x)\,dx= O\({\log (1/h) \over \psi (1/h)} \);\label{bl.4}
\end{equation} 
\begin{equation} 
  \int  w^{2}(x)\,dx= O\({1 \over h\psi^{2} (1/h)}\);\label{bl.5}
\end{equation}
\begin{equation} 
  \int_{|x|\geq u} w^{2}(x)\,dx\leq O\({1 \over u\psi^{2} (1/h)}\),\label{bl.5b}
\end{equation}
as $h\to 0.$
\end{lemma}

\begin{lemma}\label{lem-h3} Under the same hypotheses as Lemma \ref{lem-vproprvt},   set
\be c_{\psi,h,1} = \int \(\int \(\De^{ h}\De^{ -h}\,p_{s}(x)\) \,ds \)^{2}\,   \,dx. \label{jay.31}
\ee 
Then 
\begin{equation} 
\lim_{h\rar 0}h\psi^{2}(1/h)c_{\psi,h,1}=   c_{\bb, 1}  \label{lim.1}
\end{equation}
  and 
  \be 
h\psi^{2}(1/h) \(  c_{\psi,h,1} - \int  \(\int_{0}^{\sqrt{h}}\De^{ h}\De^{ -h}\,\,p_{s }(x)\,ds\)^{2}\,dx \)   =  O( {h^{1/2 }}).\label{lim.1r}
\ee 
\end{lemma}
 
\begin{remark} \label{rem-1}{\rm We allow $\psi$ to be regularly varying at infinity with index 2, but note that  that because $\psi$ is the L\'evy exponent of a symmetric L\'evy process
\begin{equation}
   \psi(\la)=O (\la^{2})  \qquad\mbox{as}\,\,\la\to\ff.\label{2.12}
   \end{equation}
(See, e.g.,   \cite[Lemma 4.2.2]{book}   and then include Brownian motion.)
 }\end{remark}

\begin{lemma}\label{lem-4.2} Under the same hypotheses as Lemma \ref{lem-vproprvt},
\be  
  \sup_{\de\le r\le 1}p_{r }(0)  \le C\(\psi^{-1}(1/\de)\vee 1\)\label{h.101xx.1};
 \ee 
 \bea   
 \sup_{\de\le r\le 1}|\De^{h}p_{r }(0)|& \le & {C\over {\de^{3}}} h^{2};\label{4.4}
 \eea
 and
  \bea  
 \sup_{\de\le r\le 1}|\De^{h}\De^{-h}p_{r }(0)| 
 &\le  &  {C\over {\de^{3}}} h^{2}.\label{4.5}
 \eea 
 \el

 \bl  \label{lem-4.3} Let $
0<\de<1$, then, under the hypotheses of Theorem \ref{theo-clt2r}, for \begin{equation}
\ov   u_{\de}(x):=  \sup_{\de\le r\le 1}p_{r }(x), \,\, \qquad \ov    v_{\de}(x):=  \sup_{\de\le r\le 1}|\De^{h}p_{r }(x)|,   
   \end{equation}
   \[ 
 \mbox{and }\qquad\ov    w_{\de}(x):=  \sup_{\de\le r\le 1}|\De^{h}\De^{-h}p_{r }(x)|,
   \]
we have    
      \begin{equation} 
\ov u_\de(x) \leq C\psi^{-1}(1/\de)    \(1\wedge {1\over x^{2}} \) \label{b.l.2wa},
   \end{equation}      \begin{equation} 
\ov v_\de(x) \leq { C\over\de^{3}} h  \(1\wedge {1\over x^{2}} \) \label{b.l.2w},
   \end{equation}
   and
    \begin{equation} 
\ov   w_\de(x)  \leq { C\over\de^{3}} h^{2} \(1\wedge {1\over x^{2}} \).\label{bl.3q}
\end{equation}

In addition 
 \be 
\int \ov   u_{\de}(x)\,dx\leq  C\psi^{-1}(1/\de)  ,\qquad\int (\ov   u_{\de}(x))^{2}\,dx\leq C\(\psi^{-1}(1/\de) \)^{2}, \label{bl.1x}
   \end{equation}
   \begin{equation} 
\int  \ov v_\de  (x)\,dx\le {C\over  { \de^{3} } }h, \qquad \label{last.f} 
\int  \ov v_\de^{2}(x)\,dx\le  {C\over  { \de^{6} } }h^{2} ,\label{1.32}
\end{equation}
\begin{equation} 
\int  \ov w_\de (x)\,dx\leq  {C\over  { \de^{3} } }h^{2} ,\qquad\int  \ov w_\de^{2}(x)\,dx\leq  {C\over  { \de^{6} } }h^{4} ,\label{bl.5ba}
\end{equation}
as $h\to 0$. 
 \end{lemma}

\section{Moments of  increments of local times.}\label{sec-3}

Let
\begin{equation}
I_{j,k,h}:=\int ( L^{ x+h}_{1}- L^{ x}_{1})\circ\th_{j}\,\, (   L^{ x+h}_{1}-  L^{ x}_{ 1})\circ  \th_{k} \,\, \,dx\label{m.1}
\end{equation}
and
\begin{equation}
  \al_{j,k  } := \int  L^{ x}_{1 }\circ \th_{j}\,\,  L^{ x}_{1}\circ  \th_{k}\,\,  \,dx.\label{m.1a}
\end{equation}

The goal of this section is to establish the following lemma.      
\bl\label{lem-multiple}
 Let 
$m_{j,k}$, $0\leq j<k \leq K$ be positive integers with  $\sum_{ j,k=0, j<k }^{K}\newline  m_{j,k}=m$. 
 If all the   integers $m_{j,k}$ are even, then  for some $\ep>0$
\begin{eqnarray} && 
E\(\prod_{\stackrel{j,k=0}{j< k}}^{K} \(I_{j,k,h} \)^{m_{j,k}}\)
\label{m.2}\\ 
&&\qquad
=  \prod_{\stackrel{j,k=0}{j< k}}^{K}{( 2n_{j,k})!\over 2^{ n_{j,k}}(n_{j,k}!)}\( 4c_{\psi,h,1}\)^{ n_{j,k}} E\(\prod_{\stackrel{j,k=0}{j< k}}^{K}\(\al_{j,k}\)^{
n_{j,k}}\) +O\(h^{(2\bb-1)n+\ep}\) ,
\nn
\end{eqnarray}
where  $n_{j,k}=m_{j,k}/2$  and $n=m/2$.
If any of the $m_{j,k}$
are odd, then  
\begin{equation}
E\(\prod_{\stackrel{j,k=0}{j< k}}^{K} \(I_{j,k,h} \)^{m_{j,k}}\)=   O\(h^{(2\bb-1)m/2+\ep}\).\label{m.2a}
   \end{equation} 
   In (\ref{m.2}) and (\ref{m.2a}) the error terms  may depend on $m$, but not on the 
individual terms $m_{j,k} $. 
\el
 
 \Proof 
We can write
\begin{eqnarray} \lefteqn{
E\(\prod_{\stackrel{j,k=0}{j< k}}^{K} \(I_{j,k,h} \)^{m_{j,k}}\)\label{m.3}}\\
&&=E\(\prod_{\stackrel{j,k=0}{j< k}}^{ } \prod_{i=1}^{m_{j,k}}\(\int ( \De^{h}_{x_{j,k,i}} L^{ x_{j,k,i}}_{1}\circ\th_{j})\,\, (  \De^{h}_{x_{j,k,i}}L^{ x_{j,k,i}}_{ 1}\circ  \th_{k}) \,\, \,dx_{j,k,i} \) \)\nn\\
&&=\int  \lc \prod_{\stackrel{j,k=0}{j< k}}^{K} \prod_{i=1}^{m_{j,k}}\De^{h,j}_{x_{j,k,i}}  \De^{h,k}_{x_{j,k,i}} \rc E\(\prod_{\stackrel{j,k=0}{j< k}}^{K} \prod_{i=1}^{m_{j,k}}\(  (  L^{ x_{j,k,i}}_{1}\circ\th_{j})\,\, (  L^{ x_{j,k,i}}_{ 1}\circ  \th_{k})  \) \)\nonumber\\
&&\hspace{3.5 in}\prod_{\stackrel{j,k=0}{j< k}}^{K} \prod_{i=1}^{m_{j,k}}\,dx_{j,k,i}\nn,
\end{eqnarray}
where the notation $\De^{h,j}_{x_{j,k,i}}$   indicates that we apply the difference operator $\De^{h}_{x_{j,k,i}}$ to    $ L^{ x_{j,k,i}}_{1}\circ\th_{j}$.  Note that there are $2m$ applications of the difference operator $\De^{h}$.  

Consider 
\begin{equation}
   E\(\prod_{\stackrel{j,k=0}{j< k}}^{K} \prod_{i=1}^{m_{j,k}}\(  (  L^{ x_{j,k,i}}_{1}\circ\th_{j})\,\, (  L^{ x_{j,k,i}}_{ 1}\circ  \th_{k})  \) \).
   \end{equation}
 We collect  all the factors containing $\th_{l}$ and write
\begin{eqnarray}
&&E\(\prod_{\stackrel{j,k=0}{j< k}}^{K} \prod_{i=1}^{m_{j,k}}\(  (  L^{ x_{j,k,i}}_{1}\circ\th_{j})\,\, (  L^{ x_{j,k,i}}_{ 1}\circ  \th_{k})  \) \)
\label{reo.1}\\
&&\qquad=E\(\prod_{l=0}^{K}\lc \(\prod_{j=0}^{l-1} \prod_{i=1}^{m_{j,l}}   L^{ x_{j,l,i}}_{1}\) 
 \(\prod_{k=l+1}^{K} \prod_{i=1}^{m_{l,k}}L^{ x_{l,k,i}}_{1}\)   \rc \circ  \th_{l}\)   \nonumber\\
&&\qquad=E\(\prod_{l=0}^{K}H_{l} \circ  \th_{l}\)   \nonumber,
\end{eqnarray}
where
\begin{equation}
H_{l}=\(\prod_{j=0}^{l-1} \prod_{i=1}^{m_{j,l}}    L^{ x_{j,l,i}}_{1}\) 
 \(\prod_{k=l+1}^{K} \prod_{i=1}^{m_{l,k}}    L^{ x_{l,k,i}}_{1}\).\label{reo.2}
\end{equation}
By the Markov property
\begin{equation}
E\(\prod_{l=0}^{K}H_{l} \circ  \th_{l}\)=E\(H_{0} \,E^{X_{1}}\(\prod_{l=1}^{K}H_{l} \circ  \th_{l-1}\)\)  .\label{reo.3}
\end{equation}
Let 
\be
m_{l}=\sum_{k=l+1}^{K}\,m_{l,k}+\sum_{j=0}^{l-1}\,m_{j,l},\qquad l=0,\ldots,K-1,\label{m4.9}
\ee
and note that  $m_{l}$ is the number of local time factors in $H_{l}$.  

Let 
\begin{equation}
   f(y)=E^{y}\(\prod_{l=1}^{K}H_{l} \circ  \th_{l-1}\).
   \end{equation}
   It follows from  Kac's Moment Formula, Theorem \ref{KMT},   for any $z\in R^{1}$
\begin{eqnarray}
\lefteqn{E^{z}\(\prod_{l=0}^{K}H_{l} \circ  \th_{l}\)}\\
&=& E^{z}\(H_{0} \,f(X_{1})\)\nn
\label{reo.4} \\
&=&   \sum_{\pi_{0}}\int_{\{\sum_{q=1}^{m_{0}}r_{0,q}\leq 1\}}
p_{r_{0,1}}(x_{\pi_{0}(1)}-z) \prod_{q=2}^{m_{0}}p_{r_{0,q}}(x_{\pi_{0}(q)}-x_{\pi_{0}(q-1)}) \nonumber\\
&&\hspace{1 in}\(\int p_{(1-\sum_{q=1}^{m_{0}}r_{0,q}) }(y-x_{\pi_{0}(m_{0})})f(y)\,dy\)\prod_{q=1}^{m_{0}}\,dr_{0,q} \nonumber,
\end{eqnarray}
where the sum runs over all 
  bijections $\pi_{0} $ from $[1,m_{0}]$ to
\begin{equation}
I_{0}=\bigcup_{k=1}^{K}\{(0,k,i),\,1\leq i\leq m_{0,k}\}.\label{reo.5}
\end{equation}
Clearly, $I_{0}$ is the set of subscripts  of  the terms $x_{\cd} $ appearing in the local time factors in $H_{0}$. 

By the Markov property 
\bea
   f(y)&=&E^{y}\(H_{1}E^{X_{2}}\(\prod_{l=2}^{K}H_{l} \circ  \th_{l-2}\)\)\label{4.13}\\
   &:=&E^{y}\(H_{1}g(X_{2}) \).\nn
   \eea
 Therefore, by (\ref{reo.3})--(\ref{4.13}), for any $z'\in R^{1}$ 
 \bea
\lefteqn{E^{z'}\(\prod_{l=0}^{K}H_{l} \circ  \th_{l}\)\label{reo.6}}\\
&&=E^{z'}\(H_{0}E^{X_{1}}\(H_{1} \,g(X_{2})\)\)
\nn\\
&&=  \sum_{\pi_{0}}\int_{\{\sum_{q=1}^{m_{0}}r_{0,q}\leq 1\}}
p_{r_{0,1}}(x_{\pi_{0}(1)}-z') \prod_{q=2}^{m_{0}}p_{r_{0,q}}(x_{\pi_{0}(q)}-x_{\pi_{0}(q-1)}) \nonumber\\
&&\hspace{.5 in}\(\int p_{(1-\sum_{q=1}^{m_{0}}r_{0,q}) }(y-x_{\pi_{0}(m_{0})})
E^{y}\(H_{1} \,g(X_{2})\)\,dy\)\prod_{q=1}^{m_{0}}\,dr_{0,q} \nonumber\\
&&=  \sum_{\pi_{0}}\int_{\{\sum_{q=1}^{m_{0}}r_{0,q}\leq 1\}}
p_{r_{0,1}}(x_{\pi_{0}(1)}-z') \prod_{q=2}^{m_{0}}p_{r_{0,q}}(x_{\pi_{0}(q)}-x_{\pi_{0}(q-1)}) \nonumber\\
&&\hspace{.5 in}  p_{(1-\sum_{q=1}^{m_{0}}r_{0,q}) }(y-x_{\pi_{0}(m_{0})})\nn\\
&&\hspace{.3 in} \sum_{\pi_{1}}\int_{\{\sum_{q=1}^{m_{1}}r_{1,q}\leq 1\}}
p_{r_{1,1}}(x_{\pi_{1}(1)}-y) \prod_{q=2}^{m_{1}}p_{r_{1,q}}(x_{\pi_{1}(q)}-x_{\pi_{1}(q-1)}) \nonumber\\
&&\hspace{.3 in}\(\int p_{(1-\sum_{q=1}^{m_{1}}r_{1,q}) }(y'-x_{\pi_{1}(m_{1})})g(y')\,dy'\)\prod_{q=1}^{m_{1}}\,dr_{1,q} \nonumber
\,dy\prod_{q=1}^{m_{0}}\,dr_{0,q} \nonumber
\end{eqnarray}
where the second sum   runs over all 
  bijections $\pi_{1} $ from $[1,m_{1}]$ to
\begin{equation}
I_{1}=\{(0,1,i),\,1\leq i\leq m_{0,1}\}\bigcup_{k=2}^{K}\{(1,k,i),\,1\leq i\leq m_{1,k}\}\label{reo.7}
\end{equation}

As above, $I_{1}$ is the set of subscripts  of  the terms $x_{\cd} $ appearing in the   local time factors in $H_{1}$. 

  We now use the Chapman-Kolmogorov equation to integrate with respect to $y$ to get
\begin{eqnarray}
\lefteqn{E^{z'}\(H_{0}E^{X_{1}}\(H_{1} \,g(X_{1})\)\)
\label{reo.7a}}\\
&&=  \sum_{\pi_{0},\pi_{1}}\int_{\{\sum_{q=1}^{m_{0}}r_{0,q}\leq 1\}}
p_{r_{0,1}}(x_{\pi_{0}(1)}-z') \prod_{q=2}^{m_{0}}p_{r_{0,q}}(x_{\pi_{0}(q)}-x_{\pi_{0}(q-1)}) \nonumber\\ 
&&\hspace{.3 in}  \int_{\{\sum_{q=1}^{m_{1}}r_{1,q}\leq 1\}}
p_{(1-\sum_{q=1}^{m_{0}}r_{0,q})+r_{1,1}}(x_{\pi_{1}(1)}-x_{\pi_{0}(m_{0})}) \nonumber\\
&&\hspace{2 in}\prod_{q=2}^{m_{1}}p_{r_{1,q}}(x_{\pi_{1}(q)}-x_{\pi_{1}(q-1)}) \nonumber\\
&&\hspace{.3 in}\(\int p_{(1-\sum_{q=1}^{m_{1}}r_{1,q}) }(y'-x_{\pi_{1}(m_{1})})g(y')\,dy'\)\prod_{q=1}^{m_{1}}\,dr_{1,q} \nonumber
 \prod_{q=1}^{m_{0}}\,dr_{0,q}. \nonumber
\end{eqnarray}

Iterating this procedure, and recalling (\ref{reo.1}) we see that
\begin{eqnarray}
\lefteqn{E\(\prod_{\stackrel{j,k=0}{j< k}}^{K} \prod_{i=1}^{m_{j,k}}\( (  L^{ x_{j,k,i}}_{1}\circ\th_{j})\,\, ( L^{ x_{j,k,i}}_{ 1}\circ  \th_{k})   \) \)
\label{m.4}}\\
&&=  \sum_{\pi_{0},\ldots, \pi_{K}}  \prod_{l=0}^{K}\int_{\{\sum_{q=1}^{m_{l}}r_{l,q}\leq 1\}}
p_{(1-\sum_{q=1}^{m_{l-1}}r_{l-1,q})+r_{l,1}}(x_{\pi_{l}(1)}-x_{\pi_{l-1}(m_{l-1})})  \nonumber\\
&&\hspace{1.8 in}\prod_{q=2}^{m_{l}}p_{r_{l,q}}(x_{\pi_{l}(q)}-x_{\pi_{l}(q-1)})\prod_{q=1}^{m_{l}}\,dr_{l,q},\nonumber
\end{eqnarray}
  where $\pi_{-1}(m_{-1}):=0$ and   $1-\sum_{q=1}^{m_{ -1}}r_{ -1,q}:=0$. In (\ref{m.4})
  the sum runs over all 
$\pi_{0},\ldots, \pi_{K}$ such that each $\pi_{l}$ is a bijection from $[1,m_{l}]$ to  
\begin{equation}
I_{l}=\bigcup_{j=0}^{l-1}\{(j,l,i),\,1\leq i\leq m_{j,l}\}\bigcup_{k=l+1}^{K}\{(l,k,i),\,1\leq i\leq m_{l,k}\}.\label{m.5}
\end{equation}
As in the observations about $I_{0}$ and $I_{1}$, we see  that  $I_{l}$ is the set of subscripts of  the terms $x_{\cd} $ terms appearing in the   local time factors in   $H_{l}$. Since there are   $2m$ local time factors  we have that   $\sum_{l=0}^{K}m_{l}=2m$.

\medskip	We  now use (\ref{m.4}) in (\ref{m.3}) and continue to develop an expression for the left-hand side of   (\ref{m.3}). Let $\mathcal{B}$ to denote the set of $K+1$ tuples $\pi=(\pi_{0},\ldots, \pi_{K})$ of bijections described in (\ref{m.5}). Clearly 
\be
|\mathcal{B}|=\prod_{l=0}^{K}m_{l}!\leq (2m)!.
\ee
  Also, similarly to the way we obtain the first equality in (\ref{reo.1}), we see that  
\begin{equation}
\prod_{\stackrel{j,k=0}{j< k}}^{K} \prod_{i=1}^{m_{j,k}}\De^{h,j}_{x_{j,k,i}}  \De^{h,k}_{x_{j,k,i}} = \prod_{l=0}^{K}\prod_{q=1}^{m_{l}}\De^{h,l}_{x_{\pi_{l}(q)}}\label{m.5a}.
\end{equation}
Consequently  
\be  
E\(\prod_{\stackrel{j,k=0}{j< k}}^{K} \(I_{j,k,h} \)^{m_{j,k}}\)
 \label{m.6a} = \sum_{\pi_{0},\ldots, \pi_{K}}\int \wt\mathcal{T}_{h}( x;\,\pi  )\prod_{j,k,i}\,dx_{j,k,i} \label{4.22}
\ee  
where  we take the product over  $ \{  0\leq j<k\leq K,\, 1\leq i\leq m_{j,k}\}$,   $\pi\in \mathcal{B}$ 
and
\begin{eqnarray}
\lefteqn{ \wt\mathcal{T}_{h}( x;\,\pi  )
\label{m.4a}}\\
&&=    \prod_{l=0}^{K}\prod_{q=1}^{m_{l}}\De^{h}_{x_{\pi_{l}(q)}}\int_{\{\sum_{q=1}^{m_{l}}r_{l,q}\leq 1\}}
p_{(1-\sum_{q=1}^{m_{l-1}}r_{l-1,q})+r_{l,1}}(x_{\pi_{l}(1)}-x_{\pi_{l-1}(m_{l-1})})  \nonumber\\
&&\hspace{1.8 in}\prod_{q=2}^{m_{l}}p_{r_{l,q}}(x_{\pi_{l}(q)}-x_{\pi_{l}(q-1)})\prod_{q=1}^{m_{l}}\,dr_{l,q}.\nonumber
\end{eqnarray}
\medskip	We continue to rewrite the right-hand side of (\ref{m.6a}). 

  In (\ref{m.4a}),   each   difference operators, say $\De^{h}_{u} $ is applied to the product of two terms, say $p_{\cd}(u-a)p_{\cd}(u-b)$, using the product rule for difference operators  we see that
\bea 
   && \De^{h}_{u}\{p_{\cd}(u-a)p_{\cd}(u-b)\}\label{3.24q}\\&&\qquad= \De^{h}_{u}\,p_{\cd}(u-a) p_{\cd}(u+h-b)+p_{\cd}(u-a) \De^{h}_{u}\, p_{\cd}(u-b)\nn 
 \eea
Consider an example of how the term $\De^{h}_{a}\De^{h}_{u}p_{\cd}(u-a)$ may appear. It could be  by the application
\begin{equation}
  \De^{h}_{a}  \(\De^{h}_{u}\,p_{\cd}(u-a) p_{\cd}(v-a)\),\label{3.25}
   \end{equation}
in which we take account of the two terms to which $\De^{h}_{a} $ is applied. Using the product rule in (\ref{3.24q}) we see that (\ref{3.25})
\begin{equation}
   =  \(\De^{h}_{a}   \De^{h}_{u}\, p_{\cd}(u-a) \)p_{\cd}(v -(a+h))+ \De^{h}_{u}\,p_{\cd}(u-a) \De^{h}_{a} p_{\cd}(v-a).\label{3.26}
   \end{equation}
Consider one more example
\bea 
  &&  \De^{h}_{a}  \(\De^{h}_{u}\,p_{\cd}(u-a)\De^{h}_{v}\, p_{\cd}(v-a)\)\label{3.27}\\
  &&\qquad\nn = \(\De^{h}_{a}   \De^{h}_{u}\, p_{\cd}(u-a) \)\De^{h}_{v}\, p_{\cd}(v -(a+h))\nn\\
  &&\hspace{1in}\nn + \De^{h}_{u}\,p_{\cd}(u-a) \De^{h}_{a}\De^{h}_{v} p_{\cd}(v-a).
 \eea
Note that in both examples the arguments of probability densities  with two difference operators applied to it does not contain an $h$. This is true in general because the difference 
formula, (\ref{3.24q}), does not add an  $  h$ to the argument of a term to which a difference operator is applied. Otherwise we may have a $\pm h$ added to the arguments of probability densities to which one  difference operator is applied, as in (\ref{3.27}), or to the arguments of probability densities to which no difference operator is applied, as in (\ref{3.26}).

  Based on the argument of the preceding paragraph we write (\ref{m.4a}) in the form
\be  
 E\(\prod_{\stackrel{j,k=0}{j< k}}^{K} \(I_{j,k,h} \)^{m_{j,k}}\)=\label{m.6} \sum_{a}\sum_{\pi_{0},\ldots, \pi_{K}}\int \mathcal{T}'_{h}( x;\,\pi ,a )\prod_{j,k,i}\,dx_{j,k,i},\label{4.25}
\ee  
where  
\begin{eqnarray}
\lefteqn{\mathcal{T}'_{h}( x;\,\pi ,a )=   \prod_{l=0}^{K}\int_{\RR_{l}}
\(\(\De^{ h}_{ x_{\pi_{l}(1)}}\)^{a_{ 1}(l,1)}
\(\De^{ h}_{ x_{\pi_{l-1}(m_{l-1})}}\)^{a_{ 2}(l,1)}\right.
\label{m.7}}\\
&&\hspace{1.5 in}
\left.\,p^{\sharp}_{(1-\sum_{q=1}^{m_{l-1}}r_{l-1,q})+r_{l,1}}(x_{\pi_{l}(1)}-x_{\pi_{l-1}(m_{l-1})})\)  \nonumber\\
&&\hspace{.1 in}\prod_{q=2}^{m_{l}}\(\(\De^{ h}_{ x_{ \pi_{l}( q)}}\)^{a_{ 1}(l,q)}
\(\De^{ h}_{ x_{ \pi_{l}( q-1)}}\)^{a_{ 2}(l,q)}\,p^{\sharp}_{r_{l,q}}(x_{\pi_{l}(q)}-x_{\pi_{l}(q-1)})\)
\prod_{q=1}^{m_{l}}\,dr_{l,q},\nonumber
\end{eqnarray}
  and where  
  $\RR_{l}=\{\sum_{q=1}^{m_{l}}r_{l,q}\leq 1\}$. In (\ref{4.25})
   the first sum is taken  over    all  
\be
a=(a_{ 1},a_{ 2})\,:\,\{(l,q),\, 0\leq l\leq K,\, 1\leq q\leq m_{l}\}\mapsto \{ 0,1\}\times \{ 0,1\}
\ee
  with the
restriction that for each triple $j,k,i$, there are exactly two  factors  of the form $\De^{
h}_{ x_{j,k,i}}$, each of which is applied to one of the  terms  $p^{\sharp}_{r_{\cdot}}(\cdot)$
that contains $x_{j,k,i}$ in its argument. This   condition can be stated more formally by saying that 
  for each $l$ and $q=1,\ldots, m_{l}-1$, if $\pi_{l}(q)=(j,k,i)$, then   $\{a_{ 1}(l,q), a_{ 2}(l,q+1) \}=\{0,1\}$ and if $q=m_{l}$ then   $\{a_{ 1}(l,m_{l}), a_{ 2}(l+1,1) \}=\{0,1\}$.    (Note that  when we write $\{a_{ 1}(l,q), a_{ 2}(l,q+1) \}=\{0,1\}$ we mean as two  sets, so, according to what $a$ is, we may have  $a_{ 1}(l,q)=1$ and $ a_{ 2}(l,q+1)=0$ or $a_{ 1}(l,q)=0$ and $a_{ 2}(l,q+1)=1$ and similarly for $\{a_{ 1}(l,m_{l}), a_{ 2}(l+1,1) \}$.)
  Also, in (\ref{m.7}) we define $(\De_{x_{i}}^{h})^{0}=1 $ and $(\De_{0}^{h}) =1 $.

  In (\ref{m.7}),    $p^{\sharp}_{r_{\cdot}}(z)$ can take any of the three values $p_{r_{\cdot}}(z)$, $p_{r_{\cdot}}(z+h)$,   or $p_{r_{\cdot}}(z- h)$. (We   must  consider all three possibilities.) Finally, it is important to emphasize that in (\ref{m.7}) each of  the difference operators is applied to only one of the  terms  $p^{\sharp}_{r_{\cdot}}(\cdot)$.

\medskip	  To avoid  confusion caused by the ambiguity of $p^{\sharp}_{\cd}$,  we first analyze 
\begin{equation}
\sum_{a}\sum_{\pi_{0},\ldots, \pi_{K}}\int \mathcal{T}_{h}( x;\,\pi ,a )\prod_{j,k,i}\,dx_{j,k,i}\label{m.8},
\end{equation}
where  
\begin{eqnarray}
&&\mathcal{T}_{h}( x;\,\pi ,a )=   \prod_{l=0}^{K}\int_{\RR_{l}}
\(\(\De^{ h}_{ x_{\pi_{l}(1)}}\)^{a_{ 1}(l,1)}
\(\De^{ h}_{ x_{\pi_{l-1}(m_{l-1})}}\)^{a_{ 2}(l,1)}\right.
\label{m.9}\\
&&\hspace{1.5 in}
\left.\,p_{(1-\sum_{q=1}^{m_{l-1}}r_{l-1,q})+r_{l,1}}(x_{\pi_{l}(1)}-x_{\pi_{l-1}(m_{l-1})})\)  \nonumber\\
&&\hspace{.1 in}\prod_{q=2}^{m_{l}}\(\(\De^{ h}_{ x_{ \pi_{l}( q)}}\)^{a_{ 1}(l,q)}
\(\De^{ h}_{ x_{ \pi_{l}( q-1)}}\)^{a_{ 2}(l,q)}\,p_{r_{l,q}}(x_{\pi_{l}(q)}-x_{\pi_{l}(q-1)})\)\prod_{q=1}^{m_{l}}\,dr_{l,q}.\nonumber
\end{eqnarray}
The difference between   $\mathcal{T} _{h}( x;\,\pi ,a )$ and $\mathcal{T}'_{h}( x;\,\pi ,a )$ is that in the former we replace $p^{\sharp}$ by $p$.  (I.e. we set $h=0$ in the arguments of the $p^{\sharp}$ terms in  (\ref{m.7}).) At the conclusion of this proof  we show  that both (\ref{m.8}) than (\ref{4.25}) have the same asymptotic limit as $h$ goes to zero.

 \medskip	We first obtain (\ref{m.2}). Let  $m=2n$, since $m_{j,k}=2n_{j,k}$,  $m_{l}=2n_{l}$ for some integer $n_{l}$. (Recall (\ref{m4.9})).  To begin we consider the case in which  $a=e$, where
 \be
 e(l, 2q)=( 1,1)\quad \mbox{and}\quad e(l, 2q-1)=( 0,0) \qquad \forall q.
 \ee
  When $a =e$ we have
\begin{eqnarray}
&&\mathcal{T}_{h}( x;\,\pi ,e )=   \prod_{l=0}^{K}\int_{\RR_{l}}p_{(1-\sum_{q=1}^{m_{l-1}}r_{l-1,q})+r_{l,1}}(x_{\pi_{l}(1)}-x_{\pi_{l-1}(m_{l-1})}) 
\nn\\
&&
\hspace{1.2 in}\prod_{q=2}^{n_{l}} p_{r_{l,2q-1}}(x_{\pi_{l}(2q-1)}-x_{\pi_{l}(2q-2)})    \label{m.10}\\
&&\hspace{1.2 in}\prod_{q=1}^{n_{l}} \De^{ h} 
\De^{ -h} \,p_{r_{l,2q}}(x_{\pi_{l}(2q)}-x_{\pi_{l}(2q-1)})\prod_{q=1}^{m_{l}}\,dr_{l,q}.\nn
\end{eqnarray}
  Here we use the following notation: $\De^{h} p (u-v) =p (u-v+h)-p (u-v)$, i.e., when $\De^{h}$ has no subscript, the difference operator is applied to the whole argument of the function. In this notation,
\begin{equation}
  \De_{u}^{h}  \De_{v}^{h} p (u-v) =  \De^{h}  \De^{-h} p (u-v) .
   \end{equation}

 \subsection{${\bf a =e}$, with all cycles of order two}\label{assinged}
 
    Consider the multigraph $G_{\pi }$ with vertices   
$\{(j,k,i),\,0\leq j<k\leq K,\, 1\leq i\leq m_{j,k}\}$.  Assign an edge 
between the vertices $\pi_{l } (2q-1)$ and $ \pi_{l}(2q)$  for each $ 0\leq l\leq K$ and $1\leq q\leq n_{l}$.  Each vertex is connected to two edges.  To see this suppose that $  \pi_{l}(2q)=\{(j,k,i)\}$,  with $j=l$ and    $ k=l'\neq l$, then there is a unique $q'$ such that   $  \pi_{l'}(2q')$ or $\pi_{l'}(2q'-1)$  is equal to  $ \{(j,k,i)\}$.  Therefore all the vertices lie in some cycle. Assume that there are $S$ cycles. We denote them by $C_{s}$, $s=1,\ldots, S$. 
Clearly, it is possible to have cycles of order two, in which case two vertices are connected by two edges. 

  It is important to note that the graph $G_{\pi}$ does not assign edges between $\pi_{l } (2q)$ and $ \pi_{l}(2q+1)$, although these vertices may be connected by   the edge assigned between $\pi_{l '} (2q'-1)$ and $ \pi_{l'}(2q')$ for some $l'$ and $q'$.

 	   We   estimate (\ref{m.9})  by breaking the calculation into two cases. In this section we consider the case when $a=e$ and all the cycles of $G_{\pi}$ are of order two. In Section \ref{sec-3.2} we consider the cases when $a=e$ and not all the cycles of $G_{\pi}$ are of order two,  and  when $a\ne e$.

 \medskip	
Let $\mathcal{P}=\{(\ga_{2v-1},\ga_{2v})\,,\,1\leq v\leq n\}$ be a pairing of the $m$ vertices   
\[\{(j,k,i),\,0\leq j<k\leq K,\, 1\leq i\leq m_{j,k}\}\] of $G_{\pi }$,  that satisfies the following special property: whenever $(j,k,i)$ and $(j',k',i')$ are paired together,   $j=j'$  and $k=k'$.    Equivalently,
\begin{equation}
\mathcal{P}=\bigcup_{\stackrel{j,k=0}{j< k}}^{K} \mathcal{P}_{j,k}\label{pair}
\end{equation}
where each $\mathcal{P}_{j,k}$ is a pairing of the $m_{j,k}$ vertices
\[\{(j,k,i),\, 1\leq i\leq m_{j,k}\}.\] 
We refer to such a pairing $\mathcal{P}$ as a special pairing and  denote the set of special pairings by $\mathcal{S}$.

Given a special  pairing $\mathcal{P}\in \mathcal{S}$, let $\pi$   be   such that for each $ 0\leq l\leq K$ and $1\leq q\leq n_{l}$, 

\be
\{\pi_{l}(2q-1), \pi_{l}(2q)\}=\{\ga_{2v-1},\ga_{2v}\}\label{4.34}
\ee
 for some, necessarily unique,   $ 1\leq v\leq n_{l}$. In this case we say that $\pi$ is compatible with the pairing $\mathcal{P}$ and write this  as $ \pi \sim \mathcal{P}$. (Recall that  when we write $\{\pi_{l}(2q-1), \pi_{l}(2q)\}=\{\ga_{2v-1},\ga_{2v}\}$, we mean as two  sets, so, according to what $\pi_{l}$ is, we may have  $\pi_{l}(2q-1)=\ga_{2v-1}$ and $ \pi_{l}(2q)=\ga_{2v}$ or $\pi_{l}(2q-1)=\ga_{2v}$ and $ \pi_{l}(2q)=\ga_{2v-1}$.)  Clearly
\begin{equation}
| \mathcal{S} |\leq {(2n)! \over 2^{n}n!}\label{ord.1}
\end{equation}
the number of pairings of $m=2n$ objects.

	Let $\pi\in \mathcal{B}$   be  such that $G_{\pi}$ consists of cycles of order two. It is easy to see that $ \pi \sim \mathcal{P}$ for some $\mathcal{P}\in \mathcal{S}$. To see this note that if $\{(j,k,i),(j',k',i')\}$  form a cycle of order two, there must exist   $l$ and $l'$ with $l\neq l'$ and $q$ and $q'$ such that both
$\{(j,k,i),(j',k',i')\}=\{\pi_{l} (2q-1), \pi_{l} (2q)\}$ and $\{(j,k,i),(j',k',i')\}=\{\pi_{l'} (2q'-1), \pi_{l'} (2q')\}$.   This implies that   $j=j'$, $k=k'$ and 
$\{j,k\}=\{l,l'\}$.  Furthermore, by (\ref{4.34}) we have
\begin{equation}
   \{\pi_{l}(2q-1), \pi_{l}(2q)\}=\{\pi_{l'} (2q'-1), \pi_{l'} (2q')\}=\{\ga_{2v-1},\ga_{2v}\}\label{4.36}
   \end{equation}

  When $\pi\sim \PP$ and all cycles are of order two we can write
\bea
  && \prod_{l=0}^{K}\prod_{q=1}^{n_{l}}  \De^{ h} 
\De^{ -h} \,p_{r_{l,2q}}(x_{\pi_{l}(2q)}-x_{\pi_{l}(2q-1)})\label{4.37} \\
&&\qquad=\prod_{v=1}^{n}\De^{h}\De^{ -h} \,p_{r_{2\nu}}( x_{\ga_{2v}}-x_{\ga_{2v-1}})\De^{h}\De^{ -h} \,p_{r'_{2\nu}}( x_{\ga_{2v}}-x_{\ga_{2v-1}})\nn,
   \eea
where $r_{2\nu}$ and $r'_{2\nu}$ are the rearranged indices $r_{l,2q}$ and $r_{l',2q'}$. We also use the fact that $\sum_{l=0}^{K}n_{l}=2n$.

For use in (\ref{91.4})   below we note that  
\bea 
&&
 \int_{0}^{1}\!\! \int_{0}^{1}| \De^{h}\De^{ -h} \,p_{r_{2\nu}}( x_{\ga_{2v}}-x_{\ga_{2v-1}})|\,|\De^{h}\De^{ -h} \,p_{r'_{2\nu}}( x_{\ga_{2v}}-x_{\ga_{2v-1}})|\,dr_{2\nu}\,dr'_{2\nu}\nn\\ && =\(\int_{0}^{1}|\De^{h}\De^{ -h} \,p_{r }( x_{\ga_{2v}}-x_{\ga_{2v-1}})|\,dr\)^{2}= w^{2}(x_{\ga_{2v}}-x_{\ga_{2v-1}}),\label{4.37j}
 \eea
(see (\ref{bl.3}).)  

\medskip	 We want to estimate the integrals in (\ref{m.8}). However, it is difficult   to integrate $\mathcal{T}_{h}( x;\,\pi ,e )$ directly, because   the variables, 
\bea
&&\{x_{\pi_{l}(1)}-x_{\pi_{l-1}(m_{l-1})},\,x_{\pi_{l}(2q-1)}-x_{\pi_{l}(2q-2)},\,x_{\pi_{l}(2q)}-x_{\pi_{l}(2q-1)};\nn \\&& \hspace{1.5in}
\,l\in [0,K],\, q\in [1,n_{l}]\}\nn,
\eea
are not independent.   We begin the estimation by showing that over much of the   domain of integration, the integral is negligible, asymptotically, as $h\to 0$. To begin, we write 
\begin{equation}
1=\prod_{ v=1}^{ n}\(1_{\{|x_{\ga_{2v}}-x_{\ga_{2v-1}}|\leq \sqrt{h}\}}+
1_{\{|x_{\ga_{2v}}-x_{\ga_{2v-1}}|\geq \sqrt{h}\}}\)\label{expm.1}
\end{equation}
and expand it as a sum of $2^{n}$ terms and use it to  write
\begin{eqnarray}
&&
\int \mathcal{T}_{h}( x;\,\pi,e)\prod_{j,k,i}\,dx_{j,k,i}\label{91.3a}\\
&&\qquad =\int \prod_{ v=1}^{ n}\(1_{\{|x_{\ga_{2v}}-x_{\ga_{2v-1}}|\leq \sqrt{h}\}}\)\mathcal{T}_{h}( x;\,\pi,e)\prod_{j,k,i}\,dx_{j,k,i}+E_{1,h}.\nonumber 
\end{eqnarray}
We now show that 
\begin{equation}
E_{1,h}=O\(h^{1/2}\(\frac{1}{h\psi^{2}(1/h)}\)^{n}\).\label{91.3}
\end{equation}
Note that every term in $E_{1,h}$   can be written in the form 
\begin{equation}
B_{h}(\pi,e, D ):=\int \prod_{v=1}^{n} 1_{ D_{v} } 
\mathcal{T}_{h}( x;\,\pi,e)\prod_{j,k,i}\,dx_{j,k,i}\label{f9.33}
\end{equation}
where each   $D_{v}$ is either $\{|x_{\ga_{2v}}-x_{\ga_{2v-1}}|\leq \sqrt{h}\}$ or $\{|x_{\ga_{2v}}-x_{\ga_{2v-1}}|\geq \sqrt{h}\}$, and at least one of the $D_{v}$ is of the second type.  

 Consider (\ref{f9.33}) and the representation of $\mathcal{T}_{h}( x;\,\pi,e)$ in (\ref{m.10}). We  
 take absolute values in the integrand in  (\ref{m.10})   and
  take all the integrals with $r_{\cdot}$  between 0 and 1  and use (\ref{4.37j}) followed by (\ref{2.2})  
   to get   
\begin{eqnarray} 
 | B_{h}(\pi,e, D )|&\leq &
 \int \prod_{v=1}^{n} 1_{ D_{v} }  w^{2}( x_{\ga_{2v}}-x_{\ga_{2v-1}})  \prod_{l=0}^{K}u(x_{\pi_{l}(1)}-x_{\pi_{l-1}(m_{l-1})})   \nonumber\\
&& \qquad
  \prod_{q=2}^{n_{l}} u(x_{\pi_{l}(2q-1)}-x_{\pi_{l}(2q-2)}) \prod_{j,k,i}\,dx_{j,k,i} \label{91.4}.
\end{eqnarray}
 We now take 
\be
\{x_{\ga_{2v}}-x_{\ga_{2v-1}},\,v=1,\ldots,n \}
\ee
 and an additional $n$ 
variables from  the $2n $ arguments of the $u$ terms,
\be \cup_{l=0}^{K} \{x_{\pi_{l}(1)}-x_{\pi_{l-1}(m_{l-1})}, x_{\pi_{l}(2q-1)}-x_{\pi_{l}(2q-2)},\,q=2,\ldots,n_{l} \}\label{kik}\ee
so that the chosen $2n$ variables generate   the space spanned by the 2n variables $\{x_{j,k,i} \}$. 
There are $n$ variables in (\ref{kik}) that are not used.  We bound  the functions $u$  of these variables by their sup norm, which by (\ref{2.2}) is finite. Then we make a change of variables and get that  
\bea 
  |B_{h}(\pi,e, D )|&\le& \int \prod_{v=1}^{n} 1_{ D_{v} }  w^{2}( y_{v})\prod_{v=n+1}^{2n}u(y_{v})\prod_{v=1}^{2n}\,dy_{v}
\label{91.oo}\\
&\le&C\int \prod_{v=1}^{n} 1_{ D_{v} }  w^{2}( y_{v}) \prod_{v=1}^{n}\,dy_{v}\nn,\\
&=&O\(h^{1/2}\(\frac{1}{h\psi^{2}(1/h)}\)^{n}\). \nonumber
\eea
  Here we use (\ref{2.2}) to see that the integrals of the $u$ terms is finite.
 Then we use (\ref{bl.5}) and  (\ref{bl.5b}) to obtain (\ref{91.3}).   (Note that it is because at least one of the $D_{v}$ is of the second type that we can use (\ref{bl.5b}).)

\medskip	We now study  
\begin{equation} \hspace{.4 in}
\int \prod_{ v=1}^{ n}\(1_{\{|x_{\ga_{2v}}-x_{\ga_{2v-1}}|\leq \sqrt{h}\}}\)\mathcal{T}_{h}( x;\,\pi,e)\prod_{j,k,i}\,dx_{j,k,i}. \label{f91.36}
\end{equation}
 Recall that for each $ 0\leq l\leq K$ and $1\leq q\leq n_{l}$, 
$\{\pi_{l}(2q-1), \pi_{l}(2q)\}=\{\ga_{2v-1},\ga_{2v}\}$ for some $ 1\leq v\leq n$.  We identify these relationships by setting $v=\si_{l} (q) $ when  $\{\pi_{l}(2q-1), \pi_{l}(2q)\}=\{\ga_{2v-1},\ga_{2v}\}$,   and sometimes write this last term as $\{\ga_{2\si_{l} (q)-1},\ga_{2\si_{l} (q)}\}$.

 For $q\geq 2$ we write 
\begin{eqnarray}
\lefteqn{p_{r_{l,2q-1}}(x_{\pi_{l}(2q-1)}-x_{\pi_{l}(2q-2)})
  \label{f91.37}}\\
&&=p_{r_{l,2q-1}}(x_{\ga_{2\si_{l} (q)-1}}-x_{\ga_{2\si_{l} (q-1)-1}})+\De^{h_{l,q}} p_{r_{l,2q-1}}(x_{\ga_{2\si_{l} (q)-1}}-x_{\ga_{2\si_{l} (q-1)-1}}), \nn
\end{eqnarray}
where $h_{l,q}=(x_{\pi_{l}(2q-1)}-x_{\ga_{2\si_{l} (q)-1}})+(x_{\ga_{2\si_{l} (q-1)-1}}-x_{\pi_{l}(2q-2)})$.  When $q=1$  we can make a similar decomposition
\begin{eqnarray}
&&p_{(1-\sum_{q=1}^{m_{l-1}}r_{l-1,q})+r_{l,1}}(x_{\pi_{l}(1)}-x_{\pi_{l-1}(m_{l-1})})
\label{f91.37.1}\\
&&\qquad= p_{(1-\sum_{q=1}^{m_{l-1}}r_{l-1,q})+r_{l,1}}(x_{\ga_{2\si_{l} (1)-1}}-x_{\ga_{2\si_{l-1} (n_{l-1})-1}})  \nonumber\\
&&\qquad\quad+\De^{h_{l,1}} p_{(1-\sum_{q=1}^{m_{l-1}}r_{l-1,q})+r_{l,1}}(x_{\ga_{2\si_{l} (1)-1}}-x_{\ga_{2\si_{l-1} (n_{l-1})-1}}) , \nonumber
\end{eqnarray}
where $h_{l,1}=(x_{\pi_{l}(1)}-x_{\ga_{2\si_{l} (1)-1}})+(x_{\ga_{2\si_{l-1} (n_{l-1})-1}})-x_{\pi_{l-1}(m_{l-1})})$.
Note that because of the presence of  the term $\prod_{ v=1}^{ n}\(1_{\{|x_{\ga_{2v}}-x_{\ga_{2v-1}}|\leq \sqrt{h}\}}\)$ in the integral in (\ref{f91.36})  we need only be concerned with values of $|h_{l,q}|\leq 2\sqrt{h}$,  for $ 0\leq l\leq K$ and   $1\leq q\leq n_{l}$.

  For $q=1,\ldots,n_{l}$, $l=0\ldots,K$, we substitute   (\ref{f91.37}) and (\ref{f91.37.1})   into the term $\mathcal{T}_{h}( x;\,\pi ,e )$ in  (\ref{f91.36}), (see also (\ref{m.10})),   and expand the products so that we can write 
(\ref{f91.36}) as a sum of $2^{\sum_{l=0}^{K}n_{l} }$  terms, which we write as  
\begin{eqnarray}
&&\int \prod_{ v=1}^{ n}\(1_{\{|x_{\ga_{2v}}-x_{\ga_{2v-1}}|\leq \sqrt{h}\}}\)\mathcal{T}_{h}( x;\,\pi,e)\prod_{j,k,i}\,dx_{j,k,i}
\label{91.6}\\
&&\qquad=\int \prod_{ v=1}^{ n}\(1_{\{|x_{\ga_{2v}}-x_{\ga_{2v-1}}|\leq \sqrt{h}\}}\)\mathcal{T}_{h,1}( x;\,\pi,e)\prod_{j,k,i}\,dx_{j,k,i} + E_{2,h},\nonumber
\end{eqnarray}
  where
\begin{eqnarray}
&&\mathcal{T}_{h,1}( x;\,\pi ,e )=   \prod_{l=0}^{K}\int_{\RR_{l}}p_{(1-\sum_{q=1}^{m_{l-1}}r_{l-1,q})+r_{l,1}}(x_{\ga_{2\si_{l} (1)-1}}-x_{\ga_{2\si_{l-1} (n_{l-1})-1}}) 
\nn\\
&&
\hspace{1.4 in}\prod_{q=2}^{n_{l}} p_{r_{l,2q-1}}(x_{\ga_{2\si_{l} (q)-1}}-x_{\ga_{2\si_{l} (q-1)-1}})  \label{m.11}\\
&&\hspace{1.6 in} \prod_{q=1}^{n_{l}} \De^{ h} 
\De^{ -h} \,p_{r_{l,2q}}(x_{\pi_{l}(2q)}-x_{\pi_{l}(2q-1)})\prod_{q=1}^{m_{l}}\,dr_{l,q}.\nn
\end{eqnarray} 
 Using  (\ref{4.37}) we can rewrite this as  
 \begin{eqnarray}
&&\mathcal{T}_{h,1}( x;\,\pi ,e )\label{m.11c}\\
&&\qquad
=   \int_{\RR_{0}\times\cdots\times \RR_{K}}\(\prod_{l=0}^{K} p_{(1-\sum_{q=1}^{m_{l-1}}r_{l-1,q})+r_{l,1}}(x_{\ga_{2\si_{l} (1)-1}}-x_{\ga_{2\si_{l-1} (n_{l-1})-1}})\right. 
\nn\\
&&\left. 
\hspace{2 in}\prod_{q=2}^{n_{l}} p_{r_{l,2q-1}}(x_{\ga_{2\si_{l} (q)-1}}-x_{\ga_{2\si_{l} (q-1)-1}})\)  \nn\\ 
&&\hspace{.1 in}\( \prod_{v=1}^{n}\De^{h}\De^{ -h} \,p_{r_{2\nu}}( x_{\ga_{2v}}-x_{\ga_{2v-1}})\De^{h}\De^{ -h} \,p_{r'_{2\nu}}( x_{\ga_{2v}}-x_{\ga_{2v-1}}) \)\nn\\
&&\hspace{3.5in}\prod_{l=0}^{K}\prod_{q=1}^{m_{l}}\,dr_{l,q}\nn,
\end{eqnarray}
where $r_{2\nu}$ and $r'_{2\nu}$ are the rearranged indices $r_{l,2q}$ and $r_{l',2q'}$. Since the variables $x_{\ga_{2v}},\,v=1,\ldots, n$, occur only in the last line of (\ref{m.11c}),    we make the change of variables $ x_{\ga_{2v}}-x_{\ga_{2v-1}}\to  x_{\ga_{2v}} $ and $x_{\ga_{2v-1}}\to x_{\ga_{2v-1}}$ and get that 
\begin{eqnarray}
\lefteqn{\int  \mathcal{T}_{h,1}( x;\,\pi,e)\prod_{j,k,i}\,dx_{j,k,i}
\label{m.11d}}\\
&&
= \int  \int_{\RR_{0}\times\cdots\times \RR_{K}}
\(\prod_{l=0}^{K} p_{(1-\sum_{q=1}^{m_{l-1}}r_{l-1,q})+r_{l,1}}(x_{\ga_{2\si_{l} (1)-1}}-x_{\ga_{2\si_{l-1} (n_{l-1})-1}}) \right.\nn\\
&&\hspace{1.9in}\left. \prod_{q=2}^{n_{l}} p_{r_{l,2q-1}}(x_{\ga_{2\si_{l} (q)-1}}-x_{\ga_{2\si_{l} (q-1)-1}})\)  \nn\\
&&\qquad \( \prod_{v=1}^{n}\De^{h}\De^{ -h} \,p_{r_{2\nu}}( x_{\ga_{2v}})\De^{h}\De^{ -h} \,p_{r'_{2\nu}}( x_{\ga_{2v}}) \)\prod_{l=0}^{K}\prod_{q=1}^{m_{l}}\,dr_{l,q}
\prod_{j,k,i}\,dx_{j,k,i}.\nonumber
\end{eqnarray}
 Since the variables $x_{\ga_{2v}},\,v=1,\ldots, n$ occur only in the last line of (\ref{m.11d}) and the variables $x_{\ga_{2v-1}},\,v=1,\ldots, n$ occur only in the second and third lines of (\ref{m.11d}), we can  write (\ref{m.11d}) as
\begin{eqnarray}
\lefteqn{\int \mathcal{T}_{h,1}( x;\,\pi,e)\prod_{j,k,i}\,dx_{j,k,i}
\label{m.11e}}\\
&&
=   \int_{\RR_{0}\times\cdots\times \RR_{K}} \int\(\prod_{l=0}^{K} p_{(1-\sum_{q=1}^{m_{l-1}}r_{l-1,q})+r_{l,1}}(x_{\ga_{2\si_{l} (1)-1}}-x_{\ga_{2\si_{l-1} (n_{l-1})-1}})\right. 
\nn\\
&&\left. 
\hspace{1.5 in}\prod_{q=2}^{n_{l}} p_{r_{l,2q-1}}(x_{\ga_{2\si_{l} (q)-1}}-x_{\ga_{2\si_{l} (q-1)-1}})\) 
 \prod_{v=1}^{n} \,dx_{\ga_{2v-1}}\nn\\
&&\hspace{.5 in}\( \prod_{v=1}^{n}\int \De^{h}\De^{ -h} \,p_{r_{2\nu}}( x_{\ga_{2v}})\De^{h}\De^{ -h} \,p_{r'_{2\nu}}( x_{\ga_{2v}}) \,dx_{\ga_{2v}}\)\prod_{l=0}^{K}\prod_{q=1}^{m_{l}}\,dr_{l,q}.\nn
\end{eqnarray}
 Note that we also use Fubini's Theorem which is    justified since the    absolute value of the integrand is integrable, (as we point out in the argument preceding (\ref{91.4})).   (In the rest of this section use Fubini's Theorem frequently for  integrals like (\ref{m.11e})  without repeating the explanation about why it is justified.)

 Analogous to (\ref{91.3a}) we note that 
\begin{eqnarray}
&&
\int \mathcal{T}_{h,1}( x;\,\pi,e)\prod_{j,k,i}\,dx_{j,k,i}\label{91.3ab}\\
&&\qquad =\int \prod_{ v=1}^{ n}\(1_{\{|x_{\ga_{2v}}-x_{\ga_{2v-1}}|\leq \sqrt{h}\}}\)\mathcal{T}_{h,1}( x;\,\pi,e)\prod_{j,k,i}\,dx_{j,k,i}+\wt E_{1,h},\nonumber 
\end{eqnarray}
where $\wt E_{1,h}=O\(h^{1/2}\(\frac{1}{h\psi^{2}(1/h)}\)^{n}\)$. The proof of (\ref{91.3ab}) is the same as the proof of  (\ref{91.3}). 

\medskip	
 We now show that   \be
 E_{2,h}=O\(\({1 \over h\psi (1/h)}\)^{1/2}\(\frac{1}{h\psi^{2}(1/h)}\)^{n}\).\label{3.24}
 \ee
 To see this note that   the terms in $E_{2,h}$ are of the form   
\begin{eqnarray}
&& 
\int  \prod_{ v=1}^{ n}\(1_{\{|x_{\ga_{2v}}-x_{\ga_{2v-1}}|\leq \sqrt{h}\}}\) \label{f9.40}\\
&&\hspace{.5in} \prod_{l=0}^{K}\int_{\RR_{l}}\wt p_{(1-\sum_{q=1}^{m_{l-1}}r_{l-1,q})+r_{l,1}}(x_{\ga_{2\si_{l} (1)-1}}-x_{\ga_{2\si_{l-1} (n_{l-1})-1}}) 
\nn\\
&&
\hspace{.7 in}\prod_{q=2}^{n_{l}}\wt  p_{r_{l,2q-1}}(x_{\ga_{2\si_{l} (q)-1}}-x_{\ga_{2\si_{l} (q-1)-1}})    \nonumber\\
&&\hspace{.9 in}\prod_{q=1}^{n_{l}} \De^{ h} 
\De^{ -h} \,p_{r_{l,2q}}(x_{\pi_{l}(2q)}-x_{\pi_{l}(2q-1)})\prod_{q=1}^{m_{l}}\,dr_{l,q}\prod_{j,k,i}\,dx_{j,k,i}\nn,
\end{eqnarray}
where $\wt p_{r_{l,2q-1}}$ is either $  p_{r_{l,2q-1}}$ or $\De^{h_{l,q}} p_{r_{l,2q-1}}$. Furthermore,   at least one of the terms  $\wt p_{r_{l,2q-1}}$  is of the form   a $\De^{h_{l,q}} p_{r_{l,2q-1}}$.  

As in the transition from (\ref{f9.33}) to (\ref{91.4}) we  bound the absolute value of (\ref{f9.40})    by  
\begin{eqnarray}
\lefteqn{\int \prod_{ v=1}^{ n}\(1_{\{|x_{\ga_{2v}}-x_{\ga_{2v-1}}|\leq \sqrt{h}\}}\)w^{2}( x_{\ga_{2v}}-x_{\ga_{2v-1}}) 
\label{m.12}}\\
&&
\prod_{l=0}^{K}\wt u(x_{\ga_{2\si_{l} (1)-1}}-x_{\ga_{2\si_{l-1} (n_{l-1})-1}}) \prod_{q=2}^{n_{l}} \wt u(x_{\ga_{2\si_{l} (q)-1}}-x_{\ga_{2\si_{l} (q-1)-1}}) 
  \prod_{j,k,i}\,dx_{j,k,i},\nn
\end{eqnarray} where each $\wt u$ is either of the form $u$ or $v$, in Lemma \ref{lem-vproprvt}, and where, obviously, the $h$ in (\ref{bl.2}) is $h_{l,q}$. Furthermore, we have $J$ terms of the type $v$, for some   $J\geq 1$. It follows from (\ref{bl.2}), the regular variation of $\psi$   and the fact that $|h_{l,q}|\leq 2\sqrt{h}$, that
\begin{equation}
   v(\cd)\le  C\({1 \over h\psi (1/h)}\)^{1/2}\frac{1}{1+x^{2}}
   \end{equation}
 Using this and (\ref{2.2}) we can bound the  integral in (\ref{m.12}) by 
\begin{eqnarray}
\lefteqn{
 C\({1 \over h\psi (1/h)}\)^{J/2}\int \prod_{ v=1}^{ n} w^{2}( x_{\ga_{2v}}-x_{\ga_{2v-1}}) 
	\label{m.13}}\\
&&   \hspace{.1 in}\prod_{l=0}^{K}\bar u(x_{\ga_{2\si_{l} (1)-1}}-x_{\ga_{2\si_{l-1} (n_{l-1})-1}}) \prod_{q=2}^{n_{l}} \bar u(x_{\ga_{2\si_{l} (q)-1}}-x_{\ga_{2\si_{l} (q-1)-1}}) 
 \prod_{j,k,i}\,dx_{j,k,i} \nonumber
\end{eqnarray}
where all the terms  $\bar u(y)=(1+y^{2})^{-1}$.

 Since   the variables $x_{\ga_{2\nu}}$, $\nu=1,\ldots, n$, occur only in the $w$ terms in (\ref{m.13}) and the variables $x_{\ga_{2v-1}},\,v=1,\ldots, n$ occur only in the $\ov u$ terms in (\ref{m.13})  , (refer to the change of variables arguments in (\ref{m.11d}) and (\ref{m.11e})), we can  write (\ref{m.13}) as 
\begin{eqnarray}
\lefteqn{ C\({1 \over h\psi (1/h)}\)^{J/2}\int \(\prod_{l=0}^{K}\ov u(x_{\ga_{2\si_{l} (1)-1}}-x_{\ga_{2\si_{l-1} (n_{l-1})-1}}) \right.
	\label{m.13ww}}\\
&&  \left. \hspace{.1 in}\prod_{q=2}^{n_{l}} \ov u(x_{\ga_{2\si_{l} (q)-1}}-x_{\ga_{2\si_{l} (q-1)-1}}) \)\prod_{v=1}^{n}\,dx_{\ga_{2v-1}}\prod_{v=1}^{n}w^{2}(x_{\ga_{2v}} )\prod_{v=1}^{n}\,dx_{\ga_{2v}}.
 \nn 
\end{eqnarray}
As we have been doing we extract a linearly independent set of variables from the arguments of the $\ov u$ terms. The other $\ov u$ terms we bound by one. Then we make a change of variables and integrate the remaining $\ov u$ terms and the $w^{2}$ terms  using    (\ref{2.2}) and  (\ref{bl.5}). 
Since $J\geq 1$, we get (\ref{3.24}).

\medskip	
Since $\psi$ is regularly varying with index $\bb>1$ we see that there exists an $\ep>0$ such that
\begin{equation}
   E_{1,h}+   E_{2,h}+ \wt E_{2,h}=O\(h^{(2\bb-1)n+\ep}\).
   \end{equation}
Therefore, it follows from (\ref{91.3a}),  (\ref{91.6})   and (\ref{91.3ab})
  that 
\begin{eqnarray}
&&\int  \mathcal{T}_{h}( x;\,\pi,e)\prod_{j,k,i}\,dx_{j,k,i}
\label{m.13q}\\
&&\qquad=\int \mathcal{T}_{h,1}( x;\,\pi,e)\prod_{j,k,i}\,dx_{j,k,i}+O\(h^{(2\bb-1)n+\ep}\).\nonumber
\end{eqnarray}

\medskip	 
Let $\wt\RR_{l}(s )=\{\sum_{q=1}^{n_{l}}r_{l,2q-1}\leq 1-s \}$ and $\wt  \si_{l} (q):=\ga_{2\si_{l} (q)-1}$.   We define
\begin{eqnarray} \lefteqn{
F (\wt  \si,s_{0},\ldots, s_{K}) 
\label{m.14}}\\
&& =\int\(\int_{\wt\RR_{0}(s_{0})\times\cdots \times\wt\RR_{K}(s_{K})} \prod_{ l=0}^{ K}\, p_{(1-\sum_{q=1}^{n_{l-1}}r_{l-1,2q-1}-s_{l-1})+r_{l,1}}\right.
\nn\\
&&\left.
\hspace{.2 in}(x_{\wt  \si_{l} (1) }-x_{\wt  \si_{l-1} (n_{l-1}) }) \prod_{q=2}^{n_{l}} p_{r_{l,2q-1}}(x_{\wt  \si_{l} (q) }-x_{\wt  \si_{l} (q-1) })
\prod_{q=1}^{n_{l}}\,dr_{l,2q-1}\)  \,d x    \nonumber,
\end{eqnarray}
where $(1-\sum_{q=1}^{n_{-1}}r_{-1,2q-1}-s_{-1}):=0$ and $ \wt  \si_{ -1} (n_{ -1}):=0$. 
Here the generic term $dx$ indicates integration with respect to all the variables $x_{\cdot}$ that appear in the integrand.

Since   $\wt  \si_{l} (q)=\ga_{2\si_{l} (q)-1}$ we can also write (\ref{m.14}) as
\begin{eqnarray} \lefteqn{
F (\wt\si,s_{0},\ldots, s_{K})\label{m.14a}}\\
&& =\int\(\int_{\wt\RR_{0}(s_{0})\times\cdots \times\wt\RR_{K}(s_{K})} \prod_{ l=0}^{ K}\, \right.
p_{(1-\sum_{q=1}^{n_{l-1}}r_{l-1,2q-1}-s_{l-1})+r_{l,1}}
\nn\\
&& \hspace{.3in}
(x_{\ga_{2\si_{l} (1)-1}}-x_{\ga_{2\si_{l-1} (n_{l-1})-1}}) \prod_{q=2}^{n_{l}} p_{r_{l,2q-1}}(x_{\ga_{2\si_{l} (q)-1}}-x_{\ga_{2\si_{l} (q-1)-1}})
  \nonumber\\
&&\hspace{3
in}\left.\prod_{q=1}^{n_{l}}\,dr_{l,2q-1}\)  \,dx, \nn
\end{eqnarray} 
  $x_{\ga_{2\si_{-1} (n_{-1})-1}}:=0$.

\medskip	  Consider the mappings  $\wt  \si_{l}$  that are used in (\ref{m.14}).   Recall that $\si_{l} (q) $ is defined by the relationship  $\{\pi_{l}(2q-1), \pi_{l}(2q)\}=\{\ga_{2\si_{l} (q)-1},\ga_{2\si_{l} (q)}\}$.   Therefore,    since $\wt  \si_{l}(q)=\ga_{2\si_{l} (q)-1}$
  we can have that either  $\wt\si_{l} (q) =\pi_{l}(2q-1)$ or
$\wt\si_{l} (q) =\pi_{l}(2q)$.  
 However, since the terms $\wt  \si_{l}(q)$ are subscripts of the terms $x$, all of which are integrated, it is more convenient to define   $\wt  \si_{l}$ differently. 

Recall that $\mathcal{P}$,  (see (\ref{pair})), is a union of pairings $\mathcal{P}_{j,k}$
 of the $m_{j,k}$ vertices
\[\{(j,k,i),\, 1\leq i\leq m_{j,k}\}.\] 
Each $\mathcal{P}_{j,k}$ consists of $n_{j,k}$ pairs, that  can ordered arbitrarily. If 
$\{\ga_{2\si_{l} (q)-1},\newline \ga_{2\si_{l} (q)}\}$ is the $i$-th pair in $\mathcal{P}_{j,k}$, we set 
$\wt  \si_{l}(q)=(j,k,i)$. (Necessarily, $l$ will be either $j$ or $k$, as we point out in the paragraph containing (\ref{4.36})). Thus, each $\wt  \si_{l}$ is a bijection from $[1,n_{l}]$ to
\begin{equation}
\wt I_{l}=\bigcup_{k=l+1}^{K}\{(l,k,i),\,1\leq i\leq n_{l,k}\}\bigcup_{j=0}^{l-1}\{(j,l,i),\,1\leq i\leq n_{j,l}\}.\label{m.5kp}
\end{equation}
Let $\wt \mathcal{B}$  denote the set of $K+1$ tuples, $\wt\si=(\wt\si_{0},\ldots, \wt\si_{K})$ of such bijections.  Note that with this definition of $\wt  \si_{l}(q)$  (\ref{m.14})  remains  unchanged since we have simply renamed the variables of integration.

 By interchanging the elements in any of the $2n$ pairs $\{\pi_{l}(2q-1),\pi_{l}(2q)\}$ 
 we obtain a new $\pi' \sim \mathcal{P}$. In fact we obtain $2^{2n}$ different permutations $\pi$, in this way,  all of which are compatible with $\mathcal{P}$, and all of which   give   the same $\wt\si$ in (\ref{m.14}). Furthermore, by permuting the 
 pairs $\{\pi_{l}(2q-1),\pi_{l}(2q)\}$, $1\leq q\leq n_{l}$, for each $l$, we get all the possible permutation $\wt\pi\sim \mathcal{P}$,  and these give   all possible  mappings $\wt\si\in \wt \mathcal{B}$. Note that $|\wt\mathcal{B}|= \prod_{l=0}^{K}n_{l}!\leq (2n)!$.
 
 \medskip	Consider (\ref{m.14a}). Since $x_{\ga_{2\si_{-1} (n_{-1})-1}}=0$, $x_{\ga_{2\si_{0} -1}}$ appears alone as the argument of one of the density functions. Therefore we can extract a linearly independent set from the arguments of the densities that spans the space spanned by all the arguments of the densities. We use 
  (\ref{2.dens}) to bound the density functions with arguments that are not in the spanning set by $C\psi^{-1}(1/s)$.  We then integrate them with respect to the time variables. Since the time variables are bounded,  all this contributes only some constant.  With what is left we can make a change of variables and use 
  (\ref{2.dens}) again to see that
   \begin{equation}
   F (\wt\si,s_{0},\ldots, s_{K})\label{m.14p}\le C,\label{4.60}
   \end{equation}
 for some constant depending only on $m$.

\medskip	
  Let $\wh\RR_{l} =\{\sum_{q=1}^{n_{l}}r_{l,2q}\leq 1\}$.  We break up the integration over $\RR_{l}$ into integration over $\wt\RR_{l}(s )$ and $\wh\RR_{l}$ in  (\ref{m.11e})
and  (\ref{m.14a}).  If one carefully examines  the time indices in (\ref{m.9}) and (\ref{m.14}) and uses Fubini's Theorem, one sees that  
  \begin{eqnarray} 
&& \int \mathcal{T}_{h}( x;\,\pi,e)\prod_{j,k,i}\,dx_{j,k,i}\label{m.15}\\ 
&&\qquad=
 \int_{\wh\RR_{0}\times\cdots \times\wh\RR_{K}}    
F (\wt\si,\sum_{q=1}^{n_{0}}r_{0,2q}\,,\ldots, \sum_{q=1}^{n_{K}}r_{K,2q})  \nn\\ 
  &&
 \hspace{.6 in} \prod_{i=1}^{n} \(\int \(  \De^{ h}\De^{ -h}
\,p_{r_{i} }(x ) \)  \(\De^{ h}\De^{ -h}
\,p_{r'_{i} }(x )\) \,dx\) \prod_{i=1}^{n}\,dr_{i}\,dr'_{i}.\nn
\end{eqnarray}
The variables $\{r_{i},r'_{i}\,|\,i=1,\ldots,n\}$ are simply a relabeling of the variables $\{r_{l,2q}\,|\,0\leq l\leq K,\,1\leq q\leq n_{l}\}$.
(The exact form of this relabeling  does not matter in what follows.) Here, as always, we set $p_{r}(x)=0$, if $r\leq 0$.

\medskip	 By   Parseval's Theorem,
\begin{eqnarray}
&&\int \(\De^{ h}\De^{ -h}\,p_{r  }(x)\)
\(\De^{ h}\De^{ -h}\,p_{r' }(x)\) \,dx
\label{f9.2}\\
&&\qquad ={1 \over 2\pi}\int |2-e^{iph}-e^{-iph}  |^{2}e^{-r\psi (p) } e^{-r'\psi (p) }\,dp\geq 0.  \nonumber
\end{eqnarray}
  Using this,   (\ref{4.60}) and Fubini's Theorem, we see  that  \begin{eqnarray}
\lefteqn{ \int_{\(\wh\RR_{0}\times\cdots \times\wh\RR_{K}\)\cap ([0,\sqrt{h}]^{2n})^{c}}    
F (\wt\si,\sum_{q=1}^{n_{0}}r_{0,2q}\,,\ldots, \sum_{q=1}^{n_{K}}r_{K,2q})
\label{parm.1}}\\
&&\qquad   \prod_{i=1}^{n} \(\int \(  \De^{ h}\De^{ -h}
\,p_{r_{i} }(x ) \)  \(\De^{ h}\De^{ -h}
\,p_{r'_{i} }(x )\) \,dx\)\prod_{i=1}^{n}\,dr_{i}\,dr'_{i}\nn\\
&& \leq C\int_{ ([0,\sqrt{h}]^{2n})^{c}}    
   \prod_{i=1}^{n} \(\int \(  \De^{ h}\De^{ -h}
\,p_{r_{i} }(x ) \)  \(\De^{ h}\De^{ -h}
\,p_{r'_{i} }(x )\) \,dx\)\prod_{i=1}^{n}\,dr_{i}\,dr'_{i}\nn\\
&& \leq C  \(\int \(\int ( \De^{ h}\De^{ -h}
\,p_{r  }(x ) )  \,dr \)^{2}  \,dx \)^{n-1}\nn\\
&&\qquad\quad\int \bigg\{\int_{0}^{\ff}\int_{\sqrt h}^{\ff} \(  \De^{ h}\De^{ -h}
\,p_{r_{i} }(x ) \)  \(\De^{ h}\De^{ -h}
\,p_{r'_{i} }(x )\)\,dr_{i}\,dr'_{i} \bigg\}\,dx
\nn\\
&&= C c^{n-1}_{\psi,h,1} \int \bigg\{\int_{0}^{\ff}\int_{\sqrt h}^{\ff} \(  \De^{ h}\De^{ -h}
\,p_{r_{i} }(x ) \)  \(\De^{ h}\De^{ -h}
\,p_{r'_{i} }(x )\)\,dr_{i}\,dr'_{i} \bigg\}\,dx,
\nn
\eea by (\ref{jay.31}).
The integral in the final line of  (\ref{parm.1})
\bea 
   \le  c_{\psi,h,1} - \int  \(\int_{0}^{\sqrt{h}}\De^{ h}\De^{ -h}\,\,p_{s }(x)\,ds\)^{2}\,dx.
 \eea
Therefore, it follows from  Lemma \ref{lem-h3} that the first integral in (\ref{parm.1})
is  
   $O(h^{(2\bb-1)n+\ep}),
$
for some $\ep>0$.

 \medskip	
Since  $\(\wh\RR_{0}\times\cdots \times\wh\RR_{K}\)\supseteq [0,\sqrt{h}]^{2n}$ for $2n\sqrt{h}\leq 1$, it follows from (\ref{m.15})  and the preceding sentence, that  
 \begin{eqnarray} 
&&\lefteqn{ \int \mathcal{T}_{h}( x;\,\pi,e)\prod_{j,k,i}\,dx_{j,k,i}\label{m.15s}}\\ 
&&\qquad=
 \int_{[0,\sqrt{h}]^{2n}}    
F (\wt\si,\sum_{q=1}^{n_{0}}r_{0,2q}\,,\ldots, \sum_{q=1}^{n_{K}}r_{K,2q})  \prod_{i=1}^{n} \bigg(\int \(  \De^{ h}\De^{ -h}
\,p_{r_{i} }(x ) \)   \nn\\ 
&&\hspace{.7in}
  \(\De^{ h}\De^{ -h}
\,p_{r'_{i} }(x )\) \,dx\bigg)\prod_{ l=0}^{ K}\prod_{q=1}^{n_{l}}\,dr_{l,2q}+O(h^{(2\bb-1)n+\ep}).\nn
\end{eqnarray}
 
We use the next lemma which is proved in  Subsection \ref{sec-m4}. 

\begin{lemma}\label{lem-m4}  For any fixed $m$ and  $s_{0},\ldots, s_{K}\leq m\sqrt{h}$,  there exists an $\ep>0$ such that for all $h> 0$, sufficiently small, 
\begin{equation}
|F (\wt\si,s_{0},\ldots, s_{K})-F (\wt\si,0,\ldots, 0)|\leq Ch^{\ep}.\label{m4}
\end{equation}
  \end{lemma} 
  
 \noindent{\bf Proof of Lemma \ref{lem-multiple} continued }  It follows from  (\ref{m.15s}) and Lemmas  \ref{lem-m4} and   \ref{lem-h3}, that 
   \begin{eqnarray} 
&&\lefteqn{ \int \mathcal{T}_{h}( x;\,\pi,e)\prod_{j,k,i}\,dx_{j,k,i}\label{m.15sw}}\\ 
&&\qquad=F (\wt\si, 0,\ldots,0) 
 \int_{[0,\sqrt{h}]^{2n}}    
 \prod_{i=1}^{n} \bigg(\int \(  \De^{ h}\De^{ -h}
\,p_{r_{i} }(x ) \)   \nn\\ 
&&\hspace{.7in}
  \(\De^{ h}\De^{ -h}
\,p_{r'_{i} }(x )\) \,dx\bigg)\prod_{ l=0}^{ K}\prod_{q=1}^{n_{l}}\,dr_{l,2q}+O(h^{(2\bb-1)n+\ep})\nn\\
&&\qquad=\(c_{\psi,h,1}\)^{n}
  F (\wt\si,0,\ldots, 0)+O(h^{(2\bb-1)n+\ep}).\label{m5}
\end{eqnarray}
We now use  the notation introduced in the paragraph containing   (\ref{m.5kp}),  and the fact that there are $2^{2n}$ permutations that are compatible with $\PP$, to see that      \bea
      && 
 \sum_{\pi\sim\mathcal{P}} \int \mathcal{T}_{h}( x;\,\pi,e)\prod_{j,k,i}\,dx_{j,k,i}\label{m6a}\\
 &&\qquad=\(4c_{\psi,h,1}\)^{n}
  \sum_{\wt\si\in\wt\mathcal{B}}F (\wt\si,0,\ldots, 0)+O(h^{(2\bb-1)n+\ep}).\nn
  \eea
   Since   $|\wt \mathcal{B}|\leq (2n)!$,  we see that the error term only depends on $m$, (recall that  $m=2n$).
 Consider (\ref{m6a}) and the definition of $F (\wt\si,0,\ldots, 0)$  in (\ref{m.14}) and use
  (\ref{m.4}),  with $m_{j,k}$ replaced by $n_{j,k}$, to see that
        \bea
      && 
 \sum_{\pi\sim\mathcal{P}} \int \mathcal{T}_{h}( x;\,\pi,e)\prod_{j,k,i}\,dx_{j,k,i}\label{m6b}\\
 &&\qquad=\(4c_{\psi,h,1}\)^{n}
  E\(\prod_{\stackrel{j,k=0}{j< k}}^{K}\(\al_{j,k}\)^{
n_{j,k}}\) +O(h^{(2\bb-1)n+\ep}).\nn
  \eea
  
  Recall the definition of $\mathcal{S}$, to set of special pairings, given in the first paragraph of this subsection.
  Since   there are ${( 2n_{j,k})!\over 2^{ n_{j,k}}n_{j,k}!}$ pairings of the $2n_{j,k}$
elements $\{1,\ldots, m_{j,k}\}$, (recall that    $m_{j,k}=2n_{j,k}$), we see that when we sum over all the special pairings we get
 \begin{eqnarray} && 
\sum_{\mathcal{P}\in \mathcal{S}}\sum_{\pi\sim \mathcal{P}}\int \mathcal{T}_{h}( x;\,\pi,e)\prod_{j,k,i}\,dx_{j,k,i}
\label{m.25}\\ 
&&\quad
=  \prod_{\stackrel{j,k=0}{j< k}}^{K}{( 2n_{j,k})!\over 2^{ n_{j,k}}n_{j,k}!}\(4c_{\psi,h,1}\)^{ n_{j,k}} E\lc\prod_{\stackrel{j,k=0}{j< k}}^{K}\(\al_{j,k}\)^{
n_{j,k}}\rc  +O\(h^{(2\bb-1)n+\ep}\). 
\nn
\end{eqnarray}  
  It follows from (\ref{ord.1})   that the error term,   still,  only depends on $m$.

\medskip	The right-hand side of (\ref{m.25}) is precisely the desired expression in (\ref{m.2}). Therefore, to complete the proof  of Lemma \ref{lem-multiple}, we show that for all the other possible values of $a$, the integral in (\ref{4.25}) can be absorbed in the error term.

 \subsection{  $\bf a=e$ but not all cycles are   of order two or $\bf a\ne e$}  \label{sec-3.2}  
 
\begin{lemma} \label{3.31}
Suppose that $ a=e$ but not all cycles are   of order two or $  a\ne e$.
  Then 
\begin{equation}
 \int \mathcal{T}_{h}( x;\,\pi ,a )\prod_{j,k,i}\,dx_{j,k,i}=O\(\frac{h^{\ep}}{h\psi^2(1/h)}\)^n\label{m.8d1},
\end{equation}
 for some $\ep>0$.
\el
   
   In the rest of this  section we ignore all factors of  $\log 1/h$.

\medskip	
\Proof  Consider the basic formula  (\ref{m.9}).  Since we    only need an upper bound, we  take absolute values in the integrand and extend the time integral to $[0,1]$, as we have done several times above.  We take the time integral and get an upper bound for (\ref{m.9}) involving the terms  $u$, $v$ and $w$.
  Since ${  a\ne e}$, the number of $w$ terms is less than  $2n$.  
  
 We obtain (\ref{m.8d1}) by dividing the 
$u$, $v$ and $w$ terms in $ \mathcal{T}_{h}( x;\,\pi ,a )$ into sets.  Clearly, if a set contains $k$  terms of the form $w$   and  $k'$  terms of the form $v$,  there are $2k+k'$   difference operators $\De^{ h}_{\cd}$   associated with this set.  There are no difference operators associated with sets of $u$ terms.
  
Consider a set of  two $w$  terms   that   lies in a cycle of order two. There are four  difference operators $\De^{ h}_{\cd}$   associated with this set.  We show this set contributes a bound to (\ref{m.8d1})
that is 
\be
O \(1\over h\psi^{2}(1/h) \)  .
\ee
(By contributes a bound we mean that this is what we get after we make an appropriate change of variables and integrate out the $w$ terms in this set.) Thus we may say that each difference operator in a cycle of order two contributes a bound of 
\be
O\(\(1\over h\psi^{2}(1/h))\)^{1/4}\) .
\ee
 
We show that   any set that has   $k>0$ associated difference operators    except for    a set of  two $w$  terms   that  forms a cycle of order two  contributes a bound that is 
   \be
O\(\(1\over h\psi^{2}(1/h))\)^{k/4}\)h^{\ep},\label{sit}
\ee
 for some $\ep>0$. 
 
 There are $4n $  difference operators $\De^{ h}_{\cd}$, in $ \mathcal{T}_{h}( x;\,\pi ,a )$. 
 Consequently unless the graph associated with $ \mathcal{T}_{h}( x;\,\pi ,a )$ consists solely of cycles of order two, we obtain (\ref{m.8d1}). 

 As we construct the sets of  $u$, $v$ and $w$ terms, we  also choose a  collection $\II\cup\II' $ of   $m$ terms with arguments  that are linearly independent.   To bound the contribution of each set we bound all the terms not in $\II\cup\II' $ by their supremum, and, after changing variables, integrate the terms in $\II\cup\II' $. Using (\ref{bl.2}), (\ref{bl.01}), (\ref{bl.3}) and (\ref{bl.4}) we  verify the bounds given in the preceding paragraph. (Actually, there   is an exceptions to this rule which we also deal with.)

This is  how we divide the $u$, $v$ and $w$ terms into  sets.  
For each $\pi$ and $ a$ we define  a multigraph $G_{\pi ,a}$ with vertices   
$\{(j,k,i),\,0\leq j<k\leq K,\, 1\leq i\leq m_{j,k}\}$, and  an edge 
between the vertices $\pi_{l } (q-1)$ and $ \pi_{l}(q)$  whenever $ a(l,q)=(1,1)$. This graph divides the $w$ terms into cycles and chains.  Assume that there are 
  $S$ cycles.  We denote them by 
   $C_{s}=\{\phi_{s,1},\ldots, \phi_{s,l(s)}\}$, written in cyclic order,  where  the cycle length $l(s)=|C_{s} |\geq 1$ and $1\le s \le S$.  For each  $1\le s \le S$ we take the set of $l(s)$ terms
   \begin{equation}
  \mathcal{G}^{\mbox{\scriptsize cycle}} _{s}=\{w(x_{\phi_{s,2}}-x_{\phi_{s,1}}),\cdots,    w(x_{ \phi_{s,l(s)}}-x_{ \phi_{s,l(s)-1}}), w(x_{\phi_{s,1}}-x_{ \phi_{s,l(s)}}) \} \label{gr.1}.
   \end{equation}
Let 
\be
y_{\phi_{s,i}}=x_{\phi_{s,i}}-x_{\phi_{s,i-1}},\qquad i=2,\ldots,l(s).\label{4.71v}
\ee
  It is easy to see that $\{y_{\phi_{s,i}}\,|  i=2,\ldots, l(s)\}$, are linearly independent. We put the corresponding $w$ terms, $w(x_{\phi_{s,2}}-x_{\phi_{s,1}}),\cdots,    w(x_{ \phi_{s,l(s)}}-x_{ \phi_{s,l(s)-1}})$ into $\mathcal{I}$.
(On  the other hand, since
\be
\sum_{i=2}^{l(s)}y_{\phi_{s,i}}=-(x_{\phi_{s,1}}-x_{ \phi_{s,l(s)}}),\label{4.71s}
\ee   
we see that we can only extract $l(s)-1$ linearly independent variables from the $l(s) $ arguments of $w$ for a given $s$.)
  A cycle of length $1$ consists of a single point  $\phi_{s,1}=\phi_{l(s),1}$ in   the graph, so in this case
   \begin{equation}
  \mathcal{G}^{\mbox{\scriptsize cycle}} _{s}=\{w(0) \} \label{gr.1z}.
   \end{equation}
 We explain below how this can occur.  Obviously, $w(0)$ is not put into $\II$.

Next, suppose there are
  $S'$ chains. We denote them by 
   $C'_{s}=\{\phi'_{s,1},\ldots, \newline \phi'_{s,l'(s)}\}$, written in order,  where $l'(s)=|C'_{s} |\geq 2$ and $1\le s \le S'$.    Note that there are $l'(s)-1$, $w$ terms corresponding to  $C'_{s}$.  
Then for each  $1\le s \le S'$ we form the set of $l'(s)+1$ terms
   \bea
   &&
  \mathcal{G}^{\mbox{\scriptsize chain}}_{s}=\{v(x_{ \phi'_{s,1}}-x_{a(s)}), w(x_{\phi'_{s,2}}-x_{\phi'_{s,1}}),\cdots, \label{gr.2}\\
  &&\hspace{1 in}\cdots,    w(x_{ \phi'_{s,l(s)}}-x_{ \phi'_{s,l(s)-1}}), v(x_{b(s)}-x_{ \phi'_{s,l(s)}}) \}\nn
   \eea
   where $v(x_{ \phi'_{s,1}}-x_{a(s)})$ is the unique $v$ term associated with $\De^{x_{ \phi'_{s,1}}}_{h}$, and similarly, $v(x_{b(s)}-x_{ \phi'_{s,l(s)}}) $ is the unique $v$ term associated with $\De^{x_{ \phi'_{s,l(s)}}}_{h}$. 
  (This deserves further clarification. There may be other $v$ terms containing the variable  $x_{ \phi'_{s,1}}$. But there is only one $v$ term   of the form
 \begin{equation}
   \int_{0}^{1}\Big|\De^{x_{ \phi'_{s,1}}}_{h}p_{s}(x_{ \phi'_{s,1}}-u)\Big|\,ds
   \end{equation}
 where $u$ is some other $x_{\cd}$ variable which we denote by $x_{a(s)}$. This is because one operator $\De^{x_{ \phi'_{s,1}}}_{h}$     is associated with $w(x_{\phi'_{s,2}}-x_{\phi'_{s,1}})$ and there are precisely two operators $\De^{x_{ \phi'_{s,1}}}_{h}$ in (\ref{m.8d1})).
   
It is easy to see that variables $ y_{\phi'_{s,i}}=x_{\phi'_{s,i}}-x_{\phi'_{s,i-1}},\,  i=2,\ldots, l(s) $, are linearly independent. 
 We put the  $w$ terms, $w(x_{\phi'_{s,2}}-x_{\phi'_{s,1}}),\cdots,    w(x_{ \phi'_{s,l(s)}}-x_{ \phi'_{s,l(s)-1}})$ into $\mathcal{I}$.  We leave the $v$ terms in $\mathcal{G}^{\mbox{\scriptsize chain}}_{s}$ out of $\mathcal{I}$.

 \medskip	
 At this stage we emphasize that the terms we have put in $\mathcal{I}$ from all cycles and chains have linearly independent arguments. If fact, the set of $x$'s appearing in the different chains and the cycles are disjoint. This is obvious for the cycles and the interior of the chains 
 since there are exactly two difference operators  
 $\De^{x}_{h}$ for each $x$. It also must be true for the endpoints of the chains, since if this is not the case they could be made into larger chains or cycles.
 
For the same reason,   if a $v$  term involving  $\De_{ x'}^{h}$ is not in any of the  sets of chains, then  $x'$ will not appear in the arguments of the terms that  are put in $\mathcal{I}$ from all the cycles and chains.

Suppose, after considering the $w$ terms and the $v$ terms associated with the chains of $w$ terms, that there are $p$ pairs of   $v$ terms left, each pair corresponding to difference operators $\De^{h}_{z_{j}}$, $j=1,\ldots,p$. ($p$ may be $0$). Let
\be
\ZZ:=\{z_{1},\ldots,z_{p}\}
\ee
A typical $v$ term is of the form  
\begin{equation}
 v^{(j)}(z_{j}-u_{j'}):=v (z_{j}-u_{j'})=\int_{0}^{1}|\De^{h}_{z_{j} } p_{t}(z_{j}-u_{j'})|\,dt.\label{3.93qq}
   \end{equation}
where $u_{j'}$ is some  $x_{\cd}$ term. We use the superscript $(j)$ is  to keep track of the fact that this $v$ term is associated with the  difference operator $\De^{h}_{z_{j}}$. We distinguish between  the variables $z_{j} $ and $u_{j'}$ by referring to $z_{j}$ as a marked variable.   Note that if $u_{j'}$ is also in $\ZZ$, say $u_{j'}=z_{k}$,  then  $u_{j'}$ is also a marked variable but in a different $v$ term. (In this case, in $v^{(k)}(z_{k}-u_{k'})$, where $u_{k'}$ is some other $x_{\cd}$ variable.)

Thus $\ZZ$ is the collection of marked variables.
Consider the corresponding terms
\begin{equation}
  v^{(j)}(z_{j}-u_{j})\quad\mbox{and}\quad   v^{(j)}(z_{j}-v_{j}), \hspace{.2 in}j=1,\ldots,p\label{114.0}
   \end{equation}
  where $u_{j}$ and $v_{j}$ represent whatever terms $x_{\cd}$ and $x'_{\cd}$ are coupled with the two variables $z_{j}$. 
   
 There may be some $j$ for which $u_{j}$ and $v_{j}$ in (\ref{114.0}) are   both in $\ZZ$. Choose such a $j$.
  Suppose  $u_{j}=v_{j}=z_{k}$.  We set  
   \bea
   &&
  \mathcal{G}^{\mbox{\scriptsize $\ZZ,1$}}_{j}=\{  v^{(j)}(z_{j}-z_{k}),   v^{(j)}(z_{j}-z_{k}), \label{gr.3}\\
   &&\hspace{1 in}    v^{(k)}(z_{k}-u_{k}),  v^{(k)}(z_{k}-v_{k}) \}\nn
   \eea
  and put $v^{(j)}(z_{j}-z_{k})$ into $\mathcal{I}$.  Here $u_{k}$ and $v_{k}$ are whatever two variables appear with the two marked variables $z_{k}$.   
 
  On the other hand, suppose $u_{j}$ and $v_{j}$ are   both in $\ZZ$ but $u_{j}= z_{k}$ and $v_{j}= z_{l}$ with $k\neq l$. We set 
 \bea
   &&
  \mathcal{G}^{\mbox{\scriptsize $\ZZ,2$}}_{j}=\{  v^{(j)}(z_{j}-z_{k}),   v^{(j)}(z_{j}-z_{l}), \label{gr.4}\\
  &&\hspace{1 in}    v^{(k)}(z_{k}-u_{k}),  v^{(k)}(z_{k}-v_{k}),v^{(l)}(z_{l}-u_{l}),  v^{(l)}(z_{l}-v_{l}) \}\nn
   \eea
  and put both  $v^{(j)}(z_{j}-z_{k})$ and $v^{(j)}(z_{j}-z_{l})$ into $\mathcal{I}$.

  We then turn to the  elements in $\ZZ$  which have not yet appeared in the arguments of the terms that have  been put into  $\mathcal{I}$.  If there is another $j'$ for which $u_{j'}$ and $v_{j'}$ are   both in $\ZZ$, choose such a $j'$ and proceed as above. If there are no longer any such elements in $\ZZ$, choose some remaining element, say, $z_{i}$.   
 Set 
 \bea
   &&
  \mathcal{G}^{\mbox{\scriptsize $\ZZ,3$}}_{i}=\{  v^{(i)}(z_{i}-u_{i}),   v^{(i)}(z_{i}-v_{i})\} \label{gr.5} 
   \eea
   and if $u_{i}\not\in \ZZ$, place $ v^{(i)}(z_{i}-u_{i})$ into $\mathcal{I}$. If $u_{i} \in \ZZ$, so that 
 $v_{i}\not\in \ZZ$,   place $ v^{(i)}(z_{i}-v_{i})$ into $\mathcal{I}$.    
    
We then continue until we have exhausted $\ZZ$.  We  form a final set $ \mathcal{G}^{u}$ which contains all the $u$ terms, so that  all $u$, $v$ and $w$ terms have been divided into sets. 

It is possible that there are no cycles of  length one.  We show how we get (\ref{m.8d1}) in this case.

We have constructed  $\mathcal{I}$ so that all its members have linearly independent arguments. However, $\mathcal{I}$ may contain less than $m $ terms. We simply add to $\mathcal{I}$ a set $\mathcal{I'}$ of enough of the remaining $u$ and $v$ terms so that $\mathcal{I}\cup \mathcal{I'}$  has $m $ terms, whose arguments span  $R^{2n} $,  the space spanned by the original $x_{\cd}$ terms. (It follows from (\ref{4.71s})   that no further $w$ terms can be added to $\mathcal{I'}$). We  bound the $v$ terms  in $\mathcal{I'}$ as follows:  
\begin{equation}
   |v(x'-x'')|\le   \frac{C}{h\psi(1/h)(1+(x'-x'')^{2})}.  \label{3.97} 
   \end{equation}
We then make a change of variables setting the arguments of the terms in $\mathcal{I}\cup \mathcal{I'}$ equal to $y_{1},\ldots,y_{m}$ and bound the  $v$ terms not in $\mathcal{I}\cup \mathcal{I'}$ by $C({h\psi(1/h)})^{-1}$ and the $u$ terms not in $\mathcal{I}\cup \mathcal{I'}$ by $ C$. Finally we   integrate. We have $m$ one dimensional integrals which we bound by   (\ref{bl.01}) for the $v$ terms in $\II$,   by $C({h\psi(1/h)})^{-1}$  for the $v$ terms in $\II'$, and by (\ref{bl.4}) for $w$ terms in $\II$. The integrals of the $u$ terms in $\II$ we bound by a constant; (see (\ref{2.2})).

  Clearly  $ \mathcal{G}^{u}$ gives a bounded contribution. We now show that (\ref{sit}) holds for all other sets of $v$ and $w$ terms, with the exception of sets of $w$ terms in cycles of length $2$.

Consider first  $ \mathcal{G}^{\mbox{\scriptsize cycle}} _{s}$ for a 
  cycle of lengths $l(s)$. We integrate the $l(s)-1$,  $w$ terms which were put in $\mathcal{I}$ and  bound the remaining $w$ term by $C({h\psi(1/h)})^{-1}$ to obtain the bound 
\be
  C  \(\frac{1}{\psi(1/h)}\)^{l(s)-1}\frac{1}{h\psi(1/h)}  = C  \(\frac{1}{\psi(1/h)}\)^{l(s)-2}\frac{1}{h\psi^{2}(1/h)}. \label{nm.1}
   \ee
   Since
\begin{equation}
 \frac{1}{\psi(1/h)} =h^{1/2} \(\frac{1}{h\psi^{2}(1/h)}\)^{1/2}\label{110}
   \end{equation}
(\ref{nm.1}) is bounded by 
\begin{equation}
C  \Big\{h^{(l(s)-2)/2}\Big\}  \(\frac{1}{h\psi^{2}(1/h)}\)^{l(s)/2}. \label{110j}
\end{equation}
Since a cycle of length $l(s)$ involves $2l(s)$  difference operators  $\De_{h}$,   and $l(s)/2=2l(s)/4$, we are in the situation of (\ref{sit}), unless all cycles are of order two.  (This shows, incidentally, that when $a=e$, (\ref{m.8d1}) holds unless all cycles are of order two.)

Consider next  $ \mathcal{G}^{\mbox{\scriptsize chain}} _{s}$. Recall that there are   $l'(s)-1$, $w$ terms in a chain, where   $l'(s)\geq 2$. We have put all $l'(s)-1$  terms  $w $ in   $\mathcal{I}$, and we     can bound their integrals  by 
\begin{equation}
   C \(\frac{1}{\psi(1/h)}\)^{l'(s)-1}.\label{nm.2}
   \end{equation} 
    In addition there are two $v$ terms in $ \mathcal{G}^{\mbox{\scriptsize chain}} _{s}$. 
 The ones not in $\II'$ can be bounded by $C({h\psi(1/h)})^{-1}$ and the ones in $\II'$ are bounded by (\ref{3.97}), which after integration also contributes $C({h\psi(1/h)})^{-1}$. Thus we obtain the following bound for    for  $ \mathcal{G}^{\mbox{\scriptsize chain}} _{s}$:
\begin{eqnarray}
&&C  \(\frac{1}{\psi(1/h)}\)^{l'(s)-1}\(\frac{1}{h\psi(1/h)}\)^{2} 
\label{nm.3}\\
&&\qquad= C \(\frac{1}{\psi(1/h)}\)^{l'(s)-3}\(\frac{1}{h\psi^{2}(1/h)}\)^{2} \nonumber\\
&&\qquad\leq   C h^{(l'(s)-3)/2}\(\frac{1}{h\psi^{2}(1/h)}\)^{1/2}  \(\frac{1}{h\psi^{2}(1/h)}\)^{l'(s) /2}. \nonumber
\end{eqnarray}
Note that each chain of length  $l'(s)$ together with the two $v$ terms associated with the end points involves $2l'(s) $  difference operators  $\De_{h}$. Clearly if $l'(s)\geq 3$ we are in the situation of (\ref{sit}). This holds even for chains of length $l'(s)=2$ since 
\begin{equation}
h^{(2-3)/2}\(\frac{1}{h\psi^{2}(1/h)}\)^{1/2}=\frac{1}{h\psi(1/h)}.\label{nm.4}
\end{equation}

Note that the $v$ terms that were not initially in $\II$ contribute a bound of 
$C({h\psi(1/h)})^{-1}$, whether or not they are placed in $\II'$. We continue to  use this fact below without commenting on it further.

We next consider $\mathcal{G}^{\mbox{\scriptsize $\ZZ,1$}}_{j}$. We   integrate the one $v$ term in $\II$   and any that are in $\II'$ and bound the remaining ones. This gives a bound of 
   \begin{equation}
 \frac{1}{h^{2}\psi^{3}(1/h)} =\frac{1}{h\psi (1/h)} \(\frac{1}{h\psi^{2}(1/h)}\). \label{110s}
   \end{equation}
Since $\mathcal{G}^{\mbox{\scriptsize $\ZZ,1$}}_{j}$ involves  four $\De^{h}_{\cd} $ operators   we are in the situation of (\ref{sit}).

For $\mathcal{G}^{\mbox{\scriptsize $\ZZ,2$}}_{j}$  we integrate two $v$ terms in $\II$ and   any that are in $\II'$ and bound the remaining ones. This gives a bound of 
   \begin{equation}
\( \frac{1}{h \psi^{2}(1/h)}\)^{2} =\(\frac{1}{h\psi^{2} (1/h)}\)^{1/2} \(\frac{1}{h\psi^{2}(1/h)}\)^{3/2}; \label{110r}
   \end{equation}
Since $\mathcal{G}^{\mbox{\scriptsize $\ZZ,2$}}_{j}$ involves  six $\De^{h}_{\cd} $ operators   we are in the situation of (\ref{sit}).

Finally, for $\mathcal{G}^{\mbox{\scriptsize $\ZZ,3$}}_{j}$  we integrate the one $v$ term in $\II$ and   the other if it is in $\II'$. Otherwise we bound it.   This gives a bound of 
   \begin{equation}
  \frac{1}{ \psi (1/h)}  =h^{1/2} \(\frac{1}{h\psi^{2}(1/h)}\)^{1/2}; \label{110rqq}
   \end{equation}
Since $\mathcal{G}^{\mbox{\scriptsize $\ZZ,3$}}_{j}$ involves  two $\De^{h}_{\cd} $ operators   we are in the situation of (\ref{sit}). 

This shows that if $a$ and the partition $\pi$ does not generate exclusively $w$ terms in cycles of order two and are such that there are no cycles of length  one, then (\ref{m.8d1}) holds.

\medskip	
We now remove the restriction that $a$ and   $\pi$ does not give rise to cycles of length  one. The only way  this anomaly  can occur is in  terms of the type   
\begin{equation}
\De^{ h}\De^{- h}p_{(1-\sum_{q=1}^{m_{l-1}}r_{l-1,q})+r_{l,1}}(x_{\ga_{2\si_{l} (1)-1}}-x_{\ga_{2\si_{l-1} (n_{l-1})-1}})\label{ex.1}
\end{equation}
when $\ga_{2\si_{l} (1)-1}=\ga_{2\si_{l-1} (n_{l-1})-1}$.    Note that in this case
\bea
   \int_{0}^{t}\De^{ h}\De^{- h}p_{s }(x_{\ga_{2\si_{l} (1)-1}}-x_{\ga_{2\si_{l-1} (n_{l-1})-1}})\,ds&=&w(0). \label{3.110}
   \eea
This is what we call a cycle of length one. In this case we have
 \begin{equation}
\De^{ h}\De^{- h}p_{(1-\sum_{q=1}^{m_{l-1}}r_{l-1,q})+r_{l,1}}(0)=-2\De^{ h}p_{(1-\sum_{q=1}^{m_{l-1}}r_{l-1,q})+r_{l,1}}(0).\label{ex.2}
\end{equation}

 We  now show how to deal with (\ref{ex.2}).  We return to the  basic formulas   (\ref{m.8}) and (\ref{m.9}). We obtain an upper bound for (\ref{m.9})  by taking the absolute value of the integrand. However, we do not, initially extend  the region of  integration with respect to time.  Instead  we proceed as follows: Let  $l'$ be the largest value of $l$   for which (\ref{ex.2}) occurs.   We extend the integral with respect to $r_{l,q}$  for all $l>l'$, and also  for $l=l'$ and $q>1$, and  bound these integrals with terms of the form $u$, $v$ and $w$. We then consider  the integral of the term in (\ref{ex.2}) with respect to $r_{l',1}$.

Clearly
\be
  \int _{0}^{1}|\Delta^{h}p_{(1-\sum_{q=1}^{m_{l'-1}}r_{l'-1,q})+r_{l',1}}(0)| \,dr_{l',1}\label{ex.48} 
  \leq    \int _{1-\sum_{q=1}^{m_{l'-1}}r_{l'-1,q}}^{2-\sum_{q=1}^{m_{l'-1}}r_{l'-1,q}}|\Delta^{h}p_{ s}(0)|   \,ds 
\ee
If  $\sum_{q=1}^{m_{l'-1}}r_{l'-1,q}\leq 1/2$ this last integral 
\bea
  && \le   \int _{1/2}^{2 }|\Delta^{h}p_{ s}(0)| \,ds \leq Ch^{2} \label{3.93}
\eea
by (\ref{4.4}).  Since we have only used two $\De^{h}_{\cd}$ operators  we are in the situation of (\ref{sit}).
 
If  $\sum_{q=1}^{m_{l'-1}}r_{l'-1,q}\geq 1/2$ then for some $q'$ we have $r_{l'-1,q'}\geq 1/2m$. 
  Note that the variable $r_{l'-1,q'}$  
appears in (\ref{ex.1}) and in only one other term.  If $q'>1$,  then using the fact that $r_{l'-1,q'}\geq 1/2m$,   we use one of  the bounds in     Lemma \ref{lem-4.3},  to bound a term which in the non-exceptional case would be  $u$, $v$ or $w$, or their integrals with respect to $x$, by 
$\ov u_{1/2m}$, $ \ov v_{1/2m}$ or $\ov w_{1/2m}$, or their integrals with respect to $x$.  One sees from Lemma  \ref{lem-vproprvt} that  
we don't loose anything in comparison with the non-exceptional case. The case
$r_{l'-1,1}\geq 1/2m$ and  $\ga_{2\si_{l'-1} (1)-1}\neq \ga_{2\si_{l'-2} (n_{l'-2})-1}$ is handled the same way.

On the other hand if $r_{l'-1,1}\geq 1/2m$
  and $\ga_{2\si_{l'-1} (1)-1}= \ga_{2\si_{l'-2} (n_{l'-2})-1}$, we use Lemma \ref{lem-4.2} to get  the same bound of $Ch^{2}$.

 After completing the procedure  described in the previous two paragraphs
we integrate in (\ref{m.9}) with respect to $r_{l'-1,q'}$ and $r_{l',1}$, since these variables now  appear only in the term in (\ref{ex.1}). What we are left with is bounded by
\begin{equation}
   \int_{1/2m}^{1-\sum_{q\neq q'}  r_{l'-1,q}}\int _{0}^{1}|\Delta^{h}p_{(1-\sum_{q=1}^{m_{l'-1}}r_{l'-1,q})+r_{l',1}}(0)| \,dr_{l',1}\,dr_{l'-1,q'}\label{3.94}
   \end{equation}
  Let $\al=1-\sum_{q\neq q'}  r_{l'-1,q}$. We make the change of variables $r=r_{l',1}$ and $s=-r_{l'-1,q'}+\al$ to get that (\ref{3.94})
\begin{eqnarray}
&&\nonumber\\
&&
 \leq   \int_{0}^{1 } \int _{0}^{1}|\Delta^{h}p_{r+s}(0)|\,dr \,ds\label{exblm.3x}w\\
&&  \leq  \int _{0}^{2}r |\Delta^{h}p_{r }(0)|\,dr\nn   \le C\int_{0}^{2} r \(\int \sin^{2}(ph)e^{-r\psi (p)}\,dp\) \,dr \\
&&   \leq Ch^{\bb} \int_{0}^{2} r \(\int p^{\bb}e^{-r\psi (p)}\,dp\) \,dr \nn\\
&&  \leq Ch^{\bb}  \int {p^{\bb} \over 1+\psi^{2} (p)} \,dp =O(h^{\bb}).\nn
\end{eqnarray}
Since
\begin{equation}
   h^{\bb}=h^{\bb+1/2}\psi(1/h)\(\frac{1}{h\psi^{2}(1/h)}\)^{1/2}
   \end{equation} 
 we are once again in the situation of (\ref{sit}).

We then apply a similar procedure for each $l$ in decreasing order, skipping those for which
   (\ref{ex.2}) occurs, if they were already bounded by the procedure described in the paragraph preceding the one containing (\ref{3.94}).    
 Thus we see  that  cycles of length one  are  in the situation of (\ref{sit}). We proceed to deal with remaining terms as we did when we assumed that there were no cycles of length one and see that (\ref{m.8d1}) holds. This completes the proof of Lemma \ref{3.31}.\qed

It follows from (\ref{m.25}) and Lemma \ref{3.31} that  when $m$ is even
 \begin{eqnarray} && 
\sum_{a}\sum_{\pi_{0},\ldots, \pi_{K}}\int \mathcal{T}_{h}( x;\,\pi,a)\prod_{j,k,i}\,dx_{j,k,i}
\label{m.25h}\\ 
&&\quad
=  \prod_{\stackrel{j,k=0}{j< k}}^{K}{( 2n_{j,k})!\over 2^{ n_{j,k}}n_{j,k}!}\(4c_{\psi,h,1}\)^{ n_{j,k}} E\lc\prod_{\stackrel{j,k=0}{j< k}}^{K}\(\al_{j,k}\)^{
n_{j,k}}\rc  +O\(h^{(2\bb-1)n+\ep}\). 
\nn
\end{eqnarray}  
  	 We  now show  that we get the same estimates when  $ \mathcal{T}_{h}( x;\,\pi ,a )$   is replaced by $ \mathcal{T}'_{h}( x;\,\pi ,a )$; (see   (\ref{m.7}) and (\ref{m.8})).

We point out, in the paragraph containing (\ref{3.27}) that terms of the form  $\De^{h}\De^{-h}p_{\cdot}^{\sharp}$ in (\ref{m.7}) are always of the form $\De^{h}\De^{-h}p_{\cdot}$. Therefore, in showing that (\ref{m.7}) and (\ref{m.8}) have the same asymptotic behavior as $h\to 0$ we need only consider how the proof of (\ref{m.25h}) must be modified when the arguments of the density functions with one or no difference operators applied is effected by adding $\pm h$.

 \label{page24}It is easy to see that the presence of these terms   has no effect on the  integrals that are $O\(h^{(2\bb-1)n+\ep}\)$ as $h\to 0$.    This is because in evaluating these expressions we either integrate over all of $R^{1}$ or else use bounds that hold on all of  $R^{1}$.  Since terms with one difference operator only occur in these estimations, we no longer need to be concerned with them.

   Consider the terms with no difference operators applied to them, now denoted by    $p^{\sharp}$. So, for example,  instead of $F (\wt\si, 0,\ldots,0) $ on  the right-hand side of (\ref{m.15sw}), we now have
 \begin{eqnarray} 
&& \int\(\int_{\wt\RR_{0}(0)\times\cdots \times\wt\RR_{K}(0)} \prod_{ l=0}^{ K}\, p^{\sharp}_{(1-\sum_{q=1}^{n_{l-1}}r_{l-1,2q-1}-s_{l-1})+r_{l,1}}\right.
\label{m.14k}\\
&&\left.
\hspace{.2 in}(x_{\wt  \si_{l} (1) }-x_{\wt  \si_{l-1} (n_{l-1}) }) \prod_{q=2}^{n_{l}} p^{\sharp}_{r_{l,2q-1}}(x_{\wt  \si_{l} (q) }-x_{\wt  \si_{l} (q-1) })
\prod_{q=1}^{n_{l}}\,dr_{l,2q-1}\)  \,d x    \nonumber.
\end{eqnarray} 
  Suppose that   $p^{\sharp}_{r}(y_{\si (i)}-y_{\si (i-1)})=p_{r} (y_{\si (i)}-y_{\si (i-1)}\pm h)$.  We write this term as
\begin{equation}\qquad
p^{\sharp}_{r}(y_{\si (i)}-y_{\si (i-1)})=p_{r} (y_{\si (i)}-y_{\si (i-1)} ) +\De^{\pm  h}p_{r}(y_{\si (i)}-y_{\si (i-1)}). \label{4.52}
   \end{equation}
   Substituting all such terms  into  (\ref{m.14k}) and expanding we get (\ref{m.25h}) and many other terms with at least one $p_{r} (y_{\si (i)}-y_{\si (i-1)} )$ replaced by $\De^{\pm  h}p_{r}(y_{\si (i)}-y_{\si (i-1)})$. In this case simply take these terms, extend their integrals to $[0,1]$ and bound them as in (\ref{bl.2}). Then follow the procedure in the paragraph containing (\ref{4.60}) to deal with the remaining terms and the functions $1/(1+ (y_{\si (i)}-y_{\si (i-1)} )^{2})$. In this the integral in (\ref{m.14k}) is bounded by $C(1/(h\psi(1/h)))^{j}$, where $j$ is the number of terms that have the difference operator applied. 
   Thus we see that replacing   $ \mathcal{T}_{h}( x;\,\pi ,a )$       by  $ \mathcal{T}'_{h}( x;\,\pi ,a )$
   does not change (\ref{m.25h}) when $m$ is even.
  
     \medskip   When $m$ is odd we can not construct a graph with all cycles of order 2. Therefore, we are not in the situation covered by Section \ref{assinged}. Moreover, in Section \ref{sec-3.2}  we never use the fact that $m$ is even. We actually obtain (\ref{m.8d1}) with $n$ replaced by $m/2$, which is what we assert in (\ref{m.2a}).  This also holds when when $p_{\cdot}$     is replaced  by $p_{\cdot}^{\sharp}$
for the reasons given in the preceding two paragraphs. \qed

 \subsection{Proof of Lemma \ref{lem-m4}}\label{sec-m4}
 
For any $A\subseteq [0,3]^{n}$ we set 
\begin{eqnarray}
&&F_{A}
 =\int\lc \int_{A} \prod_{ l=0}^{ K}\, p_{r_{l,1}}(x_{\wt  \si_{l} (1) }-x_{\wt  \si_{l-1} (n_{l-1}) }) \right.
\label{vm4.1}\\
&&
\left. \hspace{1 in}\prod_{q=2}^{n_{l}} p_{r_{l,2q-1}}(x_{\wt  \si_{l} (q) }-x_{\wt  \si_{l} (q-1) })
\prod_{ l=0}^{ K}\prod_{q=1}^{n_{l}}\,dr_{l,2q-1}\rc\prod_{q=1}^{n_{l}} \,dx_{\wt  \si_{l} (q) }.    \nn
\end{eqnarray} 
Then by H\"{o}lder's inequality, for any $1/a+1/b=1$
\begin{eqnarray}
&&\lc \int_{A} \prod_{ l=0}^{ K}\, p_{r_{l,1}}(x_{\wt  \si_{l} (1) }-x_{\wt  \si_{l-1} (n_{l-1}) }) \right.
\nn\\
&&
\left. \hspace{1 in}\prod_{q=2}^{n_{l}} p_{r_{l,2q-1}}(x_{\wt  \si_{l} (q) }-x_{\wt  \si_{l} (q-1) })
\prod_{ l=0}^{ K}\prod_{q=1}^{n_{l}}\,dr_{l,2q-1}\rc
\label{vm4.2}\\
&&\leq |A|^{1/a}\lc \int_{[0,3]^{n}} \prod_{ l=0}^{ K}\, p^{b}_{r_{l,1}}(x_{\wt  \si_{l} (1) }-x_{\wt  \si_{l-1} (n_{l-1}) }) \right.
\nn\\
&&
\left. \hspace{1 in}\prod_{q=2}^{n_{l}} p^{b}_{r_{l,2q-1}}(x_{\wt  \si_{l} (q) }-x_{\wt  \si_{l} (q-1) })
\prod_{ l=0}^{ K}\prod_{q=1}^{n_{l}}\,dr_{l,2q-1}\rc^{1/b},   \nonumber
\end{eqnarray} 
where $|A| $ denotes the volume of $A$ in $R_+^n$.

Since $\bb>1$ we can choose a $1<b<\bb$ such that  
\begin{equation}
   \int_{0}^3\({\psi^{-1}(1/s)}\)^{b}\,ds\le C.
   \end{equation}
Therefore it follows from  (\ref{2.dens}) that 
\begin{equation}
\int^{3}_{0}p^{b}_{r}(x)\,dr\leq C\frac{1}{1+x^{2}}.\label{vm4.3}
\end{equation}
Thus there exists a finite constant  $C$, depending only on $n$ and $b$, that is  independent of
 $A $, such that
\begin{equation}
F_{A}\leq C|A|^{1/a}.\label{vm4.4}
\end{equation}

It follows from   (\ref{m.14}), paying special attention to the time variable of    $p_{\cd}$ in the second line,  that 
\begin{equation}
F (\si,s_{0},\ldots, s_{K})=F_{A_{s_{0},\ldots, s_{K}}}\label{vm4.5}
\end{equation}
where
\begin{eqnarray}
 A_{s_{0},\ldots, s_{K}}&=&\lc r\in R_{+}^{n}\,\Bigg| \sum_{\la=0}^{l-1}(1-\sum_{q=1}^{n_{\la}}r_{\la,2q-1}-s_{\la})\leq \sum_{q=1}^{n_{l}}\, r_{l,2q-1}\right.    \label{vm4.6}\\
&&   \left. \hspace{.2 in}\leq \sum_{\la=0}^{l-1}(1-\sum_{q=1}^{n_{\la}}r_{\la,2q-1}-s_{\la})+(1-s_{l});\,l=0,1,\ldots,K\rc.
\nn
\end{eqnarray}
  In particular  
\begin{eqnarray}
 A_{0,\ldots, 0} &=&\lc r\in [0,3]^{n}\,\Bigg| \sum_{\la=0}^{l-1}(1-\sum_{q=1}^{n_{\la}}r_{\la,2q-1})\leq \sum_{q=1}^{n_{l}}\, r_{l,2q-1}\right.    \label{vm4.6a}\\
&&   \left. \hspace{.2in}\leq \sum_{\la=0}^{l-1}(1-\sum_{q=1}^{n_{\la}}r_{\la,2q-1})+1);\,l=0,1,\ldots,K\rc.
\nn
\end{eqnarray}
Let  $\phi_{l}(r)=\sum_{\la=0}^{l}(1-\sum_{q=1}^{n_{\la}}r_{\la,2q-1})$. We have
\begin{eqnarray}
\lefteqn{ A_{s_{0},\ldots, s_{K}}\De\,A_{0,\ldots, 0}
\label{4.124}}\\
&& \subseteq \bigcup_{l=1}^{K}   \lc r\in [0,3]^{n}\,\Bigg| \phi_{l-1}(r)-\sum_{\la=0}^{l-1} s_{\la}\leq \sum_{q=1}^{n_{l}}\, r_{l,2q-1}\leq  \phi_{l-1}(r)\rc \nonumber\\
&&\qquad\bigcup_{l=0}^{K}  \lc r\in [0,3]^{n}\,\Bigg| \phi_{l-1}(r)+1-\sum_{\la=0}^{l} s_{\la}\leq \sum_{q=1}^{n_{l}}\, r_{l,2q-1}\leq  \phi_{l-1}(r)+1\rc. \nonumber
\end{eqnarray}
(Note that the first union are the points in $A_{s_{0},\ldots, s_{K}}$ that are not in $A_{0,\ldots, 0}$ and the second union are the points in  $A_{0,\ldots, 0}$ that are not in $A_{s_{0},\ldots, s_{K}}$.)

Since for fixed $a\geq b\geq 0$  
\begin{equation}
\Bigg | \lc r\in [0,3]^{n_{l}}\,\Bigg| a-b\leq \sum_{q=1}^{n_{l}}\, r_{l,2q-1}\leq  a\rc \Bigg |\leq 
Cb^{n_{l}}\label{4.124a}
\end{equation}
we have that  
\bea
| A_{s_{0},\ldots, s_{K}}\De\,A_{0,\ldots, 0}  |&\leq & C K  \(\sum_{\la=0}^{K} s_{\la}\)^{2n}  \label{vm4.7}\\
&\leq &  C K^{m+1}(\max_{0\le \la\le K} s_{\la})^{m}\nn,
\eea
  when $\max_{0\le \la\le K} s_{\la}$ is sufficiently small. 
  Let $K'$ be the cardinality of $\{l\,|\, n_{l}> 0\}$. It is easy to see that we have actually   proved (\ref{vm4.7})  
 with $K$ replaced by $K'$.
Since  $K' \le m$,   the last line in (\ref{vm4.7}) can be written in terms of $\{s_{\la}\}$ and $m$.
  Lemma \ref{lem-m4} follows from  (\ref{vm4.4}) and (\ref{vm4.7}).\qed

 \section{Proof of Theorem \ref{theo-clt2r}}\label{sec-CLT}
 
 The proof of  Theorem \ref{theo-clt2r} follows easily from the preliminary material  in Lemmas \ref{lem-tri}--\ref{lem-4.5}.
 
Let  
\begin{equation}
I_{j,k,l,h}:=\int ( L^{ x+h}_{1/l}- L^{ x}_{1/l})\circ\th_{j/l}\,\, (   L^{ x+h}_{1/l}-  L^{ x}_{ 1/l})\circ  \th_{k/l} \,\, \,dx\label{e1.1}
\end{equation}
and 
\begin{equation}
\wt I_{l,h} :=
 \sum_{\stackrel{j,k=0}{j< k}}^{l-1} I_{j,k,l,h}.\label{e1.6}
\end{equation}
Using the additivity property of local time we can write
\begin{equation}
L^{ x}_{ 1} =\sum_{j=0}^{l-1}L^{ x}_{1/l}\circ \th_{j/l}\label{e1.3}.
\end{equation}
so that 
\begin{equation}
 \int ( L^{ x+h}_{1}- L^{ x}_{ 1})\,\, (   L^{ x+h}_{1}-  L^{ x}_{ 1}) \, \,dx =
 \sum_{j,k=0}^{l-1} I_{j,k,l,h}=2\wt I_{l,h}+\sum_{ j=0 }^{l-1}I_{j,j,l,h}\label{e1.4}
\end{equation}
Similarly, let
\begin{equation}
  \al_{j,k,l }:= \int  L^{ x}_{1/l }\circ \th_{j/l}\,\,  L^{ x}_{1/l}\circ  \th_{k/l}\,\,  \,dx.\label{e1.2}
\end{equation}
 and
\begin{equation}
\wt \al_{l } := \sum_{\stackrel{j,k=0}{j< k}}^{l-1}\al_{j,k,l}
,\label{e1.7}
\end{equation}
Thus
\begin{equation}
 \al_{1}=\int  L^{ x}_{1} \,\,  L^{ x}_{1} \,  \,dx=  \sum_{j,k=0}^{l-1}\al_{j,k,l}=2\wt \al_{l }+\sum_{ j=0 }^{l-1}\al_{j,j,l}.\label{e1.5}
\end{equation}

The main ingredient in the proof of  Theorem \ref{theo-clt2r} is Lemma \ref{lem-tri}  below. We use the following lemma in its proof.

  \begin{lemma} \label{lem-alp} Under the hypotheses of Theorem \ref{theo-clt2r},
  for $\al_{t}$ as defined in (\ref{1.4}),  there exists an $\ep>0$ and a $t_{0}:=t_{0}(\ep)$ such that for all integers $n>0$  
   \begin{equation}
 \| \ \al_{t}\|_{n}\leq C_{n}t^{2}\psi^{-1}(1/t)  \label{e1.59z},
  \end{equation}
for all $0<t\le t_{0}$ and some constant $C_{n}$ depending only on $n$.   (Note that  $\al_{t}$ is given in (\ref{1.4}) and   $\|\cd\|_{n}:=(E(\cd)^{n})^{1/n}$.) 

  In addition 
\begin{equation}
   \lim_{l\to\ff}E (\wt\al_{l})^{n}=E ( \al_{1}/2)^{n}.\label{4.9aq}
   \end{equation} 

 \end{lemma}

\Proof  
  By Kac's moment formula
  \begin{eqnarray}
 E\lc\( \al_{t}\)^{
n}\rc&=&E\(\(\int ( L^{ x}_{ t})^{ 2}\,dx\)^{n}\)
\label{e1.60}\\
&=&  2^{n}\sum_{\pi }  \int \int_{\{\sum_{i=1}^{2n}r_{i}\leq t\}}\prod_{i=1}^{2n}
p_{r_{i}}(x_{\pi(i)}-x_{\pi(i-1)})\,\,    \prod_{i=1}^{2n}\,dr_{i}
\,\,    \prod_{i=1}^{n}\,dx_{i} \nonumber,
\end{eqnarray}
where the sum runs over all maps $\pi:\,[1,2n]\mapsto [1,n]$ with $|\pi^{-1}(i)|=2$ for each $i$. The factor $2^{n}$ comes from the fact that we can interchange each  $x_{\pi(i)}$ and $x_{\pi(i-1)}$, $i=1,\ldots,2n$. 
  
 It is not difficult to see that we can find a subset  $J=\{i_{1},\ldots, i_{n}\}\subseteq [1,2n]$, such that each of  $x_{1},\ldots, x_{n}$ can be written as a linear combination of    $y_{j}:=x_{\pi(i_{j})}-x_{\pi(i_{j}-1)}$, $j=1,\ldots ,n$. For  $i\in J^{c}$ we use the bound  $p_{r_{i}}(x_{\pi(i)}-x_{\pi(i-1)})\leq p_{r_{i}}(0)$, then change variables and integrate  out the $y_{j}$, to see that
\begin{eqnarray}
&& \int \(\prod_{i=1}^{2n}
\int_{0}^{t}p_{r_{i}}(x_{\pi(i)}-x_{\pi(i-1)})\,dr_{i}\)
\,\,    \prod_{i=1}^{n}\,dx_{i}
\label{e1.61a}\\
&& \qquad\leq  C\(\int_{0}^{t}p_{r}(0)\,dr \)^{n} \int \(\prod_{i\in J} 
\int_{0}^{t}p_{r_{i}}(x_{\pi(i)}-x_{\pi(i-1)})\,dr_{i}\)
\,\,    \prod_{i=1}^{n}\,dx_{i}\nonumber\\
&&\qquad \leq  C\(\int_{0}^{t}p_{r}(0)\,dr \)^{n} \(\prod_{i\in J} 
\int \int_{0}^{t}p_{r_{i}}(y_{i})\,dr_{i}\,dy_{i}\)\nn\\
&&\qquad=Ct^{n} \(\int_{0}^{t}p_{r}(0)\,dr \)^{n}\le C\(t^{2}\psi^{-1}(1/t) \)^{n},
\nonumber
\end{eqnarray}
 where we use  (\ref{bl.1}) and (\ref{pkacv.19}) for the last line. This gives (\ref{e1.59z}). 
 
 To obtain (\ref{4.9aq}) we use (\ref{e1.59z}) to see that 
 for  
$l$ sufficiently  large,
  \begin{eqnarray}
  \Big | \| 2\wt\al_{l}\|_{n}-\| \al_{1}\|_{n}\Big |&\leq &\| 2\wt\al_{l}- \al_{1}\|_{n}=  \| \sum_{ j=0 }^{l-1}\al_{j,j,l}\|_{n}\label{e1.57}\\
  &\le&  l\| \ \al_{0,0,l}\|_{n}
  =l\| \ \al_{1/l}\|_{n}\nonumber\\
  &\le&  C \frac{\psi^{-1}(l)}{l }.\nn 
  \end{eqnarray} 
\qed

\begin{lemma}\label{lem-tri}  Under the hypotheses of Theorem \ref{theo-clt2r} and with $l=l(h)=[\log 1/h]$,
\begin{eqnarray}
&& \lim_{h\rar 0} 2\sqrt{h\psi^{2}(1/h)}\wt I_{l,h}  
\label{e1.8}\stackrel{\mathcal{L}}{\Longrightarrow}\(8c_{\bb,1}\)^{ 1/2}\sqrt{\al_{1}}\,\,\eta.
\end{eqnarray}
\end{lemma}

\Proof   We show that  for each $m$ 
\begin{eqnarray} &&
\lim_{ h\rar 0}E\(\( 2\sqrt{h\psi^{2}(1/h)}\wt I_{l,h}  \)^{m}\)\nn\\ 
&&\hspace{ 1in}  =\left\{\begin{array}{ll} \displaystyle
{( 2n)!\over 2^{ n}n!}\( 8c_{\bb,1}\)^{ n} E\lc\(\al_{ 1}\)^{
n}\rc &\mbox{ if }m=2n\\\\
0&\mbox{ otherwise.}
\end{array}
\right.
\label{e1.8aj}
\end{eqnarray}

  It follows from  \cite[(6.12)]{CLR.R}    that for  the  $\bb$-stable process, with $\bb>1$,        
 \begin{equation}
 E\lc\(\int (L^{ x}_{  1})^{ 2}\,dx\)^{ n}\rc\leq C^{ n}( (2n)!)^{ 1/(2\bb)}.\label{rb.1}
 \end{equation}
 When $\psi$ is regularly varying at infinity with index $\bb$,  for all $\ep>0$, there exists a constant $D=D_{\ep}$ such that
\be
\int_{0}^{\ff} e^{-s\psi (p)}\,dp \leq C\(1+\int_{1}^{\ff} e^{-sDp^{\bb-\ep}}\,dp\).
\ee
Using this, and the same   proof as in \cite{CLR.R}, one can show that (\ref{rb.1}) holds, with  $\bb$ replaced by $\bb-\ep$ for any $\ep>0$, when $\psi$ is regularly varying at infinity with index $\bb$.
Therefore, since $\sqrt{(2n)!}\le 2^{n}n!$, 
the right-hand side of (\ref{e1.8aj}), which is the $2n$--th moment of  $ \(8c_{\bb,1}\)^{ 1/2}\sqrt{\al_{ 1}}\,\eta$,
   is bounded above by $  C^{ n}( (2n)!)^{ (\bb+1-\ep)/(2(\bb-\ep))}$. This implies that the weak limit $ \(8c_{\bb,1}\)^{ 1/2}  \sqrt{\al_{ 1}}\,\eta$
  is determined by its moments; (see \cite[p. 227-228]{Feller}).
    Therefore, by the method of moments, \cite[Theorem 30.2]{B}\label{meth}), (\ref{e1.8}) follows from  (\ref{e1.8aj}).

We now obtain (\ref{e1.8aj}).  
 Considering (\ref{lim.1}) and  (\ref{4.9aq}) it suffices to show that  
\be  
 E\(\( \wt I_{l,h} \)^{m}\)\nn 
   =\left\{\begin{array}{ll}\label{e1.8ax}
 \displaystyle {( 2n)!\over 2^{ n}n!}\(4c_{\psi,h,1}\)^{ n} E\lc\(\wt  \al_{l} \)^{
n}\rc & m=2n\\\\
O(h^{\ep}(h\psi^{2}(1/h))^{-n})  & \mbox{otherwise.}
\end{array}
\right. 
\ee 
Using the multinomial theorem we have   
\begin{eqnarray}
&& E\(\( \wt I_{l,h} \)^{m}\)=\sum_{\wt m\in \MM}\({m! \over \prod_{\stackrel{j,k=0}{j< k}}^{l-1}( m_{j,k}! )}\)E\(\prod_{\stackrel{j,k=0}{j< k}}^{l-1} \(I_{j,k,l,h} \)^{m_{j,k}}\),
\label{mnt.1}
\end{eqnarray}
where 
\be \MM=\lc \wt m=\{m_{j,k}, 0\leq j<k\leq l-1\}\,\Bigg |\,
\sum_{\stackrel{j,k=0}{j< k}}^{l-1}m_{j,k}=m \rc.
\ee

We now use Lemma \ref{lem-multiple} to compute the expectation on the right-hand side of (\ref{mnt.1}).  Even though Lemma \ref{lem-multiple} is proved for time intervals of length 1, (see \ref{m.1}),  it is straight forward to check that it holds for  any fixed  time interval, if the term $\al_{j,k}$,  in $(\ref{m.1a})$, is altered to reflect the new length.  Therefore, for some $\ep>0$
\begin{eqnarray} \lefteqn{ 
E\(\( \wt I_{l,h} \)^{m}\)
\label{mn.2}}\\ 
&&
= \sum_{\wt m\in\MM}\({m! \over \prod_{\stackrel{j,k=0}{j< k}}^{l-1}( m_{j,k}!) }\) \prod_{\stackrel{j,k=0}{j< k}}^{l-1}{( 2n_{j,k})!\over 2^{ n_{j,k}}(n_{j,k}!)}\( 4c_{\psi,h,1}\)^{ n_{j,k}} E\lc\prod_{\stackrel{j,k=0}{j< k}}^{l-1}\(\al_{j,k,l}\)^{
n_{j,k}}\rc \nonumber\\
&&\hspace{2 in} +O\(l^{m}h^{(2\bb-1)n+\ep}\) .
\nn
\end{eqnarray}
when $m_{j,k}=2n_{j,k}$ for all $j$ and $k$, and is  $O\(l^{m}h^{(2\bb-1)n+\ep}\) $ if any of the $m_{j,k}$
are odd. Here we   use  the fact that
\begin{equation}
\sum_{\wt m\in\MM}\({m! \over \prod_{\stackrel{j,k=0}{j< k}}^{l-1} (m_{j,k}!) }\)=l^{m}\label{mn.3}
\end{equation}
 to compute the error term.  (Lemma \ref{lem-multiple} is for a fixed partition  of $m$. Here we include the factor $l^{m}$, to account for the number of possible partitions. Note that   $l$ is a function of $h$.)  
 
When $m_{j,k}=2n_{j,k}$ for all $j$ and $k$,  
 \begin{equation}
\({m! \over \prod_{\stackrel{j,k=0}{j< k}}^{l-1} (m_{j,k}! )}\) \prod_{\stackrel{j,k=0}{j< k}}^{l-1}{( 2n_{j,k})!\over 2^{ n_{j,k}}(n_{j,k}!)}={(2n)! \over 2^{n}n!}{n! \over \prod_{\stackrel{j,k=0}{j< k}}^{l-1}(n_{j,k}!)} \label{mn.4}.
 \end{equation}
Using this in (\ref{mn.2}) we get
 \begin{eqnarray} \lefteqn{
E\(\( \wt I_{l,h} \)^{m}\)
\label{mn.5}}\\ 
&&
= {( 2n)!\over 2^{ n}n!}\(4c_{\psi,h,1}\)^{ n}\sum_{\NN}\({n! \over \prod_{\stackrel{j,k=0}{j< k}}^{l-1} n_{j,k}! }\)  E\lc\prod_{\stackrel{j,k=0}{j< k}}^{l-1}\(\al_{j,k,l}\)^{
n_{j,k}}\rc \nonumber\\
&&\hspace{2 in} +O\(l^{m}h^{(2\bb-1)n+\ep}\) ,
\nn
\end{eqnarray}
where $\NN$ is defined similarly as $\MM$. Using the multinomial theorem as in (\ref{mnt.1}) we see that the sum in (\ref{mn.5}) is equal to $E\lc\(\wt\al_{l}\)^{
n}\rc $, which completes the proof of (\ref{e1.8ax}). 
  \qed

 The next   three lemmas give estimates for the  mean and variance of  $\int ( L^{ x+h}_{1}- L^{ x}_{ 1})^{ 2}\,dx$.  They are proved in Section \ref{sec-vare}.

\medskip	Let 
\begin{equation}
c_{\psi,h,0}:=\int_{0}^{\ff}\(p_{s}(0)-p_{s}(h)\)\,ds.\label{lim.7}
\end{equation}

\begin{lemma}\label{lem-exp}   Under the hypotheses of Theorem \ref{theo-clt2r}, 
\begin{equation}
\lim_{h\rar 0}h\psi (1/h)c_{\psi,h,0}=c_{\bb,0}. \label{lim.8}
\end{equation}  
\end{lemma}

\begin{lemma}\label{lem-varep} Under the hypotheses of Theorem \ref{theo-clt2r}; for small $h$ and   $t(h)=1/(\log 1/h)$,   
\begin{equation}
E\(\int ( L^{ x+h}_{t}- L^{ x}_{ t})^{ 2}\,dx\)=4c_{\psi,h,0}t+O\(g(h,t)\)\label{expep}
\end{equation}
as $h\to 0$, where  
\be
g(h,t)=     \left \{\begin{array}{ll}  h^{2}t^{2}\(\psi^{-1}(1/t)\)^{3}&\qquad 3/2<\bb\le2\\\\
 h^{2}L(1/h)&\qquad\bb=3/2\\\\
  (h\psi^{2}(1/h))^{-1}&\qquad1<\bb<3/2 \end{array}  \right.\label{4.9}
  \ee
and   $ L(\cd )$ is   some function that is slowly varying at infinity.
Also
\bea
&&\mbox{Var}\(\int ( L^{ x+h}_{t}- L^{ x}_{ t})^{ 2}\,dx\)   \label{varep}\\
&&\qquad\leq  C\({tg(h,t)\over h\psi(1/h)} +{ t^{2}\psi^{-1}(1/t)  \over h\psi^{2}(1/h)}+ \frac{Ct }{h^{3/2}\psi^{5/2}(1/h)}+ \frac{Ct \log 1/h}{h^{2}\psi^{3}(1/h)}\)\nn.
\eea
\el
 
The proof of this lemma shows that we can take any function  $t:=t(h)$  such that $\psi^{-1}(1/t)<<1/h$ and $ \lim_{h\to 0}t(h)=0$.
 
 \begin{lemma} \label{lem-4.5}
Under the hypotheses of Theorem \ref{theo-clt2r},  
\begin{equation}
E\(\int ( L^{ x+h}_{1}- L^{ x}_{ 1})^{ 2}\,dx\)=4c_{\psi,h,0} +O\(\ov g(h )\)\label{4.11}
\end{equation}
as $h\to 0$, where  
\be
\ov g(h )=     \left \{\begin{array}{ll}  h^{2} &\qquad 3/2<\bb\le2\\\\
 h^{2}L(1/h)&\qquad\bb=3/2\\\\
  (h\psi^{2}(1/h))^{-1}&\qquad1<\bb<3/2 \end{array}  \right.\label{4.11a}
  \ee
and   $ L(\cd )$ is  slowly varying at infinity. 
 \end{lemma}

\noindent\medskip

\noindent{\bf  Proof of Theorem \ref{theo-clt2r} }   We use (\ref{e1.4}) with $l=[\log 1/h]$. Since  $I_{j,j,l,h}$, $0\le j\le l-1$, are  independent and identically distributed, $E( I_{j,j,l,h})=E( I_{0,0,l,h})$, for all $j=0,\ldots,l-1$.
Consequently,
\bea 
&& \int ( L^{ x+h}_{1}- L^{ x}_{ 1})^{ 2}\,dx- E  \int ( L^{ x+h}_{1}- L^{ x}_{ 1})^{ 2}\,dx\\
&&\qquad= 2\wt I_{l}+    \sum_{ j=0 }^{l-1} ( I_{j,j,l,h}    - E( I_{j,j,l,h})) \nn\\
&&\qquad\qquad + \,lE( I_{0,0,l,h})-E  \int ( L^{ x+h}_{1}- L^{ x}_{ 1})^{ 2}\,dx . \nn
 \eea
 We show  immediately below that 
 \begin{equation}
   \lim_{h\to 0}\sqrt{h\psi^{2}(1/h)}\( lE( I_{0,0,l,h})-E  \int ( L^{ x+h}_{1}- L^{ x}_{ 1})^{ 2}\,dx\)=0\label{4.16}.
   \end{equation}
 In addition, using again the property that  $I_{j,j,l,h}$, $0\le j\le l-1$ are  independent and identically distributed, 
   \begin{equation}
   \mbox{Var} \,\(  \sqrt{h\psi^{2}(1/h)}  \sum_{ j=0 }^{l-1} ( I_{j,j,l,h}    - E( I_{j,j,l,h}))   \)= l h\psi^{2}(1/h)   \mbox{Var} \,\(  I_{0,0,l,h}  \)\label{4.32}
   \end{equation}
   We also show  below that 
\begin{equation}
  \lim_{h\to0} l h\psi^{2}(1/h)   \mbox{Var} \,\(  I_{0,0,l,h}  \)=0\label{4.17}.
   \end{equation}
Using   (\ref{4.16})--(\ref{4.17}) and Lemma \ref{lem-tri} we get (\ref{r5.0weaks}).  

  \medskip   	Using   (\ref{expep}) and (\ref{4.11}) on   the expectations in (\ref{4.16}),  and recalling that $l=1/t$, we see that   
\be
lE( I_{0,0,l,h})-E  \int ( L^{ x+h}_{1}- L^{ x}_{ 1})^{ 2}\,dx=O(g(h,t)/t)+O(\ov g(h ))\label{4.17j}
\ee
 It is easy to verify that  (\ref{4.16}) holds.  

\medskip	Showing that (\ref{4.17}) holds is a little more subtle so we provide some details. We first consider the last  three  terms in (\ref{varep}) and  multiply them by   $lh\psi^{2}(1/h)=h\psi^{2}(1/h)/t$ as in (\ref{4.32}). The first of these is
\begin{equation}
  {h\psi^{2}(1/h)\over t}{ t^{2}\psi^{-1}(1/t)  \over h\psi^{2}(1/h)}= t \psi^{-1}(1/t).
   \end{equation}
This last function is regularly varying as $t\to 0$ with index $1-1/\bb$ which is positive by hypothesis.  

 The next term is  
\begin{equation}
   {h\psi^{2}(1/h)\over t}\frac{ t }{h^{3/2}\psi^{5/2}(1/h)}= {1 \over  h^{1/2}\psi^{1/2}(1/h)}.
   \end{equation}
  Here $( h^{1/2}\psi^{1/2})^{-1}$   is regularly varying as $h\to 0$ with index $(\bb-1)/2$ which is positive. 
   
    The third of the last three   terms is 
\begin{equation}
   {h\psi^{2}(1/h)\over t}\frac{ t \log 1/h}{h^{2}\psi^{3}(1/h)}= {\log 1/h \over  h \psi^ (1/h)}.
   \end{equation}
  Here $ (\log 1/h)( h \psi (1/h))^{-1}$   is regularly varying as $h\to 0$ with index $(\bb-1) $ which is positive.  Thus (\ref{4.17}) holds for these three terms.
  
  We now consider  
  \begin{equation}
 {h\psi^{2}(1/h)\over t}{tg(h,t)\over h\psi(1/h)}= { g(h,t)   \psi (1/h)}.\label{ede}
   \end{equation} 
We use (\ref{4.9}) to see that when $\bb>3/2$  this is equal to 
 \begin{equation}
  t^{2}(\psi^{-1}(1/t))^{3}h^{2}\psi(1/h).
   \end{equation}
Here we note that $ t^{2}(\psi^{-1}(1/t))^{3}$ is regularly varying at zero with index $2-(3/\bb) $ which is positive since $\bb>3/2$. In addition by (\ref{2.12}), $\lim_{h\to 0}h^{2}\psi(1/h)<\ff$.

When $\bb=3/2$, (\ref{ede})    is equal to 
 \begin{equation}
 h^{2}L(1/h)\psi(1/h).
   \end{equation}
This function is regularly varying at zero with index  $2-(3/2) $.  

When $\bb<3/2$, (\ref{ede})    is equal to 
$
 (  h\psi (1/h))^{-1},
$
which  is regularly varying at zero with index $\bb -1$.  Thus we have verified (\ref{4.17}).   This completes the proof.
  \qed

\section{Proofs of Lemmas \ref{lem-vproprvt}--\ref{lem-4.3}} \label{sec-est}

Since the  L\'evy processes, $X$, that we are concerned with satisfy 
\begin{equation}
   \int\frac{1}{1+\psi(p)}\,dp<\ff\label{6.1}
   \end{equation}
 it follows from the Riemann Lebesgue Lemma that they have transition probability density functions, which we designate as $p_{s}(\cd)$. Taking the inverse Fourier transform of the characteristic function $X_{s}$, and using the symmetry of $\psi$, we see that 
 \bea
p_{s }(x)&=&{1 \over    2 \pi}\int e^{ipx} \, e^{-s\psi (p) }\,dp\label{69} \\
&=&{1 \over     \pi}\int_{0}^{\ff}\cos px\, e^{-s\psi (p) }\,dp.\nn 
\eea
 
We begin with a technical lemma. 

\begin{lemma}\label{lem-easy} Let    $X $ be a symmetric L\'{e}vy process 
with  L\'{e}vy exponent
 $\psi(\la)$  that is regularly varying at infinity with index $1<\bb\leq 2$ and satisfies (\ref{88.m}).
 Then for any $r\geq 0$   and $t>0$
\begin{equation} 
 \int_{0}^{t} s^{r}e^{-s\psi (p)} \,ds\leq C\(t\wedge {1 \over \psi  (p)}\)^{r+1}\le{2C t^{r+1}\over 1+( t \psi  (p) ) ^{r+1}}; \label{el.1}
\end{equation}
\begin{equation}
   \int _{0}^{\ff}\psi^{r} (p) \(\int_{0}^{1} s^{r}e^{-s\psi (p)} \,ds\)\,dp \le C;\label{el.2} 
\end{equation}
\begin{equation}
 \int | \sin ( hp)| \psi^{r} (p)\(\int_{0}^{1} \,s^{r} e^{-s\psi (p) } \,ds\)\,dp\leq {C \over h\psi (1/h)},\label{el.3}
\end{equation}
and
\begin{equation}
 \int_{0}^{1}   (p_{s }(0)-p_{s }(h)) \,ds  \leq 
C  {1 \over h\psi  (1/h)} \label{el.4}
\end{equation}
as $h\to 0$.

In addition   for all $t\leq 1$ and all $y\in R^{1}$
\begin{equation}
\int_{0}^{t} p_{s }(y )\,\,ds\leq Ct  \psi^{-1}(1/t).\label{pkacv.19}
\end{equation}

\end{lemma}

\Proof  
 The first part of the bound in  the first inequality in  (\ref{el.1}) comes from taking  $e^{-s\psi (p)}\leq 1$;  the second from letting $t=\ff$.  The second inequality in  (\ref{el.1}) is trivial.
 
Note  that for any $y>0$
\begin{equation}
    y^{r}\int_{0}^{1} s^{r}e^{-sy} \,ds =\frac{1}{y} \int_{0}^{y} s^{r}e^{-s} \,ds .
   \end{equation}
Consequently 
 \begin{eqnarray}
 y^{r}\int_{0}^{1} s^{r}e^{-sy} \,ds &\leq &\(\sup_{x\geq 0}x^{r}e^{-x}\)\wedge \(\frac{1}{y} \int_{0}^{\ff} s^{r}e^{-s} \,ds\)
 \label{sr.00}\\
 &\leq &  C   \(1 \wedge {1 \over y}\)\leq 2C\frac{1}{1+y} \nonumber.
 \end{eqnarray}
Using this it is easy to see that
  \begin{eqnarray}
 \int \psi^{r} (p) \int_{0}^{1} s^{r}e^{-s\psi (p)} \,ds\,dp \label{sr.0}
& \leq &C\int   \(1 \wedge {1 \over \psi(p)}\)  \,dp \\
& \leq &C \int_{0}^{1}1 \,dp + C\int_{1}^{\ff} {1 \over  \psi (p)} \,dp    \nonumber
 \end{eqnarray}
 which gives (\ref{el.2}).

Similarly we obtain (\ref{el.3}),
 \bea
 &&
 \int_{0}^{\ff}| \sin ( hp)| \psi^{r} (p)\int_{0}^{1} \,s^{r} e^{-s\psi (p) } \,ds\,dp  \label{e.2}\\
 &&\qquad
\leq C\int_{0}^{\ff}  |\sin (hp )| \(1\wedge {1 \over   \psi( p) }\)\,dp\nn\\
 &&\qquad\leq C \int_{0}^{\ff}   {hp\wedge 1\over 1+\psi(p) }  \,dp  \nn\\
 &&\qquad\leq C\(h \int_{0}^{1/h}   {p\over   1+\psi( p) } \,dp  +  \int_{1/h}^{\ff}   { 1\over   1+\psi( p) } \,dp \)\leq {C \over h\psi (1/h)}.\nn
 \eea
 (In (\ref{e.2}) we use the regular variation of  $\psi$ at infinity. We continue to do so throughout the   rest of this paper without further comment.)

For (\ref{el.4}) we first note that by (\ref{69})
\bea
   p_{s }(0)-p_{s }(h)&=&{1 \over     \pi}\int_{0}^{\ff}(1-\cos ph)\, e^{-s\psi (p) }\,dp\label{6.12}\\
   &=&{2 \over     \pi}\int_{0}^{\ff} \sin^{2} ph/2\, e^{-s\psi (p) }\,dp\nn.
   \eea
Therefore  by Fubini's Theorem and (\ref{el.1}), 
\bea 
  &&  \int_{0}^{1}   (p_{s }(0)-p_{s }(h)) \,ds\label{6.13}\\
  &&\qquad ={2 \over     \pi}\int_{0}^{\ff} \sin^{2} ph/2\,  \int_{0}^{1}  e^{-s\psi (p) }\,ds\,dp\nn\\
    &&\qquad\le C\int_{0}^{\ff}\(  1\wedge  {p^{2}h^{2}\over 2}\)\(  1\wedge{1   \over \psi (p)}\) \,dp\nn\\
        &&\qquad\le Ch^{2}\int_{0}^{1/h }    {p^{2}\over  \psi(p) }   \,dp+C\int_ {1/h }^{\ff}    {1\over  \psi(p) }   \,dp\le C\frac{1}{h\psi(1/h)}\nn .
 \eea

For (\ref{pkacv.19}) we use (\ref{el.1}) to see that
\begin{eqnarray}
\int_{0}^{t} p_{s }(y)\,\,ds&\le &{1 \over 2\pi}\int_{0}^{t} \int e^{-s\psi (p)} \,dp \,\,ds 
\nn\\
&\le& C  \int_{0}^{\ff}\(  t\wedge{1   \over \psi (p)}\) \,dp  \label{pkacv.19w}\\
&\le&C\( t \psi^{-1}(1/t)+  \int_{\psi^{-1}(1/t)}^{\ff} {1  \over \psi (p)} \,dp \)\nonumber\\
&\leq & Ct  \psi^{-1}(1/t) \nonumber.
\end{eqnarray}
\qed

 \noindent{\bf  Proof   of Lemma \ref{lem-vproprvt} }   
 We first note that 
 \begin{equation}
   p_{s }(x)\le C\(\psi^{-1}(1/s) \vee 1\)\label{8.15}.
   \end{equation}
   Refer to (\ref{69}). It is obvious that for $s\ge 1$, $p_{s }(x)\le C$, for all $x$. In addition,
  \bea
p_{s }(x) &\le &{1 \over     \pi} \psi^{-1}(1/s) +{1 \over     \pi}\int_ {\psi^{-1}(1/s)}^{\ff} \, e^{-s\psi (p) }\,dp.  
\eea
Also, for all $s$ sufficiently small,  the last integral is equal to
\be  \int_ {1}^{\ff} \, e^{-u }\,d\psi^{-1}(u/s)\label{8.17}\nonumber\\
<\int_ {1}^{\ff}  \psi^{-1}(u/s) \, e^{-u }\,du
  \ee by integration by parts, where we drop a negative term.    The final integral in (\ref{8.17})  
\bea
    &\le &\psi^{-1}(1/s) \int_ {1}^{\ff} \frac{\psi^{-1}(u/s)}{\psi^{-1}(1/s)}\, e^{-u }\,du\label{bgt}\\
     &\le &\psi^{-1}(1/s) K\int_ {1}^{\ff}  u^{1/\bb+\de} \, e^{-u }\,du\le C\psi^{-1}(1/s) \nn,
   \eea
  for all $\de>0$; where the constant $K$ depends on $\de$. (See e.g. \cite[Theorem 1.5.6]{BGT}.)
 Thus we get  (\ref{8.15}).

  By integration by parts  
 \bea 
  p_{s }(x) 
  & =&{1 \over     \pi x}   \int_{0}^{\ff}e^{-s\psi (p) } \,d(\sin px) \label{8.19} \\
  & = &-{ 1 \over     \pi x}  \int_{0}^{\ff}\sin px\,\left (\frac{d}{dp}e^{-s\psi (p) }\right )\,dp \nn\\
    & =&-\frac{ 1}{\pi x^{2}}  \int_{0}^{\ff}\cos px\left (\frac{d^{2}}{dp^{2}}e^{-s\psi (p) }\right ) \,dp.  \nn
 \eea
 Furthermore 
 \bea
    \frac{d^{2}}{dp^{2}} e^{-s\psi (p) } &=& \(s^{2}(\psi'(p))^{2}-s\psi''(p)\)e^{-s\psi (p) } .
   \eea
Therefore,  by (\ref{114})  
 \begin{equation}
    \bigg| \int_{0}^{1}\cos px\left (\frac{d^{2}}{dp^{2}}  e^{-s\psi (p) } \right ) \,dp   \bigg|\le C\(\int_{0}^{1}\((\psi'(p))^{2}+| \psi''(p) |\)\)\,dp\le C \label{7.17}.
   \end{equation}
   	(We use (\ref{114}) repeatedly in the rest of the paper without comment.)
 In addition, by (\ref{88.m}), for  all $s$ sufficiently small   
  \bea
\bigg|\lefteqn{ \int_{1}^\ff\cos px\left (\frac{d^{2}}{dp^{2}}  e^{-s\psi (p) } \right )\,dp\bigg|\label{7.18}}\\
    && \le C \int_ {1}^\ff\frac{1}{p^{2}}\( \psi^{2}(p)s^{2} e^{-s\psi (p) } +  s \psi(p) e^{-s\psi (p) }  \) \,dp \nn \\
    && \le C \int_ {1}^\ff\frac{1}{p^{2}} \,dp\le C, \nn
   \eea
    since  $\sup_{x\geq 0}x^{r}e^{-x}\leq C$.
    Using (\ref{8.15}) and (\ref{8.19})--(\ref{7.18}) we get (\ref{2.dens}).
   
 \medskip	The inequality in (\ref{2.2}) follows immediately from (\ref{2.dens}).

\medskip	The equality in  (\ref{bl.1}) is  trivial  since  $\int p_{s}(x)\,dx=1$.

\medskip	Note that
\bea
  \De^{ h} p_{s}(x)&= & p_{s }(x+h)-p_{s }(x )\label{70}
\\
  &=&{1\over  \pi}\int_{0}^{\ff} \(\cos p(x+h)- \cos px\)\, e^{-s\psi (p) }\,dp \nn\\
   &=&-{2 \over  \pi}\int_{0}^{\ff}   \cos (px) \sin^{2}(hp/2)e^{-s\psi (p) } \nn\\
   &&\hspace{.5 in}-\frac{ 1}{\pi}\int_{0}^{\ff}   \sin (px) \sin (hp)\,  e^{-s\psi (p) }\,dp \nn   
    \eea
and 
\bea
 \De^{ h}\De^{ -h} p_{s }(x) &= &2p_{s }(x)-p_{s }(x+h)-p_{s }(x-h)\label{72}\\
  &=&{4 \over  \pi}\int_{0}^{\ff}  \cos (px) \sin^{2}(hp /2) \,e^{-s\psi (p) }\,dp.\nn
   \eea
 Thus 
\be 
 \De^{ h} p_{s }(x)\label{72j} =  -{1 \over 2}\De^{ h}\De^{ -h} p_{s }(x)-\frac{ 1} {\pi}\int_{0}^{\ff}   \sin (px ) \sin ( hp) \, e^{-s\psi (p) }\,dp. \label{6.22}
 \ee

 We now note that
 \begin{equation}
 \sup_{x}\int_{0}^{1}|\De^{ h}  p_{s }(x)|\,ds\leq {C \over h\psi (1/h)}\label{6.23}
 \ee
 and
 \be  \sup_{x}\int_{0}^{1}|\De^{ h}\De^{ -h} p_{s }(x)|
 \,ds\leq {C \over h\psi (1/h)}.\label{72j.3}
 \end{equation}

 To obtain (\ref{72j.3}) we use (\ref{72}) to see that
 \begin{equation}
\sup_{x}\int_{0}^{1}|\De^{ h}\De^{ -h} p_{s }(x)|\,ds\le {4 \over  \pi}\int_{0}^{1}\int_{0}^{\ff}   \sin^{2}(hp /2) \,e^{-s\psi (p) }\,dp
 \,ds.    
   \end{equation}
 Using the calculation in (\ref{6.13}) we get (\ref{72j.3}).
  
  To obtain (\ref{6.23}) we note that by (\ref{el.1}),  similarly to (\ref{6.13})
   \bea
 &&  \sup_{x}\int_{0}^{1}\bigg|\int_{0}^{\ff}   \sin (px ) \sin ( hp) \, e^{-s\psi (p) }\,dp\bigg|\,ds\label{6.26}\\
 &&\qquad\le C\int_{0}^{\ff}\(  1\wedge  {p h \over 2}\)\(  1\wedge{1   \over \psi (p)}\) \,dp\nn\\
 &&\qquad\le C\(h\int_{0}^{1/h }    {p \over 1+\psi(p) }   \,dp+\int_ {1/h }^{\ff}    {1\over  \psi(p) }   \,dp\)\le C\frac{1}{h\psi(1/h)} \nn.
   \eea
 Thus (\ref{6.23}) follows from (\ref{6.22}), (\ref{72j.3}) and (\ref{6.26}).
  
   \medskip	
  
 We now show that 
\begin{equation}
\De^{ h}\De^{ -h} p_{s }(x)=\frac{8}{\pi}{K \over x^{2} }\label{89.1}
\end{equation}
where
\begin{equation}
K=K(s,x,h):= \int_{0}^{\ff}\sin^{2}( px/2)\( \sin^{2}(hp/2)\, e^{-s\psi (p) }\)''\,dp.\label{89.2}
\ee
To get this we integrate by parts 
in (\ref{72}),
\bea
  &&   \int_{0}^{\ff}  \cos p x \sin^{2}(hp/2) \, e^{-s\psi (p) }\,dp \label{89}\\
  &&\qquad= \frac{1}{x}   \int_{0}^{\ff} \sin^{2} (hp/2) \, e^{-s\psi (p) }\,d(\sin px)\nn \\
   &&\qquad= -\frac{1}{x}   \int_{0}^{\ff}\sin px\( \sin^{2}(hp/2)\, e^{-s\psi (p) }\)'\,dp\nn\\
   &&\qquad=  -\frac{1}{x}   \int_{0}^{\ff} \( \sin^{2}(hp/2)\, e^{-s\psi (p)} \)'\,d\( \int_{0}^{p}\sin  rx\,dr\)\nn\\
       &&\qquad=  -\frac{1}{x^{2}}   \int_{0}^{\ff} \( \sin^{2}(hp/2) \, e^{-s\psi (p) }\)'\,d\(1-\cos  px \)\nn\\
         &&\qquad= \frac{2}{x^{2}}   \int_{0}^{\ff}\sin^{2}( px/2)\( \sin^{2}(hp/2)\, e^{-s\psi (p) }\)''\,dp\nn.
   \eea

 Let $g(p)= e^{-s\psi (p) }$ and note that
\begin{equation}
 \( 2\sin^{2} (hp/2)\, e^{-s\psi (p) }\)'=  g(p)h\sin hp +2g'(p)\sin^{2} (hp/2)\label{90}
   \end{equation}
and
 \begin{eqnarray}
 &&
 \( 2\sin^{2} (hp/2)\, e^{-s\psi (p) }\)''\label{91}\\
 &&\qquad= g(p)h^{2}\cos hp +2 g'(p)h\sin hp+2g''(p)\sin^{2} (hp/2).\nn
   \end{eqnarray}

Substituting  (\ref{91}) in (\ref{89}) we  write $K= I+ II + III$. 
Using   (\ref{el.1})   we see that 
      \bea 
  \int_{0}^{1} | I|\,ds &= &h^{2} \int_{0}^{1} \Big |\int_{0}^{\ff}\cos hp \sin^{2}( px/2)e^{-s\psi (p) }\,dp\Big |\,ds \label{88.12}\\
    &\le&  h^{2} \int_{0}^{\ff} \(\int_{0}^{1}e^{-s\psi (p) }\,ds\)\,dp\nn\\
    &\le&Ch^{2}\int_{0}^{1}\frac{1}{1+\psi(p)}\,dp=O\(h^{2}\).\nn
     \eea
     Then using (\ref{88.m}), (\ref{114})  and  (\ref{el.2}) with $r=1$  we get 
         \bea
      \int_{0}^{1} | II|\,ds &=& 2h\int_{0}^{1} \Big |\int_{0}^{\ff}\sin hp \sin^{2}( px/2)g'(p)\,dp\Big |\,ds \label{88.13} \\
    &\le& 2h\int_{0}^{\ff}|\sin (hp) \, \psi '(p)|\,\(\int_{0}^{1} s e^{-s\psi (p) }\,ds \)\,dp\nn \\
     &\le&  C h^{2} \int_{0}^{\ff}| p \, \psi '(p)|\,\(\int_{0}^{1} s e^{-s\psi (p) }\,ds \)\,dp\nn\\
       &\le&  Ch^{2} \( C_{1}+\int_{1}^{\ff}  \, \psi (p)\,\(\int_{0}^{1} s e^{-s\psi (p) }\,ds \)\,dp\)=O\(h^{2}\)\nn.
          \eea
  Similarly, and also using (\ref{el.2}) with $r=1$  we get 
  
        \bea 
    \int_{0}^{1} | III|\,ds &=& 2\int_{0}^{1} \Big |\int_{0}^{\ff}\sin^{2} (hp/2) \sin^{2}( px/2)g''(p)\,dp \Big |\,ds 
    \nn\\
     &\le &Ch^{2} \int_{0}^{\ff}p^{2}\(\int_{0}^{1} \(  s |\psi''(p)|   + s^{2} |\psi'(p)|^{2}  \)e^{-s \psi (p) }\,ds \)  \,dp  
    \nn\\
    &\le  &  Ch^{2}\lc \int_{0}^{1}	p^{2}\( |\psi''(p)|   +   |\psi'(p)|^{2} \)\,dp\right.\nn\\
    && \left. +\int_{1}^{\ff} \(\int_{0}^{1} \(  s \psi(p)   + s^{2} \psi^{2}  (p)\)e^{-s \psi (p) }\,ds \)  \,dp\rc  \nn\\
&=& O\(h^{2}\).\label{6.34}
 \eea
Combining (\ref{88.12})--(\ref{6.34}) with (\ref{89.1}) we get the   third bound in (\ref{bl.3}). The first bound in (\ref{bl.3}) follows from (\ref{72j.3}).

To get the second  bound in (\ref{bl.3}) we use  (\ref{72}) and the third integral in (\ref{89})  to see that 
\begin{equation}
\De^{ h}\De^{ -h} p_{s }(x)=-\frac{4}{\pi} {L \over x  }
\label{89.1j}
\end{equation}
where
\begin{equation}
L=L(s,x,h):=      \int_{0}^{\ff}\sin px\( \sin^{2}(hp/2) \, e^{-s\psi (p) }\)'\,dp.\label{89.2k}
\end{equation}
   
  Using (\ref{90}), (\ref{88.m}), (\ref{114}) and (\ref{el.3}) with $r=0$ and 1, we see that 
        \bea 
  \lefteqn{
 \int_{0}^{1}  |L | \,ds\label{90.mm}}\\&&\le C \int_{0}^{1}\(h\int_{0}^{\ff } |  \sin hp|  g(p)\,dp+\int_{0}^{\ff }  \sin^2{(hp/2)}|  g'(p)|\,dp\)\,ds\nn\\
  & &\le Ch \int_{0}^{\ff} |  \sin hp|  \int_{0}^{1} e^{-s\psi (p) }\,ds\,dp\nn\\
  &&\qquad+Ch\(C_{1}+\int_{1}^{\ff } | \sin {(hp/2)}|   |p\psi'(p)|   
\int_{0}^{1}   s e^{-s\psi (p) }\,ds \,dp\) \nn\\
  & &  \le O\({1 \over \psi(1/h)}\)+Ch\int_{0}^{\ff}  |\sin {(hp/2)}|   \psi(p)    
\int_{0}^{1}   s e^{-s\psi (p) }\,ds \,dp\nn\\
&&\le O\({1 \over \psi(1/h)}\).\nn
   \eea
  Thus we get the second bound on the right--hand side of (\ref{bl.3}). This completes the proof of  (\ref{bl.3}).  

\medskip	
To prove (\ref{bl.2}) we first note  that by (\ref{6.22}) it is less than $w(x)/2$  plus
  \begin{equation}
 C \int_{0}^{1} \Big |\int_{0}^{\ff}  \sin (px ) \sin (hp) \, e^{-s\psi (p) }\,dp \Big |\,ds\label{5.34}
\end{equation}
Integrating by parts twice we obtain
  \bea
  &&   \int_{0}^{\ff} \sin (p x) \sin (hp)\, e^{-s\psi (p) }\,dp \label{89s}\\
  &&\qquad=- \frac{1}{x}   \int_{0}^{\ff} \sin (hp) \, e^{-s\psi (p) }\,d(\cos px)\nn \\
   &&\qquad= \frac{1}{x}   \int_{0}^{\ff}\cos px\(\sin (hp) \, e^{-s\psi (p) })\)'\,dp\nn\\
   &&\qquad=  \frac{1}{x^{2}}   \int_{0}^{\ff} \(\sin (hp) \, e^{-s\psi (p) }\)'\,d\(  \sin  px \)\nn\\
         &&\qquad= -\frac{1}{x^{2}}   \int_{0}^{\ff}\sin px\(\sin (hp) \, e^{-s\psi (p) }\)''\,dp.\nn
   \eea
 Note that
\be  \(\sin (hp) \, e^{-s\psi (p) }\)'  =\(  h \cos hp  -\sin hp(s \,\psi'(p) \)e^{-s\psi (p) } 
   \ee

 Thus   the left hand side of (\ref{89s}) is bounded by $\displaystyle {J \over x  }$
where
\begin{equation}
J=J(s,x,h):= \int_{0}^{1} \Big | \int_{0}^{\ff}\cos px\(\sin (hp) \, e^{-s\psi (p) }\)'\,dp\Big |\,ds.\label{89.2.1}
\end{equation}
We write 
\begin{equation}
   J\le J_{1}+J_{2}
   \end{equation}
   where 
   \bea 
   |J_{1}|&\le &h\int_{0}^{1}  \int_{0}^{\ff}|\cos px \cos (hp)| \, e^{-s\psi (p) }\,dp \,ds\\
   &\le &Ch \int_{0}^{\ff}\frac{1}{1+\psi(p)}\,dp\le C'h\nn
 \eea
   and  using  (\ref{88.m}), (\ref{114}) and (\ref{el.2})
   \bea 
      |J_{2}|&\le & \int_{0}^{1}  \int_{0}^{\ff}|\cos px \sin (hp)|  |\psi'(p)|\, se^{-s\psi (p) }\,dp \,ds\\
      &\le & h\int_{0}^{1}  \int_{0}^{\ff}p |\psi'(p)|\, se^{-s\psi (p) }\,dp \,ds\nn\\
&\le & h \int_{0}^{1}   |\psi'(p)|\,  \int_{0}^{1}se^{-s\psi (p) }\,ds\,dp\nn \\    & &+ C h\int_{0}^{1}  \int_{1}^{\ff} \psi(p)\, se^{-s\psi (p) }\,dp \,ds\le C'h.\nn
 \eea
 Therefore
 \begin{equation}
   {J \over |x | }\le C\frac{h}{|x |}.\label{6.46}
   \end{equation}
  
In addition    (\ref{5.34}) is $\displaystyle {G \over x^{2} }$
where
\begin{equation} 
G=G(x,h):= \int_{0}^{1} \Big | \int_{0}^{\ff}\sin px\(\sin (hp) \, e^{-s\psi (p) }\)''\,dp\Big |\,ds.\label{89.2.1.1}
\end{equation}
Since
\bea
\lefteqn{ \(\sin (hp) \, e^{-s\psi (p) }\)''}\label{6.43}\\
&&\hspace{-.1in}=\( -h^{2}\sin hp +2hs\cos hp\,\psi'(p)-\sin hp(s \,\psi''(p)-s^{2}(\psi'(p))^{2}\)e^{-s\psi (p) }\nn,
   \eea
 we can write 
\be
G\leq G_{1}+G_{2}+G_{3}.
\ee          

Using  (\ref{6.43}) and  (\ref{el.2}) we get
  \begin{eqnarray}
  |G_{1}|&=&h^{2}\int_{0}^{1} \Big |\int_{0}^{\ff}\sin px\(\sin (hp) \, e^{-s\psi (p) }\) \,dp\Big |\,ds\label{89.2e}\\
  & 
  \leq& Ch^{2} \int_{0}^{\ff} \int_{0}^{1}  e^{-s\psi (p) }\,ds \,dp\leq  Ch^{2}.  \nonumber
  \end{eqnarray}

 Using    (\ref{88.m}) , (\ref{114}) and (\ref{el.2}) we see that  
    \begin{eqnarray}
  |G_{2}|&=&2h\int_{0}^{1} \Big |\int_{0}^{\ff}\sin px \cos hp\( \psi'(p)\,s e^{-s\psi (p)}\) \,dp\Big |\,ds
  \label{89.2f}\\
  &\leq &2h\int_{0}^{1}  \int_{0}^{\ff}|\psi'(p)|\,s e^{-s\psi (p)}  \,dp \,ds\nn\\
  &\le&   2h \(C_{1}+ \int_{1}^{\ff}p|\psi'(p)|\(\,\int_{0}^{1}s e^{-s\psi (p)}\,ds\) \,dp\) \nonumber\\
  &\le&   Ch \(C_{1}+\int_{1}^{\ff}  \psi(p)\(\,\int_{0}^{1}s e^{-s\psi (p)}\,ds\)  \,dp\) \le Ch.\nn
    \eea

Similarly 
    \begin{eqnarray}
 \lefteqn{|G_{3}| \label{89.2g}}\\
  &&=\int_{0}^{1} \Big |\int_{0}^{\ff}\sin px\,\, \sin hp \( s\psi''(p)   - s^{2} (\psi'(p))^{2}    \)e^{-s\psi (p) } \,dp\Big |\,ds
  \nn\\
    &&\leq  h\int_{0}^{\ff}p \(\int_{0}^{1}\( s|\psi''(p)| + s^{2} (\psi'(p))^{2}   \)e^{-s\psi (p) } \,ds\) \,dp
  \nn\\
    &&\leq C h \lc C_{1}+ \int_{1}^{\ff}p^{2}\(\int_{0}^{1}\( s|\psi''(p)| + s^{2} (\psi'(p))^{2}   \)e^{-s\psi (p) } \,ds\) \,dp \rc
  \nn\\
  &&\leq C h \lc C_{1}+ \int_{1}^{\ff} \(\int_{0}^{1}\( s\psi (p) + s^{2} (\psi(p))^{2}   \)e^{-s\psi (p) } \,ds\) \,dp \rc
  \nn\\
  & & \le Ch\nonumber .
  \end{eqnarray}
  Thus we see that for all  $|x|> 0$ 
  \begin{equation}
  {G  }\le C{h } \label{6.54a},
   \end{equation}
   for some $C<\ff$ independent of $|x|$.
Combining (\ref{6.23}), (\ref{6.46}) and (\ref{6.54a}) and taking into account the value of $w(x)$, we get   (\ref{bl.2}).

\medskip	  For  (\ref{bl.01}) we use (\ref{bl.2}) to see that  
  \bea
&&\int  \(\,\int_{0}^{1}\,|\De^{ h}  \,p_{s }(x)|\,ds\)\,dx\label{6.59q}\\
&  &\qquad\le  C\( \int_{0}^{a}\frac{1}{h\psi(1/h)}\,dx+h\int_{a}^{1}\frac{1}{x }\,dx+h\int_{1}^{\ff}\frac{1}{x^{2} }\,dx\),\nn
\eea
  Set $a=a(h)=h^{2}\psi(1/h )$. For L\'evy processes excluding Brownian Motion,     $\lim_{h\to 0}h^{2}\psi(1/h)=0$; (see \cite[Lemma 4.2.2]{book}),  and we can estimate (\ref{6.59q}) to obtain (\ref{bl.01}). For Brownian Motion take $a=1$  in (\ref{6.59q}) to obtain (\ref{bl.01}).

   Similarly, to obtain (\ref{lastf}) we use (\ref{bl.2}) to get
  \bea
\int  \(\,\int_{0}^{1}\,|\De^{ h}  \,p_{s }(x)|\,ds\)^{p}\,dx&\le & C\( \int_{0}^{a}\frac{1}{h^{p}\psi^{p}(1/h)}\,dx+h^{p}\int_{a}^{\ff}\frac{1}{x^{p}}\,dx\)\nn\\
&\le & C\(  \frac{a}{h^{p}\psi^{p}(1/h)} +{h^{p}\over a^{p-1}} \)\label{6.59}.
\eea

\medskip	  For (\ref{bl.5}) we use  (\ref{bl.3}) to see that
  \begin{eqnarray}
  &&\int  \(\,\int_{0}^{1}\,|\De^{ h}\De^{ -h} \,p_{s }(x)|\,ds\)^{2}\,dx
  \label{72j.3n}\\
  &&\qquad =\int_{0}^{h}  \(\,\int_{0}^{1}\,|\De^{ h}\De^{ -h} \,p_{s }(x)|\,ds\)^{2}\,dx  \nonumber\\
  && \qquad\qquad+\int_{h}^{\ff}  \(\,\int_{0}^{1}\,|\De^{h}\De^{ -h} \,p_{s }(x)|\,ds\)^{2}\,dx  \nonumber\\
  && \leq {C \over h\psi^{2}(1/h)}   +{C \over \psi^{2}(1/h)}\int_{h}^{\ff}  {1\over x^{2}} \,dx =O\({1 \over h\psi^{2}(1/h)}\)\nn.
  \end{eqnarray}
The inequality in  
(\ref{bl.5b}) follows similarly,
  \be 
 \int_{u}^{\ff}  \(\,\int_{0}^{1}\,|\De^{ h}\De^{ -h} \,p_{s }(x)|\,ds\)^{2}\,dx
 \leq   {C \over \psi^{2}(1/h)}\int_{u}^{\ff}  {1\over x^{2}} \,dx ={C \over u\psi^{2}(1/h)}.  \label{72j.3nv}
  \ee

To obtain  (\ref{bl.4}) we use   (\ref{bl.3}) to see that 
\begin{eqnarray}
&&\int \int_{0}^{1} \Big |\De^{ h}\De^{ -h}\,p_{s }(x)\Big |\,ds\,dx
\label{st.30}\\
&&\qquad= \int_{0}^{h}\int_{0}^{1} \Big |\De^{ h}\De^{ -h}\,p_{s }(x)\Big |\,ds\,dx
+\int_{h}^{1}  \int_{0}^{1} \Big |\De^{ h}\De^{ -h}\,p_{s }(x)\Big |\,ds\,dx\nn\\
&&\hspace{1 in}
+\int_{1}^{\ff} \int_{0}^{1} \Big |\De^{ h}\De^{ -h}\,p_{s }(x)\Big |\,ds\,dx  \nonumber\\
&&\qquad\leq {C \over h\psi (1/h)}   \int_{0}^{h}1\,dx+{C \over \psi (1/h)}\int_{h}^{1} {1 \over |x|} \,dx
+Ch^{2}\int_{1}^{\ff}{1 \over |x|^{2}}\,dx  \nonumber\\
&&\qquad\leq {C \over \psi (1/h)}  +{C \log 1/h\over \psi (1/h)}+Ch^{2} \nn \end{eqnarray}
\qed


\noindent
{\bf  Proof of Lemma \ref{lem-h3}$\,$}  Using $2-e^{iph}-e^{-iph}=4\sin^{2}(hp/2)$ we see that  
\bea
 \int_{0}^{\ff} \De^{ h}\De^{ -h} p_{t}(x )\,dt& =&{1 \over  2\pi}\int_{0}^{\ff} \int e^{-ipx}(2-e^{iph}-e^{-iph} )e^{-t\psi(p)}\,dp\,dt\nn\\
 &=&  {4 \over  2 \pi} \int e^{-ipx}\sin^{2}(hp/2)\int_{0}^{\ff}e^{-t\psi(p)}\,dt\,dp\label{par.3} \\
&=& {4\over  2\pi}\int e^{-ipx}{\sin^{2}(hp/2) \over \psi(p)}\,dp\nn.
\eea
It follows from    Parseval's Theorem that  
\be 
c_{\psi,h,1}= \int  \(\int_{0}^{\ff} \De^{ h}\De^{ -h} p_{t}(x )\,dt\)^{2}\,dx
\label{par.4} =  {8 \over \pi} \int {\sin^{4}(hp/2) \over \psi^{2}(p)}\,dp.
\ee
 Using this we write  
\be
h\psi^{2}(1/h)c_{\psi,h,1}
=  {16 \over  \pi}\int_{0}^{\ff} \({\sin^{2}(p/2 )\over \psi  (p/h)/\psi (1/h)}\)^{2}\,dp \label{5.65}
\ee
For a fixed $0<a<1$,  
\bea 
\int_{0}^{a}\({\sin^{2}p/2 \over \psi  (p/h)/\psi (1/h)}\)^{2}\,dp  &=&h\psi ^{2}(1/h)  \int_{0}^{a/h}\({\sin^{2}(ph/2) \over \psi  (p)}\)^{2}\,dp\nn\\
&\le&{ h^{5}\psi ^{2}(1/h)\over 4}  \int_{0}^{a/h}{p^{4}  \over \psi^{2}  (p)} \,dp\label{5.60}.
\eea
For any $\ep>0$ we can find an $h_{0}>0$, such that for all $0<h\le h_{0}$, the  last line  above

\be  
\le {(1+\ep) h^{5}\psi ^{2}(1/h)\over 4(5-2\bb)}  {(a/h)^{5}  \over \psi^2({a/h}  ) }\le \frac{a^{5-2\bb}}{2(5-2\bb)} \label{6.52}.
\ee
Note that  for any $\ep'>0$  and $p\ge a>0$, we can find an $h'_{0}>0$, such that for all $0<h\le h'_{0}\le h_{0}$,   
\be
{ \psi^2({1/h}  )  \over \psi^2({p/h}  ) }\le C\max\({1\over p^{2\bb-\ep}}, {1\over p^{2\bb+\ep}}\) \label{6.53}.
 \ee
  (See   \cite[Theorem 1.5.6]{BGT}.) 
 Therefore, it follows from  the Dominated Convergence Theorem that
\begin{equation}
\lim_{h\to 0}   \int_{a}^{\ff}\({\sin^{2}p/2 \over \psi  (p/h)/\psi (1/h)}\)^{2}\,dp=     \int_{a}^{\ff} {\sin^{4}p/2 \over p^{2\bb}} \,dp\label{6.55}.
   \end{equation}
Since (\ref{5.60}), (\ref{6.52}) and (\ref{6.55}) hold for all $a>0$ sufficiently small, we get (\ref{lim.1}).

\medskip	We now consider  (\ref{lim.1r}).   Just as we obtained (\ref{par.3}) and (\ref{par.4}) we see that
  \bea 
\lefteqn{\int_{[0,\sqrt{h}]^{2}}\int \(\De^{ h}\De^{ -h}\,p_{r }(x)\)
\(\De^{ h}\De^{ -h}\,p_{r'}(x)\) \,dx \,dr \,dr'\nn} \\
&  &  ={8\over \pi}\int{\sin^{4}(ph/2) \over \psi^{2} (p)}\( 1-e^{-\sqrt{h} \psi (p)}  \)^{2}  \,dp. \label{3.41} 
 \eea
 We show below that  
 \begin{equation}
   h\psi^{2}(1/h)\int {\sin^{4}(ph/2) \over \psi^{2} (p)}e^{-\sqrt{h} \psi (p)}   \,dp
=O(h^{1/2 }),\label{6.57}
   \end{equation}
which proves (\ref{lim.1r}).
 
 To obtain (\ref{6.57}) we note that 
 \begin{eqnarray}
 \lefteqn{h\psi^{2}(1/h)\int {\sin^{4}(ph/2 )\over \psi^{2} (p)}e^{-\sqrt{h} \psi (p)}   \,dp
 \label{lim.5}}\\
 && = h\psi^{2}(1/h)\int_{0\leq |p|\leq 1}  {\sin^{4}(ph/2 )\over \psi^{2} (p)}e^{-\sqrt{h} \psi (p)}   \,dp \nonumber\\
 && \qquad + h\psi^{2}(1/h)\int_{1\leq |p|\leq 1/h} {\sin^{4}(ph/2) \over \psi^{2} (p)}e^{-\sqrt{h} \psi (p)}   \,dp \nonumber\\
 && \qquad\quad + h\psi^{2}(1/h)\int_{ |p|\geq 1/h}  {\sin^{4}(ph/2) \over \psi^{2} (p)}e^{-\sqrt{h} \psi (p)}   \,dp \nonumber\\
 && \leq C h^{5}\psi^{2}(1/h)\int_{0\leq |p|\leq 1} {p^{4}\over \psi ^{2}(p)} \,dp \nonumber\\
 && \qquad + C h^{5}\psi^{2}(1/h)\int_{1\leq |p|\leq 1/h}   {p^{4}\over \psi ^{2}(p)} {1 \over \sqrt{h} \psi (p)}  \,dp \nonumber\\
 && \qquad \quad+C h\psi^{2}(1/h)e^{-\sqrt{h} \psi (1/h)}\int_{ |p|\geq 1/h}  {1\over \psi^{2} (p)}   \,dp,\nonumber
 \eea
 where, in the next to last line of (\ref{lim.5}), we use the fact that for $s\ge0$, $e^{-s}\le (\sup_{s\ge 0} se^{-s})/s$. 
 It is obvious that the first and last integral in the last inequality in (\ref{lim.5}) is $O(\sqrt{h})$. As for the second integral, if $1<\bb<5/3$
 \begin{equation}
  h^{5}\psi^{2}(1/h) \int_{1\leq |p|\leq 1/h}   {p^{4}\over \psi ^{2}(p)} {1 \over \sqrt{h} \psi (p)}\,dp \le C \frac{1}{h^{ 1/2}\psi (1/h)}=O\(\sqrt h\,\);
   \end{equation}
    if $ 5/3<\bb\le 2$
 \begin{equation}
  h^{5}\psi^{2}(1/h) \int_{1\leq |p|\leq 1/h}   {p^{4}\over \psi ^{2}(p)} {1 \over \sqrt{h} \psi (p)}\,dp \le C \frac{  h^{5}\psi^{2}(1/h)}{h^{ 1/2} }=O\(\sqrt h\,\),
   \end{equation}
where we use Remark \ref{rem-1} when $\bb=2$.
When $\bb=5/3$ 
 \begin{equation}
  h^{5}\psi^{2}(1/h) \int_{1\leq |p|\leq 1/h}   {p^{4}\over \psi ^{2}(p)} {1 \over \sqrt{h} \psi (p)}\,dp \le  C \frac{L(1/h)}{h^{ 1/2}\psi (1/h)}\le O(h)  \end{equation}
  for some   function   $L(\cd)$ that is slowly varying at infinity. This gives us (\ref{6.57}). \qed

 \begin{lemma} For $r\ge0$
 \be
 \sup_{\de\le s\le 1}s^{r}e^{-s\psi(p)}\le C\(1\wedge \frac{1}{\psi^{r}(p)}\)\le\frac{2C}{1+\psi^{r}(p)},\label{9.1}
 \ee
and for $k>0$
  \be
 \sup_{\de\le s\le 1}s^{r}e^{-s\psi(p)}\le   \sup_{\de\le s\le 1}{s^{r+k}\over\de^{k}}e^{-s\psi(p)}\le \frac{1}{\de^{k}}\frac{2C}{1+\psi^{r+k}(p)}.\label{9.2}
  \ee
 \end{lemma}
 
 \Proof  The first inequality in (\ref{9.1}) follows from the facts that  $y^{r}e^{-y}\le C $ and, of course, $ \sup_{\de\le s\le 1}s^{r}e^{-s\psi(p)}\le 1$.  The second inequality in (\ref{9.1}) is elementary. The inequality in (\ref{9.2}) follows from (\ref{9.1}).\qed
 
  \noindent{\bf Proof  of Lemma \ref{lem-4.2}}  The inequality in (\ref{h.101xx.1}) follows immediately from (\ref{2.dens}).

By (\ref{9.2}) with $r=0$ and $k=3$  
 \bea  
 \sup_{\de\le s\le 1}|\De^{h}p_{s }(0)|&= &\sup_{\de\le s\le 1}{1 \over     \pi}  \int_{0}^{\ff}  \sin^{2}(ph/2)\,   e^{-s\psi (p) } \,dp \label{1.25}\\
 &\le  &\sup_{\de\le s\le 1}{h^{2} \over   2  \pi}  \int_{0}^{\ff} p^{2} \,   e^{-s\psi (p) } \,dp \nn\\
 &\le& C{h^{2} \over  \de^{3}}  \int_{0}^{\ff} {p^{2} \over 1+\psi^{3}(p)}  \,dp       \le {C\over {\de^{3}} }h^{2}\nn,
 \eea
  Note that $\De^{h}p_{r }(0)=p_{r }(h)-p_{r }(0)<0$
 and  $\De^{h}\De^{-h}p_{r }(0)=2(p_{r }(0)-p_{r }(h))$. Thus (\ref{4.5}) follows immediately from (\ref{4.4}).
\qed
 
  \noindent{\bf Proof  of Lemma \ref{lem-4.3}}     
 The inequality in  (\ref{b.l.2wa}) follows immediately from (\ref{2.dens}).

  To obtain (\ref{bl.3q}) consider the material in the proof of Lemma \ref{lem-vproprvt} from (\ref{89.1}) to the statement that $K= I+ II + III$.
 Now, instead of integrating $I$, $II$ and $III$ we take their supremum as $\de\le s\le 1$. We have
   \bea 
\sup_{\de\le s\le 1}   | I| &\le  &h^{2} \sup_{\de\le s\le 1}  \Big |\int_{0}^{\ff}\cos hp \sin^{2}( px/2)e^{-s\psi (p) }\,dp\Big |  \label{88.12a}\\
    &\le&  h^{2}\sup_{\de\le s\le 1} \int_{0}^{\ff} e^{-s\psi (p) } \,dp\nn\\
    &\le&h^{2}\frac{C}{\de}\int_{0}^{\ff}\frac{1}{1+\psi(p)}\,dp\le {C\over {\de}}h^{2}\nn,
     \eea
  where we use (\ref{9.2})    with $r=0$ and $k=1$.
         \bea
\sup_{\de\le s\le 1}  | II|  &=& 2\sup_{\de\le s\le 1} h   \Big |\int_{0}^{\ff}\sin hp \sin^{2}( px/2)g'(p)\,dp\Big |\  \label{88.13a} \\
    &\le& 2h\sup_{\de\le s\le 1}\int_{0}^{\ff}|\sin (hp) \, \psi '(p)|\,    e^{-s\psi (p) } \,dp\nn \\
     &\le&  C h^{2}\sup_{\de\le s\le 1} \int_{0}^{\ff}| p \, \psi '(p)|\,   s e^{-s\psi (p) } \,dp\nn\\
       &\le&Ch^{2}\sup_{\de\le s\le 1} \( C_{1}+\int_{1}^{\ff}  \, \psi (p)\,  s  e^{-s\psi (p) } \,dp\)\nn\\
       &\le&    {C\over\de }h^{2}  \( C'_{1}+\int_{1}^{\ff}  \,{ \psi (p)\over 1+\psi^{2}(p) }   \,dp\)\nn     ,
          \eea
where we use (\ref{9.2}) with $r=1$ and $k=1$. Similarly, but with $r,k=0,1$ and $r,k=2,1$  
 \bea  
\sup_{\de\le s\le 1} | III| &\le &  \sup_{\de\le s\le 1}  \Big |\int_{0}^{\ff}\sin^{2} (hp/2) \sin^{2}( px/2)g''(p)\,dp \Big | 
    \nn\\
     &\le &Ch^{2}  \sup_{\de\le s\le 1}\int_{0}^{\ff}p^{2} \(   |s\psi''(p)|   +  s^{2} |\psi'(p)|^{2}  \)e^{-s \psi (p) }  \,dp  
    \nn\\
    &\le  &  Ch^{2}\lc \int_{0}^{1}	p^{2}\( |\psi''(p)|   +   |\psi'(p)|^{2} \)\,dp\right.\label{6.34a}\\
    && \left. + \sup_{\de\le s\le 1}\int_{1}^{\ff} \(s\psi(p)   +  s^{2} \psi^{2}  (p)\)e^{-s\psi (p) } \,dp\rc  \nn\\
     &\le  &  Ch^{2}\lc C_{1}+  \sup_{\de\le s\le 1} \int_{1}^{\ff} \(s\psi(p)   +   s^{2}\psi^{2}  (p)\)e^{-s\psi (p) } \,dp\rc  \nn\\
         &\le  & { C\over \de }h^{2}\lc C_{1}+   \int_{1}^{\ff}{ \psi(p)\over  1+\psi^{2}(p)}\,dp   + \int_{1}^{\ff} { \psi^{2}  (p)\over 1+\psi^{3}(p)}  \,dp\rc   \nn
 \eea
Combining (\ref{88.12a})--(\ref{6.34a}) with (\ref{89.1}) we get the  second bound in (\ref{bl.3q}). 

The first bound on the right--hand side of (\ref{bl.3q}) follows from (\ref{4.5}) since,   
\be
\left| \sup_{\de\le r\le 1}\De^{h}\De^{-h}p_{r }(x)\right| \le \sup_{\de\le r\le 1}\De^{h}\De^{-h}p_{r }(0)\label{1.40},
\ee
(see (\ref{72}).)
 
 To get the second  bound on the right--hand side of (\ref{b.l.2w}) consider the material in the paragraph containing (\ref{5.34}). For our purposes here we need to obtain  
  \begin{equation}
\sup_{\de\le s\le 1} \Big |\int_{0}^{\ff}  \sin (px ) \sin (hp) \, e^{-r\psi (p) }\,dp \Big | \label{5.34bb}
\end{equation}
Integrating by parts twice as in (\ref{89s}) we see that (\ref{5.34bb}) is bounded by
  \be 
 \sup_{\de\le s\le 1}  \left|\frac{1}{x^{2}}   \int_{0}^{\ff}\sin px\(\sin (hp) \, e^{-s\psi (p) }\)''\,dp\right|.\nn
   \ee 
Thus we have to take $\sup_{\de\le s\le 1}$ of the terms in (\ref{89.2e})--(\ref{89.2g}), but without the integral on $s$.   It is easy to see that we get the same bounds as in  (\ref{89.2e})--(\ref{89.2g}) but with the   factor  $1/\de$ as in (\ref{88.12a})--(\ref{6.34a}),
 
  By  (\ref{70}), (\ref{b.l.2w}) is bounded by  (\ref{4.5}) plus 
\bea 
  &&C h\int p\,e^{-\de\psi(p)}\,dp\label{5.85}\\
  &  &\qquad\le  C h \(\int_{0}^{\psi^{-1}(1/\de)}  p \,dp+\frac{1}{\de^{2}}\int_ {\psi^{-1}(1/\de)}^{\ff}  \frac{p}{\psi^{2}(p)} \,dp\)\nn\\
  && \nn \qquad\le  C h  ( \psi^{-1}(1/\de) )^{2}\le  C\frac{h}{\de^{2}}  .
  \eea
(For the second integral in the middle line of (\ref{5.85}) see the comment following (\ref{lim.5}).) This gives 
  the first bound on the right--hand side of (\ref{b.l.2w}).

 The inequalities in (\ref{bl.1x})--(\ref{bl.5ba})  follow easily from 
 (\ref{b.l.2wa})--(\ref{bl.3q}). \qed

\section{Proofs of Lemmas \ref{lem-exp}--\ref{lem-4.5}}\label{sec-vare}
 \medskip	
\noindent {\bf Proof of Lemma \ref{lem-exp} }  By (\ref{69}) 
\begin{eqnarray}
h\psi (1/h)c_{\psi,h,0}
 &=& {h\psi (1/h) \over   \pi}\int_{0}^{\ff} {1-\cos (ph)  \over \psi  (p)}\,dp \label{lim.12}\\
&=& {2 h\psi (1/h) \over \pi}\int_{0}^{\ff} { \sin^{2} (ph/2)  \over \psi  (p)}\,dp \nonumber\\
&=& {2\over \pi}\int_{0}^{\ff}  {\sin^{2} (p/2) \over \psi  (p/h)/\psi (1/h)} \,dp \nonumber.
\end{eqnarray}
Compare this to (\ref{5.65}). Following the proof of (\ref{lim.1}), from (\ref{5.60}) to (\ref{6.55}), with obvious modifications, we get (\ref{lim.8}).

\begin{remark} \label{rem-6.1}{\rm We note that by (\ref{lim.12}), for symmetric stable processes of index $\bb$, we get equality in (\ref{lim.8}), namely
\begin{equation}
   c_{\psi,h,0}=h^{\bb-1}c_{\bb,0}.
   \end{equation}
 }\end{remark}
	
\noindent{\bf  Proof of Lemma \ref{lem-varep} } We show in   \cite[(8.3)--(8.4)]{CLR} that by the Kac moment formula,
for    $0<t\le 1$, 
 \bea 
 E\(\int ( L^{ x+h}_{t}- L^{ x}_{ t})^{ 2}\,dx\)\label{pkacv.4}& =& 4  \int_{0}^{t} (t-r)\(p_{r }( 0)-p_{r }( h)\) \,dr \\ 
  &=&  \frac{8}{\pi}\int_{0}^{\ff} \sin^{2}(hp/2 )\int_{0}^{t} (t-r)e^{-r\psi(p)} \,dr\,dp.\nn
\eea
(In truth we show it for $t=1$ but it is obvious that it holds for any $t$.)
Note that  
\begin{equation}
   \int_{0}^{t} (t-r)e^{-r\psi(p)}\,dr=\frac{t}{\psi(p)}-\frac{1-e^{-t\psi(p)}}{\psi^{2}(p)}.
   \end{equation}
By (\ref{lim.12})  
\begin{equation}
   {8 t\over \pi}\int_{0}^{\ff}   {\sin^{2} (ph/2) \over \psi(p)}  \,dp  =4c_{\psi,h,0}t.\label{7.6}
   \end{equation}
   This gives the dominant term in (\ref{expep}). 
The absolute value of the  remainder is  
\begin{equation}
  {8   \over \pi}\int_{0}^{\ff}   {\sin^{2} (ph/2) \over \psi^{2}(p)} \(1-e^{-t\psi(p)}\) \,dp \le    {8   \over \pi} \int _{0}^{\ff}  {\sin^{2} (ph/2) \over \psi^{2}(p)} \(1\wedge t\psi(p)\) \,dp .\label{7.7}
   \end{equation}
   We break this last integral into three parts and see that it is bounded by
   \be C\(  h^{2}t \int_{0}^{\psi^{-1}(1/t)} {p^{2}\over \psi (p)} \,dp +   h^{2}  \int_ {\psi^{-1}(1/t)} ^{1/h}{p^{2} \over \psi^{2} (p)} \,dp +\int_{1/h}^\ff  {1 \over \psi^{2} (p)} \,dp\)\label{7.8}
 \ee 
 We have
   \begin{equation}
    h^{2}t \int_{0}^{\psi^{-1}(1/t)} {p^{2}\over \psi (p)} \,dp\le Ch^{2}t^{2}\(\psi^{-1}(1/t)\)^{3}\label{7.9},
   \end{equation}
  (Since $\lim_{p\to 0}\psi(p)/p^{2}>0$    this integral is finite; (see \cite[Lemma 4.2.2]{book}).    
  
     In addition 
   \begin{equation}
  \int_{1/h}^\ff  {1 \over \psi^{2} (p)} \,dp\le C\frac{1}{h\psi^{2}(1/h)}\label{7.10}
   \end{equation}
If $\bb>3/2$
\begin{equation}
    h^{2}  \int_ {\psi^{-1}(1/t)} ^{1/h}{p^{2} \over \psi^{2} (p)} \,dp \le Ch^{2}t^{2}\(\psi^{-1}(1/t)\)^{3}.\label{7.11}
   \end{equation}
If $\bb=3/2$
\begin{equation}
    h^{2}  \int_ {\psi^{-1}(1/t)} ^{1/h}{p^{2} \over \psi^{2} (p)} \,dp \le Ch^{2}L(1/h) \label{7.12}  \end{equation}
for some function $L(\cd)$ that is slowly varying at infinity. 
If $\bb<3/2$
\begin{equation}
    h^{2}  \int_ {\psi^{-1}(1/t)} ^{1/h}{p^{2} \over \psi^{2} (p)} \,dp \le C\frac{1}{h\psi^{2}(1/h)}. \label{7.13}\end{equation}
Using (\ref{7.7})--(\ref{7.13}) we get (\ref{4.9}). 

\medskip	Let 
\be
Z=\int ( L^{ x+1}_{t}- L^{ x}_{ t})^{ 2}\,dx.\label{7.14}
\ee
 We get an upper bound for the variance of  $Z$ by finding an upper bound for $EZ^{2}$ and using  (\ref{expep}) to estimate $(EZ)^{2}$.  We proceed as in the beginning of the proof of Lemma \ref{lem-multiple}, however there are enough differences that it is better to repeat some of the arguments.
 
 By the   Kac  Moment Theorem  
\begin{eqnarray}
 & &
E\(\prod_{i=1}^{2}\(\De^{h}_{x_{i}}L_{t}^{x_{i}}\)\(\De^{h}_{y_{i}}L_{t}^{y_{i}}\)\)\label{kacv}\\
& &\quad= \prod_{i=1}^{2}\(\De^{h}_{x_{i}}\De^{h}_{y_{i}}\)\sum_{\si}\int_{\{\sum_{i=1}^{4}r_{i}\leq t\}}\prod_{i=1}^{4}p_{r_{i}}( \si(i) - \si(i-1))\,\,    \prod_{i=1}^{4}\,dr_{i} \nn
\end{eqnarray}
where the sum runs over all  bijections    $\si:\,[1,4]\mapsto \{x_{i},y_{i},\,1\leq i\leq 2\}$   and we take $\si(0)=0$.   
We rewrite (\ref{kacv}) 
 so that each   $\De^{h}_{\cd}$ applies to a single $p_{\cd}$ factor and then set  $y_{i}=x_{i}$ and then integrate with respect to $x_{1},\ldots,x_{m}$ to get 
\begin{eqnarray}
&&E\(\(\int ( L^{ x+h}_{t}- L^{ x}_{ t})^{ 2}\,dx\)^{2}\)
\label{kacv2}\\
&&\qquad =4\sum_{\pi,a}  \int \int_{\{\sum_{i=1}^{4}r_{i}\leq t\}}\prod_{i=1}^{4}
\(\De^h_{x_{\pi(i)}}\)^{{a_{1}(i)}}\(\De^h_{x_{\pi(i-1)}}\)^{{a_{2}(i)}}\nn\\
&&\hspace{1 in}p^{\sharp}_{r_{i}}(x_{\pi(i)}-x_{\pi(i-1)})\,\,    \prod_{i=1}^{4}\,dr_{i}
\,\,    \prod_{i=1}^{2}\,dx_{i}, \nonumber
\end{eqnarray}
  as in (\ref{m.7}).
  As we did following (\ref{m.7}) we continue the analysis with  $p^{\sharp}$ replaced by $p$. 
  
In (\ref{kacv2}) the sum runs over all maps $\pi:\,[1,4]\mapsto [1,2]$ with $|\pi^{-1}(i)|=2$ for each $i$ and over all    
$a=(a_{ 1},a_{ 2})\,:\,[1,\ldots, 4]\mapsto \{ 0,1\}\times \{ 0,1\}$ with the
property that for each $i$ there   are   exactly two factors of the form $\De^{
h}_{ x_{i}}$. The factor 4 comes from the fact that we can interchange each $y_{i}$ and $x_{i}$, $i=1,2$.   As usual we take $\pi(0)=0$

Note  that in (\ref{kacv2}) it is possible to have  `bound states', that is values of $i$ for which $\pi (i)=\pi (i-1)$. 
 We first consider the terms in (\ref{kacv2}) with two bound states. There are  two possible maps. They are  $(\pi(1),\pi(2),\pi(3),\pi(4))=(1,1,2,2)$ and $(\pi(1), \pi(2),  \pi(3), 	\pi(4))=(2,2,1,\newline  1)$. The terms in (\ref{kacv2}) for the  map   $(\pi(1),\pi(2),\pi(3),\pi(4))=(1,1,2,\newline	2)$ are of the form
 
 \begin{equation}
   \prod_{i=1}^{4}
\(\De^h_{x_{\pi(i)}}\)^{{a_{1}(i)}}\(\De^h_{x_{\pi(i-1)}}\)^{{a_{2}(i)}} p _{r_{i}}(x_{\pi(i)}-x_{\pi(i-1)})\label{7.16},
   \end{equation}
 where the  density terms have the form
 \begin{equation}
   p_{r_{1}}(x_{1})p_{r_{2}}(y_{1}-x_{1} )p_{r_{3}}(x_{2}-y_{1} )p_{r_{4}}(y_{2}-x_{2} ),
   \end{equation}
 and where $y_{i}-x_{i}=0$. The value of the integrals of the terms in (\ref{7.16}) depend upon how the difference operators are distributed. In many cases the integrals are equal to zero. For example suppose we have
 \begin{equation}
  \De^{h}_{x_{1}}   p_{r_{1}}(x_{1})  \De^{h}_{x_{1}}p_{r_{2}}(0 )  \De^{h}_{x_{2}}p_{r_{3}}(x_{2}-x_{1} )\De^{h}_{x_{2}}p_{r_{4}}( 0),
   \end{equation}
   which we obtain by setting  $y_{1}=x_{1}$. (Note that $\De^{h}_{x_{1}}p_{r_{2}}(0 )  $ should be interpreted as  $\De^{h}_{x_{1}}p_{r_{2}}(x_{1}-y_{1} )  $ or $\De^{h}_{x_{1}}p_{r_{2}}(y_{1} -x_{1} ) $). Written out this term is 
 \bea
&&    \(p_{r_{1}}(x_{1}+h)-p_{r_{1}}(x_{1} ) \) \De^{h}_{x_{1}}p_{r_{2}}(0 ) \\
&&\qquad \(p_{r_{3}}(x_{2}-x_{1}+h )-p_{r_{3}}(x_{2}-x_{1} )\)\De^{h}_{x_{2}}p_{r_{4}}( 0)\nn
   \eea
  By a change of variables one sees that the integral of this term with respect to $x_{1}$ and $x_{2}$ is zero. 
  
  The only non-zero integrals  in  (\ref{7.16}) comes from 
  \begin{equation}
   p_{r_{1}}(x_{1})  \De^{h}\De^{-h}p_{r_{2}}(0 )   p_{r_{3}}(x_{2}-x_{1} )\De^{h}\De^{-h}  p_{r_{4}}( 0).
   \end{equation}
  (Similar to the above $\De^h \De^{-h} p_{r_{2}}( 0)$  is  $\De^h_{x_{1}}\De^{-h}_{y_{1}}p_{r_{2}}( x_{1}-y_{1})$ where $y_{1}=x_{1}$.)   The integral of this term with respect to $x_{1}$ and $x_{2}$ is
    \begin{equation}
   \De^{h}\De^{-h}p_{r_{2}}(0 )   \De^{h}\De^{-h}  p_{r_{4}}( 0).
   \end{equation}
 We get the same contribution when $(\pi(1),\pi(2),\pi(3),\pi(4))=(2,2,1,1)$.
 Consequently, the contribution to  (\ref{kacv2}) of maps with two bound states is 
\begin{eqnarray}
&&
  8\int_{\{\sum_{i=1}^{4}r_{i}\leq t\}}\De^h_{x}\De^{-h}_{x}p_{r_{2}}( 0)\, \De^h_{x}\De^{-h}_{x}p_{r_{4}}( 0)\,\prod_{i=1}^{4}\,dr_{i}\label{kacv.8} \\
  &&\qquad=32\int_{\{\sum_{i=1}^{4}r_{i}\leq t\}}\(p_{r_{2} }( 0)-p_{r_{2} }( h)\)\, \(p_{r_{4} }( 0)-p_{r_{4} }( h)\)\,\prod_{i=1}^{4}\,dr_{i}\nn\\
  &&\qquad=16\int_{\{u+v\leq t\}}(t-u-v)^{2}\(p_{u }( 0)-p_{u }( h)\)\, \((p_{v}( 0)-p_{v }( h)\)\,du\,dv.\nn\\
     &&\qquad\le 16t^{2}\(\int_{0}^{\ff} \(p_{u }( 0)-p_{u }( h)\)\, du\)^{2}=(4c_{\psi,h,0}t)^{2}  ,\nn
\end{eqnarray}
see (\ref{lim.7}). 


%


\medskip	
We next consider the contribution from terms with exactly one bound state. These   come from maps of the form $(\pi(1),\pi(2),\pi(3),\pi(4))=(1,2,2,1)$ or $(\pi(1),\pi(2),\pi(3),\pi(4))=(2,1,1,2)$. These terms  give   non-zero contributions of the form
\begin{eqnarray}
\lefteqn{Q_{2}:=\int\int_{\{\sum_{i=1}^{4}r_{i}\leq t\}}p_{r_{1}}(x) \De^{h}_{x}p_{r_{2}}(y-x)\, \De^h_{y}\De^{-h}_{y}p_{r_{3}}( 0)\,\De^h_{x} p_{r_{4}}(x-y)
\nn}\\
&&  \hspace{3 in} \,\prod_{i=1}^{4}\,dr_{i} \,dx\,dy\label{kacv.14}\\
&&\qquad 
=\int\int_{\{\sum_{i=1}^{4}r_{i}\leq t\}} \De^{-h}_{ }p_{r_{2}}(y )\, \De^h \De^{-h} p_{r_{3}}( 0)\,\De^{-h}_{y} p_{r_{4}}(y) \,\prod_{i=1}^{4}\,dr_{i} \,dy\nonumber;\\
\lefteqn{Q_{3}:=\int\int_{\{\sum_{i=1}^{4}r_{i}\leq t\}}p_{r_{1}}(x) \De^{h}_{x} \De^h_{y}p_{r_{2}}(y-x)\,  p_{r_{3}}( 0)\,\De^h_{x} \De^h_{y} p_{r_{4}}(x-y)
\nn}\\
&&  \hspace{3 in} \,\prod_{i=1}^{4}\,dr_{i} \,dx\,dy\label{pkacv.17aa}\\
&&
=\int\int_{\{\sum_{i=1}^{4}r_{i}\leq t\}} \De^h\De^{-h}p_{r_{2}}(y )\,   p_{r_{3}}( 0)\,\De^h\De^{-h} p_{r_{4}}(y) \,\prod_{i=1}^{4}\,dr_{i} \,dy;\nonumber
\end{eqnarray} 
and 
\begin{eqnarray}
\lefteqn{Q_{4}:=\int\int_{\{\sum_{i=1}^{4}r_{i}\leq t\}}p_{r_{1}}(x) \De^{h}_{x} \De^h_{y}p_{r_{2}}(y-x)\, \De^h_{y} p_{r_{3}}( 0)\,\De^h_{x}  p_{r_{4}}(x-y)
\nn}\\
&&  \hspace{3 in} \,\prod_{i=1}^{4}\,dr_{i} \,dx\,dy\label{pkacv.17bb}\\
&&
=\int\int_{\{\sum_{i=1}^{4}r_{i}\leq t\}} \De^h\De^{-h}p_{r_{2}}(y )\, \De^h  p_{r_{3}}( 0)\,\De^{-h} p_{r_{4}}(y) \,\prod_{i=1}^{4}\,dr_{i} \,dy.\nonumber
\end{eqnarray} 

For further explanation consider  $Q_{2}$. This arrangement comes from the sequence $(x_{1},y_{2},x_{2},y_{1} )$. The expression it is equal to comes by making the change of variables, $y-x\to y$ and then integrating with respect to $x$. 

Integrating and using    (\ref{bl.3})  we see that
  \bea
   Q_{2}&\le  &t \(\int_{0}^{1}|\De^{h}\De^{-h} p_{s }(0 )|\,\,ds\)\,\,\int \( \int_{0}^{t}\De^{-h} p_{r }(y )\,\,dr \)^{2}\,dy  
\label{pkacv.15}\\
&  \leq &{Ct \over h\psi (1/h)} \int \( \int_{0}^{t}\De^{-h}p_{r }(y )\,\,dr \)^{2}\,dy.\nonumber
      \eea
      Here we use  the fact that $\int \De^{-h} p_{r_{2} }(y )\De^{-h} p_{r_{4} }(y )\,dy\geq 0$ to extend the region of integration with respect to $r_{2}$ and $r_{4}$.
By Parseval's Theorem   and (\ref{el.1})
\begin{eqnarray} 
&&\int \( \int_{0}^{t}\De^{-h} p_{r }(y)\,\,dr \)^{2}\,dy
\label{pkacv.16}\\
&&\qquad ={1 \over 2\pi}  \int |1-e^{iph}|^{2}\( \int_{0}^{t}e^{-r\psi (p) }\,\,dr \)^{2}\,dp\nonumber\\
  &&\qquad \le {8 \over \pi} \int   \sin^{2} (hp/2)\(t\wedge {1 \over \psi  (p)}\)^{2}\,dp .
\nonumber 
\end{eqnarray}
Similar to the transition between (\ref{7.7})  and (\ref{7.8}) the last integral is bounded by
   \be C\(  h^{2}t^{2} \int_{0}^{\psi^{-1}(1/t)}  p^{2} \,dp +   h^{2}  \int_ {\psi^{-1}(1/t)} ^{1/h}{p^{2} \over \psi ^{2} (p)} \,dp +\int_{1/h}^\ff  {1 \over \psi^{2}  (p)} \,dp\) .
 \ee  
 Note that
    \begin{equation}
    h^{2}t^{2} \int_{0}^{\psi^{-1}(1/t)} p^{2}  \,dp\le Ch^{2}t^{2}\(\psi^{-1}(1/t)\)^{3}.
   \end{equation}
 This bound is the right hand side of (\ref{7.9}).
 Bounds for the other integrals are given in (\ref{7.10})--(\ref{7.13}). Since the bounds in  (\ref{7.9})--(\ref{7.13}) give (\ref{4.9}), we see that
 \be 
   Q_{2}\le  {Ctg(h,t) \over h\psi (1/h)}.
         \ee

%

%

%
%


 To obtain a bound for  $Q_{3}$ we use  (\ref{bl.5}) and (\ref{pkacv.19}) to see that  it is bounded in absolute value by 
\begin{equation}
t \(\int_{0}^{t} p_{s }(0 )\,\,ds\)\,\,\int \( \int_{0}^{1}|\De^{h} \De^{-h} p_{r }(y )|\,\,dr \)^{2}\,dy  \label{pkacv.17a}\leq C{ t^{2}\psi^{-1}(1/t)  \over h\psi^{2}(1/h)}.
\end{equation}
Integrating $Q_{4}$  and using the Cauchy-Schwarz Inequality we see that  it is bounded in absolute value by  
 \begin{equation}
  t\bigg|\int_{0}^{1} \De^h  p_{r }( 0)\,dr \bigg|\(\int\left|\int_{0 }^1  \De^h\De^{-h}p_{r }(y )\, dr\right|^{2} \,dy    \int\left|\int_{0 }^1   \De^{-h}p_{r }(y )\, dr\right|^{2} \,dy \)^{1/2}.
   \end{equation}
By (\ref{bl.2}),  
 (\ref{bl.5}) and   (\ref{lastf}) we get  
\begin{equation}
   Q_{4}\le \frac{Ct }{h^{3/2}\psi^{5/2}(1/h)}.
   \end{equation}
  
  \medskip

Finally, we consider  the contribution from terms in (\ref{kacv2}) with no bound states. These have to be from $\pi$ of the form $(\pi(1),\pi(2),\pi(3),\pi(4))=(1,2,1,2)$ or of the form $(\pi(1),\pi(2),\pi(3),\pi(4))=(2,1,2,1)$. 
They give contributions of the form
\begin{eqnarray}
\lefteqn{Q_{5}:=\int\int_{\{\sum_{i=1}^{4}r_{i}\leq t\}}p_{r_{1}}(x) \De^{h}_{x}p_{r_{2}}(y-x)\, \De^h_{y}\De^{h}_{x}p_{r_{3}}( x-y)\,\De^h_{y} p_{r_{4}}(y-x)\nn}
\\
&&  \hspace{3 in} \,\prod_{i=1}^{4}\,dr_{i} \,dx\,dy\label{pkacv.20}\\
&&
=\int\int_{\{\sum_{i=1}^{4}r_{i}\leq t\}} \De^{-h}p_{r_{2}}(y )\, \De^h \De^{-h}p_{r_{3}}( y)\,\De^{h} p_{r_{4}}(y) \,\prod_{i=1}^{4}\,dr_{i} \,dy\nonumber
\end{eqnarray} 
and  
\begin{eqnarray}
\lefteqn{Q_{6}:=\int\int_{\{\sum_{i=1}^{4}r_{i}\leq t\}}p_{r_{1}}(x) \De^{h}_{x} \De^h_{y}p_{r_{2}}(y-x)\,  p_{r_{3}}( x-y)\,\De^h_{x} \De^h_{y} p_{r_{4}}(x-y)
\nn}\\
&&  \hspace{3 in} \,\prod_{i=1}^{4}\,dr_{i} \,dx\,dy\label{pkacv.21}\\
&&
=\int\int_{\{\sum_{i=1}^{4}r_{i}\leq t\}} \De^h\De^{-h}p_{r_{2}}(y )\,   p_{r_{3}}( y)\,\De^h\De^{-h} p_{r_{4}}(y) \,\prod_{i=1}^{4}\,dr_{i} \,dy\nonumber.
\end{eqnarray} 
Clearly   
\bea
&&  Q_{5}   
 \le t \int   \(\int_{0}^{1}|  \De^{-h}p_{r }(y ) |\,dr \)\\
&&\qquad\qquad\qquad  \(\int_{0}^{1}|  \De^{h}p_{r }(y ) |\,dr \) \(\int_{0}^{1}|  \De^{h}\De^{-h}p_{r }(y ) |\,dr \)\,dy.\nn
   \eea   
Using    (\ref{bl.2}),  and (\ref{bl.4}) we see that 
\begin{equation}
   Q_{5}\le \frac{Ct \log 1/h}{h^{2}\psi^{3}(1/h)}.
   \end{equation}
 The term  $Q_{6}$ is bounded the same way we bounded $Q_{3}$ and has the same bound.

\medskip	 
It follows from (\ref{expep}), Lemma \ref{lem-exp} and (\ref{kacv.8}) that 
 \begin{equation}
   \mbox{Var Z}\le C\( \sum_{j=2}^{6}|Q_{j}|+\({tg(h,t)\over h\psi(1/h)} \)\)
   \end{equation}
as $h\to 0$,  since $g(h,t)<t/(h\psi(1/h)$.
(We need a large constant because expressions for   $Q_{j}$, $j=2,\ldots,6$ occur may ways, according to combinatorics of the distribution of the difference operators.)   

  We leave it to the reader to verify that replacing $p$ by $p^{\sharp}$ only adds error terms that do not change (\ref{4.11}) and (\ref{4.11a}).
\qed

\noindent{\bf  Proof of Lemma \ref{lem-4.5} } Use (\ref{pkacv.4})--(\ref{7.8}) with $\psi^{-1}(1/t)$ replaced by 1. In place of (\ref{7.9}) we have 
   \begin{equation}
    h^{2} \int_{0}^{1} {p^{2}\over \psi (p)} \,dp\le Ch^{2}. \label{7.9a}
   \end{equation}
 (Since $\lim_{p\to 0}\psi(p)/p^{2}>0$    this integral is finite; (see \cite[Lemma 4.2.2]{book}).     
   
In place of (\ref{7.11}) we have,
 if $\bb>3/2$
\begin{equation}
    h^{2}  \int_ {1} ^{1/h}{p^{2} \over \psi^{2} (p)} \,dp \le Ch^{2}.   \end{equation}
 The statements in (\ref{7.12}) and (\ref{7.13}) remain the same when  $\psi^{-1}(1/t)$ replaced by 1. With these changes the proof of (\ref{4.9}) gives (\ref{4.11}).\qed

\noindent{\bf Proof of Theorem \ref{theo-clt2s} }
 \label{rem-6.2}{\rm We note that by Remark \ref{rem-6.1} and Lemma \ref{lem-4.5}, for symmetric stable processes of index $\bb$,   \begin{equation}
E\( \int (  L^{ x+h}_{1}- L^{ x}_{ 1})^{ 2}\,dx\)=   4 c_{\bb,0}h^{\bb-1}  +O\(\ov g(h )\).\label{asymexp}
\ee
Since  
\begin{equation}
   \lim_{h\to 0}\sqrt{h\psi^{2}(1/h)}\,\ov g(h )=0
   \end{equation}
 we get (\ref{5.0weaks}) with $t=1$.  As we remark in the paragraph proceeding (\ref{1.10}) this suffices to prove Theorem \ref{theo-clt2s} as stated.
  \qed

\section{Appendix: Kac Moment Formula}\label{sec-Kac}

Let    $X=\{X_{t},t\in R_{+} \}$    denote a symmetric L\'evy process with continuos local time $L=\{L_{t}^{x}\,;\,(x,t)\in R^{1}\times R_{+}\}$.   Since $L$ is continuous   we have the occupation density formula,  
 \begin{equation}
\int^{t}_{0} g(X_{s})\,ds=\int g(x) L_{t}^{x}\,dx\label{odf},
 \end{equation}
 for all continuous functions  $g $ with compact support. (See, e.g.  \cite[Theorem 3.7.1]{book}.)

 Let $f(x)$ be a  continuous  function on $R^{1}$ with   compact  support with $\int f(x)\,dx=1$. Let $f_{\ep,y}(x):={1 \over \ep}f\({x-y \over \ep}\)$.   I.e., $f_{\ep,y} (x)$ is an approximate $\de$-function at $x$. Set
\begin{equation}
L^{ x }_{ t,\ep}=\int_{0}^{t} f_{\ep,x }\(X_{s}\)\,ds. \label{ka.3}
\end{equation}
It follows from  (\ref{odf}) that
\begin{equation}
L^{ x }_{ t }=\lim_{\ep\rar 0}L^{ x }_{ t,\ep}\hspace{.2 in}\mbox{a. s. }\label{ka.0}
\end{equation}
Let $p_{t}(x,y)$ denote the probability density of $X_{t}$.

\begin{theorem} [Kac Moment Formula]\label{KMT}
For any fixed $0<t<\ff$, bounded continuous $g$, and any $x_{1},\ldots, x_{m},z\in R^{1}$,
\begin{eqnarray}
&& 
E^{z}\(\prod_{i=1}^{m}      L^{ x_{i}}_{ t} \,g(X_{t})\)= \sum_{\pi} \int_{\{\sum_{j=1}^{m} r_{j}\leq t\} } \prod_{j=1}^{m}p_{r_{j}}(x_{\pi(j-1)},x_{\pi(j)})\label{ka.1}\\
&& \hspace{2 in}\(\int p_{t-r_{m}}(x_{\pi(m)},y)g(y)\,dy\)\prod_{j=1}^{m}\,dr_{j},
\nn
\end{eqnarray}
where the sums run over all permutations  $\pi$ of $\{1,\ldots, m\}$ and    $\pi(0) : =0$
and $x_{0}: =z$.

 \end{theorem}
 
\Proof Let
\bea
&&
F_{t}(x_{1},\ldots, x_{m})=\int_{\{\sum_{j=1}^{m} r_{j}\leq t\} } \prod_{j=1}^{m}p_{r_{j}}(x_{ j-1 },x_{ j}) \label{ka.10}\\
&& \hspace{2 in}\(\int p_{t-r_{m}}(x_{m},y)g(y)\,dy\)\prod_{j=1}^{m}\,dr_{j}\nn
\eea
 Then 
\begin{eqnarray}
\lefteqn{ 
E^{z}\(\prod_{i=1}^{m}      L^{ x_{i}}_{ t,\ep}\,g(X_{t}) \)\label{ka.4j} }\\
&&= \sum_{\pi}\int_{\{0\leq t_{\pi(1)}\leq \ldots \leq t_{\pi(m)}\leq t\} }E^{z}\( \prod_{j=1}^{m}f_{\ep,x_{j} }\(X_{t_{\pi(j)}}\)\,g(X_{t}) \)\prod_{j=1}^{m}\,dt_{\pi(j)}
\nn\\
&&= \sum_{\pi}\int_{\{0\leq t_{1}\leq \ldots \leq t_{m}\leq t\} }E^{z}\( \prod_{j=1}^{m}f_{\ep,x_{\pi(j)} }\(X_{t_{j}}\) \,g(X_{t})\)\prod_{j=1}^{m}\,dt_{j}
\nn\\
&&= \sum_{\pi}\int \int_{\{\sum_{j=1}^{m} r_{j}\leq t\} } \prod_{j=1}^{m}f_{\ep,x_{\pi(j)} }(y_{j})p_{r_{j}}(y_{j-1},y_{j})\\
&& \hspace{1.5 in}\(\int p_{t-r_{m}}(y_{m},y)g(y)\,dy\)\prod_{j=1}^{m}\,dr_{j}\nn\,dy_{j}
\nn\\
&&= \sum_{\pi}\int F_{t}( y_{0},\ldots,  y_{m})\prod_{j=1}^{m}f_{\ep,x_{\pi(j)} }(y_{j})\,dy_{j}
\nn
\end{eqnarray}
where $ y_{0} :=z$.

 Since the integrand in (\ref{ka.10}) is dominated by $(2\pi)^{-m/2}\prod_{j=1}^{m }r_{j}^{-1/2}$
it follows from the Dominated Convergence Theorem  that $F_{t}(x_{1},\ldots, x_{m})$
is a continuous function of  $(x_{1},\ldots, x_{m}) $ for all $0\leq t<\ff$ and all $m$. 
It then follows immediately from (\ref{ka.4j}) and the fact that  $\prod_{j=1}^{m}f_{\ep,x_{\pi(j)} }(y_{j})$ has compact support that
\begin{eqnarray}
&& 
\lim_{\ep\rar 0}E\(\prod_{i=1}^{m}      L^{ x_{i}}_{ t,\ep}\,g(X_{t}) \) = \sum_{\pi}F_{t}( x_{\pi(0)},x_{\pi(1)},\ldots,  x_{\pi(m)}).
\label{ka.11}
\end{eqnarray}
 A repetition of the above proof  shows that $E\(\lc \prod_{i=1}^{m}      L^{ x_{i}}_{ t,\ep}\rc^{2} \)$
is bounded uniformly in $\ep>0$.  This fact and  (\ref{ka.0})   show that 
\begin{equation}
\lim_{\ep\rar 0}E\(\prod_{i=1}^{m}      L^{ x_{i}}_{ t,\ep}\,g(X_{t}) \)=E\(\prod_{i=1}^{m}      L^{ x_{i}}_{ t }\,g(X_{t}) \).\label{ka.11j}
\end{equation}
Obviously (\ref{ka.11}) and (\ref{ka.11j}) imply (\ref{ka.1}). \qed

  \def\noopsort#1{} \def\printfirst#1#2{#1}
\def\singleletter#1{#1}
      \def\switchargs#1#2{#2#1}
\def\bibsameauth{\leavevmode\vrule height .1ex
      depth 0pt width 2.3em\relax\,}
\makeatletter
\renewcommand{\@biblabel}[1]{\hfill#1.}\makeatother
\newcommand{\bysame}{\leavevmode\hbox to3em{\hrulefill}\,}

\bigskip
\noindent
\begin{tabular}{lll} 
      & Jay Rosen & Michael Marcus\\
      & Department of Mathematics& Department of Mathematics\\
     &College of Staten Island, CUNY& City College, CUNY\\
     &Staten Island, NY 10314& New York, NY 10031\\ &jrosen30@optimum.net &mbmarcus@optonline.net
\end{tabular}

\end{document}